\documentclass{article}
\usepackage[utf8]{inputenc}
\usepackage[a4paper, total={15.92cm, 24.37cm}, top=1in, left=1in]{geometry}
\usepackage[T1]{fontenc}
\usepackage{lmodern}
\usepackage{amsmath} \allowdisplaybreaks
\usepackage{dsfont} 
\usepackage{bm}
\usepackage{amsthm}
\newtheorem{thm}{Theorem}
\newtheorem{prp}{Proposition}[section]
\newtheorem{dfn}[prp]{Definition}
\newtheorem{lmm}[prp]{Lemma}
\newtheorem{cor}[prp]{Corollary}

\title{%
    Asymptotic exponential law for the transition time to equilibrium\\
    of the metastable kinetic Ising model\\
    with vanishing magnetic field
}
\author{%
    A. Gaudillière\thanks{%
        Aix Marseille Univ, CNRS, Centrale Marseille, I2M,
        Marseille,
        France,
        e-mail: alexandre.gaudilliere@math.cnrs.fr
    }
    \and
    P. Milanesi\thanks{%
        Aix Marseille Univ, CNRS, Centrale Marseille, I2M,
        Marseille,
        France,
        e-mail: paolo.milanesi@univ-amu.fr 
    }
    \and 
    M. E. Vares\thanks{%
        \hbox{%
            Instituto de Matem\'atica, UFRJ,
            Av. Athos da Silveira Ramos 149, Cidade Universit\'aria,}
        \hbox{\indent\indent
            Ilha do Fund\'ao
            21941-909~Rio de Janeiro, RJ,
            Brasil,
            e-mail: eulalia@in.ufrj.br
        }
    }
}
\date{November 18, 2019}

\begin{document}
\maketitle
\begin{abstract}
    We consider a Glauber dynamics 
    associated with the Ising model 
    on a large two-dimensional box 
    with minus boundary conditions
    and in the limit of a vanishing positive external magnetic field.
    The volume of this box increases quadratically
    in the inverse of the magnetic field.
    We show that at subcritical temperature
    and for a large class of starting measures,
    including measures that are supported by configurations
    with macroscopic plus-spin droplets,
    the system rapidly relaxes to some metastable equilibrium
    ---with typical configurations 
    made of microscopic plus-phase droplets
    in a sea of minus spins---
    before making a transition 
    at an asymptotically exponential random time
    towards equilibrium 
    ---with typical configurations made of 
    microscopic minus-phase droplets
    in a sea of plus spins
    inside a large contour that separates
    this plus phase from the boundary.
    We get this result
    by bounding from above the local relaxation times
    towards metastable and stable equilibria.
    This makes possible to give a pathwise description
    of such a transition,
    to control the asymptotic behaviour
    of the mixing time
    in terms of soft capacities
    and to give estimates of these capacities.
    \par\bigskip\noindent
    {\it MSC 2010:} primary: 82C20; secondary: 60J27, 60J45, 60J75.
    \par\smallskip\noindent
    {\it Keywords:} Metastability, Glauber dynamics, exponential law,
        relaxation time, quasi-stationary measures, potential theory.
    \par\smallskip\noindent
    {\it Acknowledgments:} 
    M. E. V. thanks Roberto Schonmann for many long conversations
    on metastability, in particular for the hospitality at UCLA
    back in 1997, when she was studying the paper~\cite{SS},
    and they discussed the difficulties to achieve a result
    along the line of the current paper;
    she also thanks Augusto Q. Teixeira for discussions
    on the subject matter of the paper,
    and Vladas Sidoravicius (in memoriam)
    for inspiring general discussions on metastability.
    A. G. and P. M. thank Julien Sohier
    for the fruitful discussions they had in Leiden 
    when studying \cite{SS},
    on which much of this work is based.
    This was possible thanks to the kind hospitality 
    of Leiden university, 
    which hosted them for two fall seasons
    through Frank den Hollander's ERC Advanced Grant 267356-VARIS. 
    A.~G. and P.~M. also thank the kind hospitality of the Universidade Federal
    do Rio de Janeiro through Faperj E26/102.338/2013.
    M.~E.~V. acknowledges partial support of CNPq (grant 305075/2016-0)
    and Faperj  E-26/203.948/2016.
\end{abstract}

\section{Model and results}

\subsection{Glauber dynamics for the Ising model}
For a finite subset $\Lambda$ of $\mathds{Z}^2$
and $\eta \in \Omega_{\mathds{Z}^2} = \{-1, +1\}^{\mathds{Z}^2}$,
the Ising model in the domain $\Lambda$,
with boundary conditions $\eta$,
at inverse temperature $\beta > 0$
and with magnetic field $h \in \mathds{R}$,
is associated with the Hamiltonian
\begin{equation}\label{clo}
    H_{\Lambda, \eta, h}(\sigma)
    = -{1 \over 2} \sum_{
        \scriptstyle \{x, y\} \subset \Lambda,
        \atop \scriptstyle \|x - y\|_1 = 1
    } \sigma(x) \sigma(y) 
    - {1 \over 2} \sum_{
        \scriptstyle x \in \Lambda,\, y \not\in \Lambda,
        \atop \scriptstyle \|x - y\|_1 = 1
    } \sigma(x) \eta(y) 
    - {h \over 2} \sum_{x \in \Lambda} \sigma(x),
    \qquad \sigma \in \Omega_{\Lambda} = \{-1, +1\}^{\Lambda},
\end{equation}
the partition function
$$
    Z_{\Lambda, \eta, h}
    = \sum_{\sigma \in \Omega_{\Lambda}} e^{-\beta H_{\Lambda, \eta, h}(\sigma)}
$$
and the Gibbs measure
$$
    \mu_{\Lambda, \eta, h}(\sigma) = {
        e^{-\beta H_{\Lambda, \eta, h}(\sigma)}
    \over
        Z_{\Lambda, \eta, h}
    }\,,
    \qquad \sigma \in \Omega_{\Lambda}.
$$
The maybe unusual factors $1 / 2$ in Equation~\eqref{clo}
are here to stick to the conventions of \cite{SS},
which is the main reference we will follow.

The associated Glauber dynamics
are irreducible continuous time Markov processes
$$
    X_{\Lambda, \eta, h} = \bigl(X_{\Lambda, \eta, h}(t)\bigr)_{t \geq 0}
$$
with a single spin flip generator
$$
    ({\cal L}_{\Lambda, \eta, h} f)(\sigma)
    = \sum_{x \in \Lambda} w(\sigma, \sigma^x) \bigl[
        f(\sigma^x) - f(\sigma)
    \bigr],
    \qquad f: \Omega_\Lambda \rightarrow \mathds{R},
    \qquad \sigma \in \Omega_\Lambda,
$$
where the configuration $\sigma^x$
is obtained from $\sigma$ by flipping the spin at $x$,
$$
    \sigma^x(y) = \left\{
        \begin{array}{ll}
            \sigma(y) & \hbox{if $x \neq y$,}\\
            -\sigma(x) & \hbox{if $x = y$,}
        \end{array}
    \right.
$$
and the transition rates $w(\sigma, \sigma^x)$
are chosen to satisfy the detailed balance equations
$$
    \mu_{\Lambda, \eta, h}(\sigma) w(\sigma, \sigma^x)
    = \mu_{\Lambda, \eta, h}(\sigma^x) w(\sigma^x, \sigma),
    \qquad \sigma \in \Omega_\Lambda,
    \qquad x \in \Lambda.
$$
One can for example consider a Metropolis dynamics 
with
$$
    w(\sigma, \sigma^x) = \exp\left\{
        -\beta \bigl[
            H_{\Lambda, \eta, h}(\sigma^x) - H_{\Lambda, \eta, h}(\sigma)
        \bigr]_+
    \right\}, 
    \qquad \sigma \in \Omega_{\Lambda},
    \qquad x \in \Lambda,
$$
where the brackets $[\cdot]_+$ stand for the positive part,
or a heat bath dynamics
$$
    w(\sigma, \sigma^x)
    = {
        \exp\bigl\{-\beta H_{\Lambda, \eta, h}(\sigma^x)\bigr\}
        \over \exp\bigl\{-\beta H_{\Lambda, \eta, h}(\sigma)\bigr\} + \exp\bigl\{-\beta H_{\Lambda, \eta, h}(\sigma^x)\bigr\}
    }\,, 
    \qquad \sigma \in \Omega_\Lambda,
    \qquad x \in \Lambda.
$$

In this paper we will consider such a dynamics $X_{\Lambda_h, -, h}$
in the limit of a vanishing positive magnetic field $h \ll 1$,
with uniform minus boundary conditions
and inside a box $\Lambda_h$,
the volume\footnote{
    Working in dimension two,
    the word ``area'' could have been more appropriate.
    We will follow the usage by referring to volumes
    and surfaces rather than areas and perimeters.
}
of which will quadratically diverge in $1 / h$. 
As far as the jump rates $w(\sigma, \sigma^x)$ are concerned,
we will only assume
that there are two positive constants
$w_{\min}$ and $w_{\max}$,
possibly depending on our fixed parameter $\beta$,
such that 
$$
    w_{\min}
    \leq w(\sigma, \sigma^x) 
    \leq w_{\max},
    \qquad \sigma \in \Omega_{\Lambda_h},
    \qquad x \in \Lambda_h,
$$
which implies in particular
that $X_{\Lambda_h, -, h}$ is irreducible.

\subsection{Metastability issues}
This kind of evolution
is used as a dynamic model
to study hysteresis phenomena.
The critical temperature of a ferromagnet
is the temperature below which,
when exposed to a strong negative external magnetic field,
it keeps a spontaneous negative magnetization
after removing this external field.
Then, by exposing the ferromagnet
to a small enough positive magnetic field
it will keep a higher, but still negative,
magnetization for a long time,
typically longer than usual experiment times.
One gets a positive magnetization
only by increasing the value of the external field,
or waiting long enough for a relaxation to equilibrium.
Then, by removing again the magnetic field 
before making it decrease back to negative values,
the same kind of picture reappears:
the ferromagnet gets a spontaneous positive magnetization,
then a smaller but still positive magnetization
before jumping to an equilibrium negative magnetization
after a long enough time
or after reaching low enough values
for the external field.
Two of the main questions associated
with such a phenomenon are those of 
i) describing such a metastable equilibrium
and in particular such a higher,
but still negative, magnetization;
ii) characterizing such a late and abrupt relaxation
to equilibrium, and in particular computing
the order of magnitude of this relaxation time.

In the fundamental paper \cite{SS},
Schonmann and Shlosman studied such a dynamics $X_\infty$ in infinite volume
and they described the state of the system 
at time $t = e^{\alpha / h}$ 
for positive $\alpha$,
with vanishing magnetic field $0 < h \ll 1$,
at any subcritical temperature $1 / \beta < 1 / \beta_c$
when starting from any initial measure $\nu$
stochastically dominated by $\mu_-$, which is the thermodynamic limit
of the Ising model in a finite box with minus boundary conditions
and zero magnetic field.
They identified a critical $\alpha_c$
such that for any $\alpha < \alpha_c$
the mean value $\mathds{E}_\nu[f(X_\infty(t))]$ of any local observable 
$f : \{-1, +1\}^{\mathds{Z}^2} \!\rightarrow \mathds{R}$
is close to the ${\cal C}^k$ continuations
of its expected values for negative values of the magnetic field
$h < 0 \rightarrow \mu_h(f)$,
with $\mu_h$ the thermodynamic limit of the Ising model
in a finite box with non-zero magnetic field $h$.
More precisely they answered the first question by proving that,
for all $k > 0$,
\begin{equation}\label{alain}
    \mathds{E}_\nu\bigl[
        f(X_\infty(t))
    \bigr] = \sum_{j < k} {h^j \over j!} \left.
        d^j \mu_h(f) \over dh^j
    \right|_{h = 0_-} + O\bigl(h^k\bigr).
\end{equation}
As far as the second question is concerned
they also proved that for any $\alpha > \alpha_c$
the mean value $\mathds{E}_\nu[f(X_\infty(t))]$
of any local observable $f$
is close to its expected value $\mu_h(f)$.
The formula they established for $\alpha_c$
is particularly remarkable:
\begin{equation}\label{noemie}
    \alpha_c = {\beta w_\beta^2 \over 12 m^*_\beta}\,,
\end{equation}
where $m^*_\beta$ is the spontaneous magnetization
at inverse temperature $\beta$,
$$
    m^*_\beta = -\mu_-(\sigma_0)
$$
with $\sigma_0$ the local observable defined by 
$\sigma_0 : \omega \in \Omega_{\mathds{Z}^2} \mapsto \omega(0)$,
and $w_\beta$ is the surface tension
of the unitary volume Wulff shape
(see Section~\ref{aldo}).

At this point it remains to describe the evolution
of the system at times of order $e^{\alpha_c / h}$,
the order of the relaxation time of this dynamics.
Since we are in the regime $h \ll 1$,
for any given $\alpha \neq \alpha_c$
the two cases $\alpha < \alpha_c$
and $\alpha > \alpha_c$
refer to very small and very large times $t = e^{\alpha / h}$
with respect to $e^{\alpha_c / h}$.
The $O(h^k)$ in formula \eqref{alain}
depends on $\alpha < \alpha_c$ just as,
in the case $\alpha > \alpha_c$, the ``small enough $h$''
from which $\mathds{E}_\nu[f(X_\infty(e^{\alpha / h}))]$
will be close to $\mu_h(f)$ depends on $\alpha$.
More precisely it holds,
for any given $\epsilon >0$,
$$
    \left|
        \mathds{E}_\nu[f(X_\infty(e^{\alpha /h}))]
        - \mu_h(h)
    \right| < \epsilon
$$
for $h < h_0(\alpha)$;
and $h_0(\alpha)$ vanishes as $h$ does.
One cannot then use these results to describe
the system at times $t$
of order $e^{\alpha_c / h}$
for small $h > 0$.
This is the goal of this paper 
in the simpler case of the dynamics $X_{\Lambda_h, -, h}$ 
on, instead of the infinite volume $\mathds{Z}^2$,
a Wulff shape domain $\Lambda_h$
containing around $(B_{\max} / h)^2$ sites
for a large enough $B_{\max} > 0$.
The box $\Lambda_h$ is formally defined by 
$$
    \Lambda_h = \left(
        {B_{\max} \over h} W
    \right) \cap \mathds{Z}^2
$$
with $W$ defined after Equation~\eqref{benjamin}
at page~\pageref{benjamin}. 
As it will be clear
from the heuristics of the next section,
that goes back to Schonmann and Shlosman indeed, 
with a small $B_{\max}$ we would not have
any metastable behaviour:
equilibrium would look like the minus phase.
On the contrary, with a large $B_{\max}$,
and with such a box shape,
the plus phase will invade the whole box at equilibrium,
due to the positivity of the magnetic field
and despite the minus boundary conditions.

\subsection{A pathwise description}
In this finite volume case,
we can give another description,
in terms of {\it restricted ensemble},
of the metastable equilibrium
by following \cite{SS}.
The configurations in $\Omega_{\Lambda_h}$,
which we identify with 
$$
    \Omega_{\Lambda_h, -}
    = \left\{
        \sigma \in \Omega_{\mathds{Z}^2} 
        : \sigma(x) = -1 \hbox{ for all $x \not\in \Lambda_h$}
    \right\},
$$
can be described as a collection 
of closed self-avoiding contours on the dual lattice,
which separate plus spins from minus spins.
In doing so we adopt a standard ``splitting rule'',
the one used in \cite{DKS} (Section~3.1 there).
We call {\it external contour} 
of a given configuration
any contour that is not surrounded
by any other contour.
We define ${\cal R}_-$
as the set of configurations in $\Omega_{\Lambda_h}$
such that the volume of each external contour,
i.e., the number of sites enclosed in it,
is smaller than $(B_c / h)^2$
with 
\begin{equation}\label{pierre}
    B_c = {w_\beta \over 2 m^*_\beta}
    \,.
\end{equation}
The expansion \eqref{alain} is actually
an expansion for $\mu_{\Lambda_h, -, h}(f| {\cal R}_-)$.
Our pathwise description will also make
use of such a restricted ensemble
$\mu_{\Lambda_h, -, h}(\cdot\,|{\cal R})$
but for another ${\cal R} \neq {\cal R}_-$.
The reader can think of ${\cal R}$
as a set that is smaller than ${\cal R}_-$,
since some configurations with limited volume
but large perimeter are allowed in the latter
and will be excluded from the former.
However ${\cal R}$ will not be a subset of ${\cal R}_-$,
since it will include slightly supercritical configurations
in the sense of the heuristics of the next paragraph,
while all configurations in ${\cal R}_-$
are subcritical.

Before describing the set ${\cal R}$ we will choose,
let us first recall the heuristics
where Formula \eqref{pierre}
comes from.
If $w_\beta$ is the surface free energy
of a unitary volume Wulff shape $W$,
then the free energy of a discrete ``plus phase'' Wulff shape
with a volume of order $(B / h)^2$ 
in a ``minus phase''
can be estimated, for $h \ll 1$
and up to an additive function
that does not depends on $B$, by
$$
    \Phi\left(
        {B \over h} W
    \right)
    = w_\beta {B \over h} - 2 {h \over 2} \left(
        B \over h
    \right)^2 m^*_\beta
    = {1 \over h} \left[
        w_\beta B - m^*_\beta B^2
    \right].
$$
We will refer to the quantity
$B / h$ as the {\it linear size\/}
of such a Wulff shape with volume
$(B / h)^2$.
The $1 / 2$ factor in the previous equation
comes from the Hamiltonian,
while the factor $2$ accounts for the volume
of the plus phase as well as the volume
of the minus phase,
which is the volume of $\Lambda_h$
minus the volume of the Wulff droplet.
Let us set
\begin{equation}\label{irene}
    \phi(B)
    = \left[
        w_\beta B - m^*_\beta B^2
    \right]
    = {w_\beta^2 \over 4 m^*_\beta}- m^*_\beta \Biggl(
        B - {w_\beta \over 2 m^*_\beta}
    \Biggl)^2
    = A - m^*_\beta \bigl(
        B - B_c 
    \bigr)^2
\end{equation}
with
\begin{equation}\label{charlotte}
    A = {w_\beta^2 \over 4 m^*_\beta}\,.
\end{equation}
This computation suggests
that a plus phase Wulff droplet 
of size $(B / h)^2$ 
will have a tendency to shrink or grow
depending on $B < B_c$ or $B > B_c$.
Being the Wulff shape a minimizer
of the surface free energy 
for a given volume,
critical Wulff droplets 
of size $B_c / h$ will indeed constitute
a bottleneck for the dynamics
and we will refer to the cases $B < B_c$
and $B > B_c$ as the {\it subcritical}
and {\it supercritical cases}.

To make rigorous such free energy estimates,
we will follow \cite{SS} 
and use the skeleton description of contours
of \cite{DKS}.
Skeletons are associated with
long enough contours only.
This motivates the following definition
inherited from \cite{SS} and extended
to all contours, external or not.
\begin{dfn}
    Let $b$ a positive number which is less than $1/4$.
    A contour is said {\em $b$-vertebrate},
    or simply {\em vertebrate},
    if it encloses more than $1 / h^{2b}$
    sites in its interior.
    A contour is said {\em $b$-invertebrate},
    or simply {\em invertebrate}
    if the number of sites that are enclosed 
    in its interior is less than
    or equal to $1 / h^{2b}$.
\end{dfn}

We are now ready to define our set ${\cal R}$.
To this end we introduce another parameter $B_+ > B_c$,
which has to be thought of
as close\footnote{
    As long as $\phi(B_+)$ is positive
    the restricted ensemble $\mu_{\Lambda_h, -, h}(\cdot | {\cal R})$
    will be concentrated on the same kind of configurations,
    but, because some dynamical quantities
    will also play a role,
    we will get stronger results 
    by taking $B_+$ close to $B_c$
    rather than only asking for the positivity of $\phi(B_+)$.
} to $B_c$, and which, just as $b$, will not depend on $h$.
\begin{dfn}
    For $0 < b < 1 / 4$ and $B_+ > B_c$,
    we call ${\cal R}$ the set 
    of all configurations $\sigma$ in $\Omega_{\Lambda_h}$
    for which one can find a collection
    of at most $1 / h^{(1 - b / 2)}$ disjoint Wulff shapes
    and with total linear size less than $B_+ / h$
    that contains all the $b$-vertebrate contours of $\sigma$.
\end{dfn}
The reader can think of the relevant configurations
in ${\cal R}$ as those with only one large contour
enclosed in a subcritical,
or slightly supercritical, Wulff shaped box.
The reason why we need an upper bound
on the number of involved boxes is technical.
At some point
(see inequality \eqref{paul} at page \pageref{paul})
we will need to upper bound
the number of such possible box arrangements,
and this restriction will help.

We define the mixing time of $X_{\Lambda_h, -, h}$ by
$$
    t_{{\rm mix},h} = \inf\left\{
        t \geq 0 :
        \forall \sigma \in \Omega_{\Lambda_h},\,
        \forall E \subset \Omega_{\Lambda_h},\,
        \Bigl|
            \mathds{P}_\sigma\bigl(
                X_{\Lambda_h, -, h}(t) \in E
            \bigr) - \mu_{\Lambda_h, -, h}(E)
        \Bigr| \leq {1 \over e}
    \right\},
$$
with $\mathds{P}_\sigma$ the probability measure
associated with $X_{\Lambda_h, -, h}$ started in $\sigma$;
so that the total variation distance
between $\mu_{\Lambda_h, -, h}$
and the law of $X_{\Lambda_h, -, h}(t)$
is exponentially small in $t$ 
for $t$ larger than $t_{{\rm mix}, h}$.
By using techniques from \cite{SS}
one could get the following proposition,
that we will obtain as a byproduct
of our main results.
\begin{prp} \label{amandine}
    For all supercritical $\beta > \beta_c$
    and any $B_{\max} > 2 B_c$ 
    it holds (recall \eqref{charlotte})
    \begin{equation}\label{gregory}
        \lim_{h \rightarrow 0} h \ln (t_{{\rm mix}, h})
        = \beta A.
    \end{equation}
\end{prp}
To describe our dynamics on this time scale $t_{{\rm mix}, h}$
we will use a suitable random time $T$
so that, starting from the restricted ensemble
$\mu_{\Lambda_h, -, h}(\cdot\,|{\cal R})$,
the rescaled time $T / t_{{\rm mix}, h}$
will converge in law to an exponential random variable
of mean one
and, for $t > T$,
the law of $X_{\Lambda_h, -, h}(t)$ will be close to $\mu_{\Lambda_h, -, h}$.
The definition of $T = T_{\lambda_{\cal S}}$
involves another set of configurations ${\cal S}$
(see Definition~\ref{stephane})
and a further randomization: 
it can be interpreted as a killing time under a killing rate $\lambda_{\cal S}$
defined below (i.e., rate $\lambda = \lambda(h)$ effective only 
when the process is in ${\cal S}$).
The idea behind the use of such a time $T_{\lambda_{\cal S}}$
comes from~\cite{BG}, which proposed the use of soft measures
and of these random times.
In comparison with the plain use of exit times
from suitable subsets of the configuration space
(approximation to a ``metastable basin'')
this gives a softer (better) way to deal 
with the escape from metastability,
also allowing a more natural use of potential theoretical tools.
For a formal statement of the mentioned convergence in law
that does not use stopping times see 
Definition~\ref{stephane}, equation~\eqref{hachem}
and formula~\eqref{fx} below, where 
$\nu$ can be taken equal to $\mu_{\Lambda_h, -, h}(\cdot | {\cal R})$
and $\lambda = \lambda(h) = e^{-\epsilon / h}$
for a small enough $\epsilon > 0$.

Now, following \cite{CGOV},
as fully detailed in \cite{OV},
we will use time averages to describe the state
of our system at earlier times.
We will identify a deterministic time scale
$\theta \ll t_{{\rm mix}, h}$ such that,
for a large class of starting measures
that will be attracted by the restricted ensemble
and for all times $t < T - \theta$,
the time averages of any observable
$f : \Omega_{\Lambda_h} \rightarrow \mathds{R}$,
$$
    A_{\theta}(t, f) 
    = {1 \over \theta} \int_t^{t + \theta} f\bigr(X_{\Lambda_h, -, h}(u)\bigl)\,du,
$$
will be close to $\mu_{\Lambda_h, -, h}(f\,|{\cal R})$
with a probability that goes to $1$ 
for a vanishing magnetic field $h$.

Before characterizing this ``large class'' of starting measures
that fall in the basin of attraction of the restricted ensemble,
we need to make precise the definitions of 
${\cal S} \subset \Omega_{\Lambda_h}$
and of the random time $T = T_{\lambda_{\cal S}}$.
The definition ${\cal S}$ is essentially symmetric
to that of ${\cal R}$
and uses 
the symmetric $B_-$ of $B_+$ 
with respect to $B_c$:
$$
    B_- = B_c - (B_+ - B_c).
$$
Note that when $B_+$ is only {\it slightly\/} supercritical
$B_-$ too is only {\it slightly\/} subcritical.
\begin{dfn}\label{stephane}
    We call ${\cal S}$ the set 
    of all configurations $\sigma$ in $\Omega_{\Lambda_h}$
    for which there is at least one external contour 
    such that a Wulff shape of volume $(B_- / h)^2$
    can fit in its interior.
\end{dfn}
\noindent We stress that,
while ${\cal R}$ refers too ``small enough'' contours
and ${\cal S}$ refers to ``large enough'' contours,
since ${\cal R}$ allows slightly supercritical contours
and ${\cal S}$ allows slightly subcritical contours,
${\cal R}$ and ${\cal S}$ do have a non-empty intersection.
These sets are actually tailored to cover
all the relevant configurations
along typical relaxation paths of the process
and allow, at the same time,
for some control of the local relaxation times
associated with the restricted processes
in ${\cal R}$ and ${\cal S}$.
Their non-empty intersection is a corollary
of such requirements.
As a consequence, we will have to use
the results of~\cite{BGM} rather than~\cite{BG};
and~\cite{BGM} will also provide,
from such bounds on local relaxation times,
the previously mentioned deterministic time scale
$\theta \ll t_{{\rm mix}, h}$.

Let now $\tau$ be a unit mean exponential time
independent of $X_{\Lambda_h, -, h}$
and let $\ell_{\cal S}(t)$ 
be the local time in ${\cal S}$ up to time $t$,
i.e., the total time spent in ${\cal S}$
by $X_{\Lambda_h, -, h}$ up to time $t$:
\begin{equation}\label{hachem}
    \ell_{\cal S}(t)
    = \int_0^t \mathds{1}\!{
        \left\{
            X_{\Lambda_h, -, h}(u) \in {\cal S}
        \right\}
    }\,du.
\end{equation}
(The law of $\ell_{\cal S}$, just as that of $X_{\Lambda_h, -, h}$,
depends on the starting distribution of $X_{\Lambda_h, -, h}$,
but, as for $X_{\Lambda_h, -, h}$, we omit it in the notation.)
$T_{\lambda_{\cal S}}$ is the time $t$
when $\ell_{\cal S}(t)$ reaches $\tau/\lambda$:
$$
    T_{\lambda_{\cal S}}
    = \min\left\{
        t \geq 0 :
        \lambda \ell_{\cal S}(t) \geq \tau
    \right\}.
$$
In other words, 
$T_{\lambda_{\cal S}}$
can be interpreted as the killing time
associated with the killing rate
defined by
$$
    \lambda_{\cal S}(\sigma)
    = \lambda \mathds{1}_{\{\sigma \in {\cal S}\}},
    \qquad \sigma \in \Omega_h.
$$
The precise value of $\lambda$
is not relevant, 
it will be enough to choose it 
in such a way to have $1 / \lambda$ large,
on the one hand,  with respect
to some ``local relaxation time in ${\cal S}$''
---more precisely, with respect
to the mixing time of the ``restricted dynamics in ${\cal S}$''---
and small,
on the other hand, with respect to the
``global mixing time'' $t_{{\rm mix}, h}$.

Let us finally introduce
two last stopping times
to state our main result.
For another parameter $\kappa > 0$
we define $T_{\kappa_{\cal R}}$
in an analogous way,
as the killing time associated
with a killing rate $\kappa_{\cal R}$,
equal to $\kappa$ in ${\cal R}$
and 0 outside of ${\cal R}$.
With $\tilde\tau$ another unit exponential time
independent of $\tau$ and $X_{\Lambda_h, -, h}$,
$T_{\kappa_{\cal R}}$ is then the time $t$
when $\ell_{\cal R}(t)$,
local time in ${\cal R}$,
reaches $\tilde \tau / \kappa$.
We call $T_{{\cal X}^c}$
the first time when $X_{\Lambda_h, -, h}$
goes outside
$$
    {\cal X} = {\cal R} \cup {\cal S}.
$$
Note that $T_{\lambda_{\cal S}}$
can also be built from a Poisson clock with rate $\lambda$
and that is independent from $X_{\Lambda_h, -, h}$:
it is the first ring time $T$ 
for which $X_{\Lambda_h, -, h}(T)$ is in ${\cal S}$.
Using another independent Poisson clock with rate $\kappa$
we can also build $T_{\kappa_{\cal R}}$ in a similar way.
$T_{\kappa_{\cal R}}$, $T_{\lambda_{\cal S}}$
and $T_{{\cal X}^c}$
are stopping times with respect to the natural filtration
associated with $X_{\Lambda_h, -, h}$ and these
two independent Poisson processes.

\begin{thm}\label{sandrine}
    For any supercritical $\beta > \beta_c$,
    any $B_{\max} > 2 B_c$, any $b < 1 / 4$
    and for all small enough $\epsilon > 0$,
    one can choose $B_+$ close enough to $B_c$
    and $\lambda = \lambda(h) = e^{- \epsilon / h}$
    for which there are
    $h_0 > 0$, $\delta > 0$ and $\delta' < \epsilon$
    such that the following holds
    for $X_{\Lambda, -, h}$ started from a probability measure $\nu$ 
    and any observable $f : \Omega_{\Lambda_h} \rightarrow \mathds{R}$.
    \begin{itemize}
    \item[i.] If $\nu = \mu_{\Lambda_h, -, h}(\cdot\,|{\cal R})$, then
        $T_{\lambda_{\cal S}} / t_{{\rm mix}, h}$
        converges in law to an exponential random variable
        of mean~$1$, i.e., for all $t > 0$,
        \begin{equation}\label{fabienne}
            \lim_{h \rightarrow 0} \mathds{P}_\nu\left(
                {T_{\lambda_{\cal S}} \over t_{{\rm mix}, h}} > t
            \right) = e^{-t}.
        \end{equation}
        Also
        \begin{equation}\label{giulio}
            \lim_{h \rightarrow 0} \mathds{P}_\nu\left(
                \theta < T_{\lambda_{\cal S}},
                \sup_{t < T_{\lambda_{\cal S}} - \theta}
                \left|
                    A_\theta(t, f) - \mu_{\Lambda_h, -, h}\bigl(f|{\cal R}\bigr)
                \right| \leq \|f\|_\infty e^{- \delta / h}
            \right) = 1,
        \end{equation}
        with
        \begin{equation}\label{cecile}
            \theta = \exp\left\{
                {1 \over h}\left(
                    {\beta A \over 2} + \delta' 
                \right)
            \right\}. 
        \end{equation}
    \item[ii.] For all $h < h_0$ it holds
        $$ 
            \Bigl|
                E_\nu\left[
                    f\Bigl(
                        X_{\Lambda_h, -, h}\bigl(
                            T_{\lambda_{\cal S}}
                        \bigr)
                    \Bigr)
                \right] - \mu_{\Lambda_h, -, h}(f)
            \Bigr| \leq \|f\|_\infty e^{- \delta / h},
        $$
        whatever the starting measure $\nu$.
    \item[iii.] If $\nu$ is such that, with $\kappa = \lambda$,
        $$
            \lim_{h \rightarrow 0} \mathds{P}_\nu\bigl(
                T_{\kappa_{\cal R}}
                < T_{\lambda_{\cal S}} \wedge T_{{\cal X}^c}
            \bigr) = 1,
        $$
        then \eqref{fabienne}--\eqref{cecile} are also in force.
    \end{itemize}
\end{thm}
\goodbreak

\par\noindent
{\bf Comments:}
\begin{itemize}
\item[i.] Equation~\eqref{fabienne} can be rewritten
    without the stopping time $T_{\lambda_{\cal S}}$,
    i.e., by referring to $X_{\Lambda_h, -, h}$ only:
    it reads
    \begin{equation}\label{fx}
        \lim_{h \rightarrow 0} \mathds{E}_\nu\left[
            e^{-\lambda \ell_S(s t_{{\rm mix},h})}
        \right] = e^{-s},
        \qquad s \geq 0.
    \end{equation}
\item[ii.] Since both $\mu_{\Lambda_h, -, h}(\cdot\,|{\cal R}_-)$,
    which does not depend on the parameters $B_+$ and $b$,
    and $\mu_{\Lambda_h, -, h}(\cdot\,|{\cal R})$ are concentrated,
    up to large deviation events, on the subset ${\cal I}$
    of ${\cal R}_-$ and ${\cal R}$ that is made of 
    configuration with invertebrate contours only,
    the same results hold with $\mu_{\Lambda_h, -, h}(\cdot\,|{\cal R}_-)$
    in place of $\mu_{\Lambda_h, -, h}(\cdot\,|{\cal R})$.
    We chose to write them with $\mu_{\Lambda_h, -, h}(\cdot\,|{\cal R})$
    for one main reason only.
    The key point of the proof will be 
    the derivation of an upper bound
    for the relaxation time of the dynamics restricted to ${\cal R}$
    (as well as the dynamics restricted to ${\cal S}$)
    and we were not able to do the same 
    with the dynamics restricted to~${\cal R}_-$.
\item[iii.] Such upper bounds will allow us to apply
    the results of \cite{BGM}.
    In particular, given a small enough $\epsilon > 0$
    we will see that one can choose
    some $B_+$ sufficiently close to $B_c$
    and $\lambda = \lambda(h) = e^{-\epsilon / h}$
    for which there are constants $C > 0$ and $\delta > 0$
    such that, if $\nu = \mu_{\Lambda_h, -, h}(\cdot\,|{\cal R})$
    or $\nu$ satisfies,
    with $\kappa = \lambda$ and for $h$ small enough,
    $$
        \mathds{P}_\nu\bigl(
            T_{\kappa_{\cal R}}
            > T_{\lambda_{\cal S}} \wedge T_{{\cal X}^c}
        \bigr) \leq  e^{-2\epsilon / h},
    $$
    then, for all $a$ such that 
    $$
        \epsilon < \beta a  < \beta A - \epsilon
    $$
    and all observable $f : \Omega_{\Lambda} \rightarrow \mathds{R}$,
    we recover
    \begin{equation}\label{bruno}
        \left|
            \mathds{E}_\nu\left[
                f\Bigl(
                    X_{\Lambda_h, -, h}\bigl(
                        e^{\beta a / h}
                    \bigr)
                \Bigr)
            \right] - \mu_{\Lambda_h, -, h}(f|{\cal R})
        \right| \leq C \|f\|_\infty e^{- \delta / h}.
    \end{equation}
    This allows, following Schonmann and Shlosman,
    for an expansion as in \eqref{alain}.
\item[iv.] The critical value for $a$ in \eqref{bruno}
    is $A$ and not $\alpha_c/\beta = A/3$
    (recall \eqref{noemie} from page \pageref{noemie}).
    The factor $1/3$ has to do with a different
    relaxation mechanism in larger boxes.
    It was first studied in \cite{DS}
    and is related both to some spatial entropy
    associated with the nucleation of a critical droplet
    and to the time needed for a supercritical droplet to invade a fixed box.
    In the infinite volume case or already in the case of a
    large domain $\Lambda$ of exponentially large volume $e^{C/h}$
    with a large enough $C$, not only the asymptotic
    value of the mean ``transition time to equilibrium'' 
    would change; it is not
    clear anymore whether we
    should expect its law to be asymptotically exponential:
    to an exponential random time needed to nucleate a critical droplet
    we should add another time of the same
    logarithmic scale order (the time needed to invade the given box),
    and prefactors enter the game at this point.
    The asymptotic exponential law would survive
    if the prefactor associated with the nucleation of the
    critical droplet is dominant.

\item[v.] The condition $B_{\max} > 2 B_c$
    ensures that the volume is large enough
    for the positive magnetic field to overcome
    the effect of the negative boundary condition,
    in such a way that the plus phase invades
    the whole box at equilibrium.
\item[vi.] The restriction on the shape of the domain
    is technical and will simplify the proof.
    It avoids in particular 
    a description of typical equilibrium configurations
    in more general domains.
\item[vii.] Theorem~\ref{sandrine} allows us
    to consider more general starting distributions
    than in \cite{SS}. This is due to the fact
    that controlling the local relaxation time 
    in ${\cal R}$ and ${\cal S}$,
    we will not have to rely on the monotonicity
    of $X$ in the same way.
\end{itemize}

Thinking of a slowly changing magnetic field
as in the hysteresis phenomena,
it is natural to consider starting distributions like
$\mu_{\Lambda_h, -, h'}(\cdot\,|{\cal R}_{h'})$
associated with a different magnetic field $h'$,
but with the same domain $\Lambda_h$.
This is one possibility considered in
the following corollary of Theorem~\ref{sandrine}.
The other possibility we consider
in this corollary is that of the canonical ensemble
associated with a small enough magnetization
$$
    M: \omega \in \Omega_{\Lambda_h} \mapsto \sum_{x \in \Lambda_h} \omega(x),
$$
namely
$\mu_{\Lambda_h, -, h}(\cdot\,|{\cal R} \hbox{ and } M > m(B_{\max} / h)^2)$
with $m < m^*_\beta[2 (B_c / B_{\max})^2 - 1]$.
This upper bound corresponds to the magnetization
of a critical Wulff shape droplet of plus phase
in the minus phase.
\begin{cor}\label{maxime}
    Let $\epsilon > 0$, $c > 0$ and
    $m < m^*_\beta[2 (B_c / B_{\max})^2 - 1]$
    associated with $\beta > \beta_c$
    and $B_{\max} > 2 B_c$.
    If $\nu = \mu_{\Lambda_h, -, h'}(\cdot\,|{\cal R}_{h'})$
    associated with $h' = c h$
    or $\nu = \mu_{\Lambda_h, -, h}(\cdot\,| {\cal R}, M > m(B_{\max} / h)^2)$,
    then there are $B_+ > B_c$, $\lambda = \lambda(h) = e^{-\epsilon / h}$,
    $\delta > 0$, $\delta' < \epsilon$ and $C > 0$
    such that \eqref{fabienne}--\eqref{cecile}
    and \eqref{bruno} hold for any observable
    $f : \Omega_{\Lambda_h} \rightarrow \mathds{R}$
    and if $\epsilon < \beta a  < \beta A - \epsilon$.
\end{cor}

In the next section we introduce
a collection of tools
for the proof of Theorem~\ref{sandrine}, Proposition~\ref{amandine}
and Corollary~\ref{maxime}.
This includes in particular
static estimates,
for which the main references are
\cite{SS}, \cite{DKS}, \cite{Pfi}, \cite{Iof1} and \cite{Iof2}
and dynamical techniques,
for which the main references are
\cite{Sin} and \cite{Mar}.
We use the former in Section~\ref{veronique}
to give lower bounds on the transition time
to equilibrium.
We use the latter in Section~\ref{julia}
to give upper bounds on local relaxation times.
This is the key point of the proof:
we show in the last part of Section~\ref{hassan}
how to use the results of \cite{BGM}
to obtain from such estimates 
an equivalent of Theorem~\ref{sandrine},
Proposition~\ref{amandine}
and Estimate~\eqref{bruno}
for the restriction $X$
of our process $X_{\Lambda_h, -, h}$
to ${\cal X} = {\cal R} \cup {\cal S}$,
and we explain how to reduce the study
of $X_{\Lambda_h, -, h}$ to that of $X$. 
We finally prove Theorem~\ref{sandrine},
Proposition~\ref{amandine}
and Corollary~\ref{maxime} in Section~\ref{erwan}.
From now on we will always assume our
fixed parameters $\beta$ and $B_{\max}$
to be respectively larger than
the critical inverse temperature $\beta_c$ and $2 B_c$.

\section{Tools, notation and strategy}\label{hassan}
\subsection{Wulff shape and surface tension}\label{aldo}
In order to define the surface tension
in a direction orthogonal to the unitary vector
${\bf n} = (\cos \theta, \sin \theta)$
for $\theta \in [0, 2\pi]$,
we have to consider the Ising model in a square box
$\Lambda(L) = [-L , L]^2$ with boundary condition 
$$
    \eta_\theta(x) = \left\{
        \begin{array}{ll}
            +1 & \hbox{if $u \cos \theta + v \sin \theta \leq 0$,}\\
            -1 & \hbox{if $u \cos \theta + v \sin \theta > 0$,}
        \end{array}
    \right.
    \qquad x = (u, v) \in \mathds{Z}^2.
$$
In a contour description of the configurations
that are associated with such a boundary condition,
one contour, on the dual lattice,
must join two points that are close
to $y(L)$ and $z(L)$,
which are the two points
where the boundary of the box $[-L , L]^2$
intersects the straight line
that goes through the origin and admits
${\bf n}$ as normal vector.
The {\it surface tension} in the direction
of this straight line is
$$
    \tau(\theta)
    = \lim_{L \rightarrow +\infty}
    - {1 \over \beta \|y(L) - z(L)\|_2}
    \ln {
        Z_{\Lambda(L), \eta_\theta, 0}
    \over
        Z_{\Lambda(L), +, 0}
    }
    \,,
$$
with $Z_{\Lambda(L), +, 0}$
the partition functions associated
with the Ising model in $\Lambda(L)$, 
with uniform plus boundary condition and
without magnetic field.
Thus, the surface tension $\tau(\theta)$ 
is the free energy per unit length
of an interface between the plus and minus phase
in the direction orthogonal to $\bf n$.
It is positive and finite 
for subcritical temperature
$1 / \beta < 1 / \beta_c$.

We then define the {\it surface free energy}
of any rectifiable $\gamma \subset \mathds{R}^2$
that is the boundary
of a simply connected domain $D \subset \mathds{R}^2$
by the quantity
\begin{equation}\label{marina}
    {\cal W}(\gamma) = \oint_\gamma \tau(\theta_s) \,ds,
\end{equation}
with $\theta_s$ the direction of the external normal,
i.e., which points outside $D$,
at the curvilinear abscissa $s$.
We will refer to ${\cal W}$
as the {\it Wulff functional}.
The {\it Wulff shape} has a boundary
that minimizes this quantity
among all the rectifiable boundaries
of domains with a given volume.
It is defined for $\rho > 0$
and up to dilatation and translation by
\begin{equation}\label{benjamin}
    W_\rho = \bigcap_{\theta \in [0, 2\pi]} \Bigl\{
        x = (u, v) \in \mathds{R}^2 :
        u \cos \theta + v \sin \theta \leq \rho \tau(\theta)
    \Bigr\}.
\end{equation}
As a consequence of the symmetries
of $\tau$ that are inherited from those of the lattice,
$W_\rho$ is invariant by rotations of angle $\pi / 2$.
We will simply write $W$,
without the index $\rho$,
when $\rho$ is chosen
in such a way that $W_\rho$ has a volume 
equal to one.

The {\it support function}
with respect to the origin 0 of
the convex set $W_\rho \ni 0$ is actually $\rho\tau$,
i.e.,
$$
    \rho \tau(\theta)
    = \max_{x = (u, v) \in W_\rho} u \cos \theta + v \sin \theta,
    \qquad \theta \in [0, 2\pi].
$$
This is a consequence of the triangular inequality:
for $x$, $y$ and $z$ in $\mathds{R}^2$, if
$$
    {\bf n}_z = (\cos \theta_z, \sin \theta_z),
    \qquad
    {\bf n}_x = (\cos \theta_x, \sin \theta_x),
    \qquad\hbox{and}\qquad
    {\bf n}_y = (\cos \theta_y, \sin \theta_y)
$$
are the external normals to the three sides
$[x,y]$, $[y, z]$ and $[z, x]$
of the triangle $xyz$, 
then 
$$
    \|x - z\|_2 \tau(\theta_y) 
    \leq \|x - y\|_2 \tau(\theta_z) 
    + \|y - z\|_2 \tau(\theta_x)
$$
(see Section~4.21 in \cite{DKS}).

Let us denote by $|D|$ 
the volume of any measurable domain
$D \subset \mathds{R}^2$.
Then Bonnesen's inequality says that 
for any such domain $D$
with a rectifiable boundary $\gamma$,
choosing $\rho$ in such a way that
$|W_\rho| = |D|$,
it holds
\begin{equation}\label{sebastien}
    {\cal W}(\gamma) 
    \geq {\cal W}(\partial W_\rho) \sqrt{
        1 + \left(
            {\alpha_{out} - \alpha_{in} \over 2}
        \right)^2
    },
\end{equation}
where $\partial W_\rho$
stands for the boundary of $W_\rho$,
and $\alpha_{out}$, respectively $\alpha_{in}$,
is the smallest, 
respectively the largest, $\alpha$
for which a translate of $\alpha W_\rho$
contains, respectively is contained in, $D$.
In the case where $D$ is a convex set,
this is proven in \cite{Fla}
by counting the mean number of intersections
between $\gamma$ and the border
of a random translate $X + \alpha W_\rho$,
for $\alpha \in [\alpha_{in}, \alpha_{out}]$
and $X$ uniformly chosen in $D - \alpha W_\rho$.
Flanders proves in this way Blaschke's inequality
$$
    \alpha^2 |W_\rho| - \alpha \rho {\cal W}(\gamma) + |D| \leq 0,
    \qquad \alpha \in [\alpha_{out}, \alpha_{in}],
$$
with equality in the case $\alpha_{out} = \alpha_{in}$.
This gives a lower bound
on the distance between the two roots
of this polynomial of degree two in $\alpha$,
i.e., a lower bound on its discriminant,
which leads, together with the equality 
for $D = W_\rho$,
\begin{equation}\label{ben}
    \rho {\cal W}(\partial W_\rho)
    = 2 |W_\rho|,
\end{equation}
to inequality~\eqref{sebastien}.
In the case where $D$ is not a convex set,
these inequalities are not a direct consequence
of those of the convex case, but the same strategy
can be followed even though the computation 
of this mean intersection number is more delicate.
In \cite{DKS} the authors adapt an argument 
from \cite{O1, O2} to cover the case
of a non-convex simply connected $D$
(see Section~2.5 in \cite{DKS}).
We will use this result,
rewriting it with the following notation.
With $\rho$ and $B$ such that $|D| = |W_\rho| = B^2$,
we set $B_{in} = \alpha_{in} B$ and $B_{out} = \alpha_{out} B$,
with $\alpha_{out}$ and $\alpha_{in}$ as above
so that $B_{in}^2$ is the volume of the largest Wulff shape
that fits in $D$ and $B_{out}^2$ that of the smallest Wulff
shape that contains it.
We denote by $w_\beta$
the {\it surface free energy of the unitary volume Wulff shape} $W$,
so that
$$
    {\cal W}(\partial W_\rho) = w_\beta B
$$
and, as a consequence of \eqref{ben},
\begin{equation}\label{charles}
    B = {w_\beta \over 2} \rho.
\end{equation}
\begin{prp}[Blaschke's inequalities \cite{DKS}]\label{william}
    For any simply connected domain $D \subset \mathds{R}^2$
    with a rectifiable boundary $\gamma$
    it holds 
    $$
        {\cal W}(\gamma) \geq {w_\beta \over 2}\left(
            {|D| \over B_{in}} + B_{in}
        \right)
        \qquad\hbox{and}\qquad
        {\cal W}(\gamma) \geq {w_\beta \over 2}\left(
            {|D| \over B_{out}} + B_{out}
        \right).
    $$
\end{prp}
\goodbreak

We will also need two simple consequences
of the Wulff construction
from the support function $\rho\tau$.
\begin{lmm}\label{raphael}
    If two translates of possibly different size Wulff shapes
    $x_1 + W_{\rho_1}$ and $x_2 + W_{\rho_2}$,
    of volume $B_1^2$ and $B_2^2$,
    have a non-empty intersection,
    then their union is contained in some Wulff shape
    $x_0 + W_{(\rho_1 + \rho_2)}$ of volume $(B_1 + B_2)^2$.
\end{lmm}
\par\noindent
{\bf Proof:} Since $x_1 + W_{\rho_1}$ and $x_2 + W_{\rho_2}$
have a non-empty intersection, 
there are $w_1$ and $w_2$ in $W_{\rho_1}$ and $W_{\rho_2}$
such that $x_1 + w_1 = x_2 + w_2$, i.e.,
$$
    x_1 - w_2 = x_2 - w_1.
$$
This means thats $x_1 - W_{\rho_2}$
and $x_2 - W_{\rho_1}$
also have a non-empty intersection.
Let us then choose
$$
    x_0 \in \Bigl(
        x_1 - W_{\rho_2}
    \Bigr) \cap \Bigl(
        x_2 - W_{\rho_1}
    \Bigl).
$$
We have $x_1 - x_0 \in W_{\rho_2}$,
then, writing $(u_1, v_1)$ and $(u_0, v_0)$
for the coordinates in $\mathds{R}^2$
of $x_1$ and $x_0$,
it holds 
$$
    (u_1 - u_0) \cos \theta + (v_1 - v_0) \sin \theta \leq \rho_2 \tau(\theta)
$$
for any $\theta \in [0, 2\pi]$.
For any $x = (u, v)$ in $x_ 1 + W_{\rho_1}$
we also have    
$$
    (u - u_1) \cos \theta + (v - v_1) \sin \theta \leq \rho_1 \tau(\theta),
$$
hence
$$
    (u - u_0) \cos \theta + (v - v_0) \sin \theta \leq (\rho_1 + \rho_2) \tau(\theta).
$$
This shows that
$x_1 + W_{\rho_1}$ is contained in $x_0 + W_{(\rho_1 + \rho_2)}$,
and we can check in the same way that 
$x_2 + W_{\rho_2}$ is contained in $x_0 + W_{(\rho_1 + \rho_2)}$.
\qed

\medskip\par
The previous proof only use the fact
that the Wulff shape is a convex set,
to which one can associate a support function to describe it.
The last lemma of this section
uses by contrast the symmetries of the lattice,
namely the fact that $W = -W$,
i.e., that $\rho\tau$ is $\pi$-periodic.
\begin{lmm}\label{kata}
    Given $B_2 > B_1$, the largest Wulff shapes
    to fit in the annulus $B_2 W \setminus B_1 W$
    of volume $B_2^2 - B_1^2$
    have a volume $B_0^2 = (B_2 - B_1)^2 / 4$.
\end{lmm}

\par\noindent
{\bf Proof:} The Wulff shape construction
from the $\pi$-periodical support function $\rho\tau$
implies that, for any positive $\rho_1$ and $\rho_0$
the union of $W_{\rho_1}$
with all the externally tangent Wulff shapes
$$
    x + W_{\rho_0},
    \qquad x \in \partial W_{\rho_1 + \rho_0},
$$
is the Wulff shape $W_{\rho_1 + 2 \rho_0}$.
We get the desired result by choosing $\rho_1$ and $\rho_0$
in such a way that, with $\rho_2 = \rho_1 + 2 \rho_0$,
$$
    W_{\rho_1} = B_1 W
    \qquad {\rm and} \qquad
    W_{\rho_2} = B_2 W,
$$
i.e.,
$$
    \rho_1 = 2 B_1 / w_\beta
    \qquad {\rm and} \qquad
    \rho_2 = 2 B_2 / w_\beta
$$
so that 
$$
    \rho_0 = {\rho_2 - \rho_1 \over 2} = {B_2 - B_1 \over w_\beta}
$$
and
$$
    B_0^2
    = \left(\rho_0 w_\beta \over 2\right)^2
    = \left(B_2 - B_1 \over 2\right)^2.
$$
\qed

\subsection{Random paths, flows and block flows}\label{nelly}
Given a generic irreducible Markov process $Y$
on a finite configuration space ${\cal Y}$
with generator\footnote{
    The index ${\cal Y}$, rather than $Y$,
    in the notation ${\cal L}_{\cal Y}$ can seem unnatural 
    since the generator depends on the whole process
    and not only on the configuration space,
    but we are foreseeing here a later more natural notation,
    in accordance with~\cite{BGM}.
}
${\cal L}_{\cal Y}$ 
$$
    ({\cal L}_{\cal Y}f)(\sigma) = \sum_{\sigma' \in {\cal Y}} w(\sigma, \sigma')\bigl[
        f(\sigma') - f(\sigma)
    \bigr],
    \qquad f: {\cal Y} \mapsto \mathds{R},
    \qquad \sigma \in {\cal Y},
$$
a {\it path} $\pi$ 
is a finite sequence $(\sigma_0, \sigma_1, \dots, \sigma_l)$
of configurations in ${\cal Y}$
such that $w(\sigma_k, \sigma_{k + 1}) > 0$
for all $k < l$.
The {\it length} $|\pi|$ of such a path $\pi$
is the integer $l$.
If $e = (\sigma, \sigma')$ belongs to the edge set ${\cal E}$
associated with $Y$,
i.e., if $\sigma$ and $\sigma'$ are distinct configurations
such that $w(\sigma, \sigma') > 0$,
we write $e \in \pi$ if there is $k < |\pi|$
such that $e = (\sigma_k , \sigma_{k + 1})$.
We will also write $\sigma \in \pi$ if there is $k \leq |\pi|$
such that $\sigma = \sigma_k$.

Random paths $\Pi$ are associated with {\it flows},
i.e., with functions $\psi : {\cal E} \rightarrow \mathds{R}$,
such that 
$$
    \psi(\sigma, \sigma') = - \psi(\sigma', \sigma),
    \qquad (\sigma, \sigma') \in {\cal E}.
$$
Indeed, with $\Pi = (Y_0, \dots, Y_{|\Pi|})$,
$\Pi_- = Y_0$ and $\Pi^+ = Y_{|\Pi|}$
we get such an antisymmetric function
by setting
$$
    \psi(\sigma, \sigma') = E\left[
        \sum_{k < |\Pi|} \mathds{1}\!\bigl\{
            (\sigma, \sigma') = (Y_k, Y_{k + 1})
        \bigr\} - \mathds{1}\!\bigl\{
            (\sigma', \sigma) = (Y_k, Y_{k + 1})
        \bigr\}
    \right]
$$
and we note that, for all $\sigma$ in ${\cal Y}$,
$$
    {\rm div}_\sigma \psi
    = \sum_{\sigma' \in {\cal Y}} \psi(\sigma, \sigma')
    = P\bigl(\Pi_- = \sigma\bigr) - P\bigl(\Pi^+ = \sigma\bigr).
$$
In particular, if there are two disjoint subsets
${\cal A}$ and ${\cal B}$ of ${\cal Y}$ 
such that $\Pi_- \in {\cal A}$ 
and $\Pi^+ \in {\cal B}$ with probability one,
then $\psi$ is a
{\it unitary flow from ${\cal A}$ to ${\cal B}$},
i.e., such that
$$
    {\rm div}_\sigma \psi > 0 \Rightarrow \sigma \in {\cal A},
    \qquad
    {\rm div}_\sigma \psi < 0 \Rightarrow \sigma \in {\cal B}
    \qquad{\rm and}\qquad
    \sum_{\sigma \in {\cal A}} {\rm div_\sigma \psi}
    = 1
    = - \sum_{\sigma \in {\cal B}} {\rm div_\sigma \psi}.
$$ 

Sinclair proved in~\cite{Sin}
that if $Y$ is reversible with respect to some probability measure $\mu_{\cal Y}$,
i.e., if the {\it conductances}
$$
    c(\sigma, \sigma') = \mu_{\cal Y}(\sigma) w(\sigma, \sigma'),
    \qquad \sigma, \sigma' \in {\cal Y},
$$
are symmetrical,
then for any random path $\Pi$ 
with starting and ending configurations
that are independently distributed
according to $\mu_{\cal Y}$,
it holds
$$
    {1 \over \gamma_{\cal Y}}
    \leq \max_{e \in {\cal E}} {1 \over c(e)} P\bigl(
        e \in \Pi
    \bigr) E\left[
        |\Pi|
        \bigm| e \in \Pi
    \right]
    \leq \max_{e \in {\cal E}} {1 \over c(e)} P\bigl(
        e \in \Pi
    \bigr) \bigl\||\Pi|\bigr\|_\infty
$$
with $1 / \gamma_{\cal Y}$ the relaxation time of $Y$,
i.e.,
\begin{equation}\label{elizabeth}
    \gamma_{\cal Y} = \min_{{\rm Var}_{\mu_{\cal Y}}(f) \neq 0} {
        {\cal D}(f)
    \over
        {\rm Var}_{\mu_{\cal Y}}(f)
    }\,,
\end{equation}
where ${\cal D}$ is the {\it Dirichlet form}
defined by
\begin{equation}\label{baptiste}
    {\cal D}(f)
    = {1 \over 2} \sum_{\sigma, \sigma' \in {\cal Y}} c(\sigma, \sigma) \bigl[
        f(\sigma) - f(\sigma')
    \bigr]^2.
\end{equation}
In particular, if there is a lower bound 
$$
    w(\sigma, \sigma') \geq w_{\min}, 
    \qquad (\sigma, \sigma') \in {\cal E},
$$
then 
\begin{equation}\label{caroline}
    {1 \over \gamma_{\cal Y}} 
    \leq {\bigl\||\Pi|\bigr\|_\infty \over w_{\min}}
    \max_{(\sigma, \sigma') \in {\cal E}} {
        P\bigl(
            (\sigma, \sigma') \in \Pi
        \bigr)
    \over
         \mu_{\cal Y}(\sigma) \vee \mu_{\cal Y}(\sigma')
    }
\end{equation}

The simplest way to obtain upper bounds for relaxation times
with a random path $\Pi$ is to build
for each $\sigma$ and $\sigma'$ in ${\cal Y}$
a deterministic path $\pi_{\sigma, \sigma'}$,
usually referred to as {\it canonical path},
and set $\Pi = \pi_{\sigma, \sigma'}$
with probability $\mu_{\cal Y}(\sigma) \mu_{\cal Y}(\sigma')$.
Martinelli gave in \cite{Mar}
an upper bound for the relaxation time
of the Glauber dynamics $X_{\Lambda(L), +, 0}$
in the square box $\Lambda(L) = [-L, L]^2$
by introducing a ``block dynamics'',
bounding its mixing time by a coupling argument
and bounding the relaxation time
of the Glauber dynamics in each block
with such canonical paths.
For a block covering of
$$
    \Lambda(L) = \bigcup_{j < k} \Lambda_j
$$
by partially overlapping rectangular blocks $\Lambda_j$
of size $L \times L^{\epsilon + 1 / 2}$,
the associated block dynamics
update at rate one the current configuration $\sigma$
according to $\mu_{\Lambda_j, \sigma, 0}$.
He bounded the mixing time of this block
dynamics by using its monotonicity properties.
And as far as the relaxation time of each 
$X_{\Lambda_j, \eta, 0}$ is concerned,
he built the canonical path $\pi_{\sigma, \sigma'}$
from any $\sigma$ in $\Omega_{\Lambda_j}$
to any $\sigma'$ in the same configuration space
by ordering,
independently of $\sigma$ and $\sigma'$,
the sites of the rectangle $\Lambda_j$
and flipping the spins from their value in $\sigma$
to their value in $\sigma'$ in this prescribed order.
This order $\preceq$ had the key property
that for any $x$ in $\Lambda = \Lambda_j$,
with
$$
    \Lambda^{\prec x}
    = \left\{
        y \in \Lambda 
        : y \prec x
    \right\},
$$
$$
    \Lambda^{\succeq x}
    = \left\{
        z \in \Lambda 
        : z \succeq x
    \right\}
$$
and
$$
    \partial\Lambda^{\prec x}
    =\left\{
        (y, z) \in \Lambda^{\prec x} \times \Lambda^{\succeq x}
        : \|y - z\| = 1
    \right\},
$$
$|\partial\Lambda^{\prec x}|$ was of the same order 
has the {\it shorter} side of the rectangle $\Lambda$.
Martinelli could then use a practical version
of the following abstract lemma.
\begin{lmm}\label{fabio}
    For any finite box $\Lambda \subset \mathds{Z}^2$,
    any order $\preceq$ on $\Lambda$,
    any boundary condition $\eta \in \Omega_{\mathds{Z}^2}$,
    any configuration $\sigma_0$ in $\Omega_\Lambda$
    and any site $x$ in $\Lambda$,
    it holds
    $$
        {1 \over \mu_{\Lambda, \eta, h}(\sigma_0)}
        \sum_{\sigma, \sigma' \in \Omega_\Lambda}
        \mu_{\Lambda, \eta, h}(\sigma) \mu_{\Lambda, \eta, h}(\sigma')
        \mathds{1}\!\left\{
            (\sigma_0, \sigma_0^x) \in \pi^\preceq_{\sigma, \sigma'}
        \right\}
        \leq \exp\left\{2\beta |\partial\Lambda^{\prec x}|\right\}
    $$
    where $\pi^\preceq_{\sigma, \sigma'}$
    stands for the canonical path from $\sigma$ 
    to $\sigma'$ associated with the order $\preceq$.
\end{lmm}

\par\noindent
{\bf Proof:} Following the computation 
made in \cite{Mar}, Section~2,
denoting, for any $\sigma^{\prec x} \in \Omega_{\Lambda^{\prec x}}$
and $\sigma^{\succeq x} \in \Omega_{\Lambda^{\succeq x}}$,
by $\sigma^{\prec x} \cdot \sigma^{\succeq x}$
the configuration of $\Omega_{\Lambda}$ 
that coincides with $\sigma^{\prec x}$ in $\Lambda^{\prec x}$
and $\sigma^{\succeq x}$ in $\Lambda^{\succeq x}$,
and recalling the presence of the somewhat unusual factor $1/2$
in our Hamiltonian definition,
we have
\begin{align*}
    &{1 \over \mu_{\Lambda, \eta, h}(\sigma_0)}
    \sum_{\sigma, \sigma' \in \Omega_\Lambda}
    \mu_{\Lambda, \eta, h}(\sigma) \mu_{\Lambda, \eta, h}(\sigma')
    \mathds{1}\!\left\{
        (\sigma_0, \sigma_0^x) \in \pi^\preceq_{\sigma, \sigma'}
    \right\}\\
    &\qquad \leq {1 \over \mu_{\Lambda, \eta, h}(\sigma_0)}
    \sum_{\sigma, \sigma' \in \Omega_\Lambda}
    \mu_{\Lambda, \eta, h}(\sigma) \mu_{\Lambda, \eta, h}(\sigma')
    \mathds{1}\!\bigl\{
        \sigma|_{\Lambda^{\succeq x}} = \sigma_0|_{\Lambda^{\succeq x}},
        \sigma'|_{\Lambda^{\prec x}} = \sigma_0|_{\Lambda^{\prec x}}
    \bigr\}\\
    &\qquad = {1 \over Z_{\Lambda, \eta, h}}
    \sum_{
        \scriptstyle \sigma^{\prec x} \in \Omega_{\Lambda^{\prec x}},
    \atop 
        \scriptstyle \sigma^{\succeq x} \in \Omega_{\Lambda^{\succeq x}}
    }{
        \exp\left\{
            -\beta H_{\Lambda, \eta, h}\left(
                \sigma^{\prec x} \cdot \sigma_0|_{\Lambda^{\succeq x}}
            \right)
            -\beta H_{\Lambda, \eta, h}\left(
                \sigma_0|_{\Lambda^{\prec x}} \cdot \sigma^{\succeq x} 
            \right)
        \right\}
    \over
        \exp\left\{
            -\beta H_{\Lambda, \eta, h}\left(
                \sigma_0|_{\Lambda^{\prec x}} \cdot \sigma_0|_{\Lambda^{\succeq x}}
            \right)
        \right\}
    }\\
    &\qquad = {1 \over Z_{\Lambda, \eta, h}}
    \sum_{
        \scriptstyle \sigma^{\prec x} \in \Omega_{\Lambda^{\prec x}},
    \atop 
        \scriptstyle \sigma^{\succeq x} \in \Omega_{\Lambda^{\succeq x}}
    }{
        \exp\left\{
            -\beta H_{\Lambda, \eta, h}\left(
                \sigma^{\prec x} \cdot \sigma_0|_{\Lambda^{\succeq x}}
            \right)
            -\beta H_{\Lambda, \eta, h}\left(
                \sigma_0|_{\Lambda^{\prec x}} \cdot \sigma^{\succeq x} 
            \right)
        \right\}
    \over
        \exp\left\{
            -\beta H_{\Lambda, \eta, h}\left(
                \sigma_0|_{\Lambda^{\prec x}} \cdot \sigma_0|_{\Lambda^{\succeq x}}
            \right)
            -\beta H_{\Lambda, \eta, h}\left(
                \sigma^{\prec x} \cdot \sigma^{\succeq x}
            \right)
        \right\}
    } e^{
        -\beta H_{\Lambda, \eta, h}\left(
            \sigma^{\prec x} \cdot \sigma^{\succeq x}
        \right)
    }\\
    &\qquad \leq {1 \over Z_{\Lambda, \eta, h}}
    \sum_{
        \scriptstyle \sigma^{\prec x} \in \Omega_{\Lambda^{\prec x}},
    \atop 
        \scriptstyle \sigma^{\succeq x} \in \Omega_{\Lambda^{\succeq x}}
    }\exp\left\{
        2 \beta |\partial\Lambda^{\prec x}|
    \right\} e^{
        -\beta H_{\Lambda, \eta, h}\left(
            \sigma^{\prec x} \cdot \sigma^{\succeq x}
        \right)
    }\\
    &\qquad = \exp\left\{
        2 \beta |\partial\Lambda^{\prec x}|
    \right\}.
\end{align*}
\qed

\medskip\par
The crucial spectral gap estimates of the present paper
(see Section~\ref{julia})
rely on the following observation:
as far as leading orders are concerned,
Martinelli's lower bound on $\gamma_{\Lambda(L), +, 0}$
can be obtained by direct application of formula \eqref{caroline}.
To do so one has to build a random path $\Pi$ 
with starting and ending configurations independently distributed
according to $\mu_{\Lambda(L), +, 0}$.
Equivalently one has to build,
for each $\sigma$ and $\sigma'$ in $\Omega_{\Lambda(L)}$,
a random path $\Pi_{\sigma, \sigma'}$ and set $\Pi = \Pi_{\sigma, \sigma'}$
with probability
$\mu_{\Lambda(L), +, 0}(\sigma) \mu_{\Lambda(L), +, 0}(\sigma')$.
Here is a block dynamic inspired way to build a suitable 
$\Pi_{\sigma, \sigma'}$ from two random paths $\Pi_{\sigma}$
and $\Pi_{\sigma'}$ starting from $\sigma$ and $\sigma'$,
respectively.
From $\sigma$ we build a sequence of $k$ random configurations
that we will call ``milestones''
$M_1$, $M_2$, \dots, $M_k$ in $\Omega_{\Lambda(L)}$.
We set $M_0 = \sigma$, call it our first milestone,
and build from each milestone $M_j$, with $j < k$,
the next milestone $M_{j + 1}$ 
by setting $M_{j + 1}|_{\Lambda_j^c} = M_j|_{\Lambda_j^c}$
and drawing $M_{j + 1}|_{\Lambda_j}$ according to $\mu_{\Lambda_j, M_j, 0}$.
Next, we use,
in each block $\Lambda_j$,
a canonical path of the single spin flip Glauber dynamics
to connect $M_j$ with $M_{j + 1}$;
this defines our random path $\Pi_\sigma$ 
and we build $\Pi_{\sigma'}$ in an analogous way from $\sigma'$.
Consider now, with obvious notation, the event
$$
    E_{\sigma, \sigma'} = \bigl\{
        M_k = M'_k
    \bigr\}.
$$
When $E_{\sigma, \sigma'}$ occurs we can build $\Pi_{\sigma, \sigma'}$
by concatenation of $\Pi_\sigma$, from $\sigma$ to $M_k$,
and the reversed path $\Pi_{\sigma'}$, from $M'_k = M_k$
to $\sigma'$.
From the conditional probability associated with $E_{\sigma, \sigma'}$
we get a random path $\Pi_{\sigma, \sigma'}$ from $\sigma$ to $\sigma'$,
then a random path $\Pi$ 
with starting and ending configurations independently distributed
according to the equilibrium distribution. 
When used in formula~\eqref{caroline}, estimating the relaxation time
$1 / \gamma_{\Lambda(L), +, 0}$ boils down, through DLR equations,
to computing a uniform lower bound on $P(E_{\sigma, \sigma'})$
that, in turns, can be obtained with the very same arguments
used by Martinelli for controlling the mixing time of the block dynamics.

This is nothing but an alternative way 
of articulating Martinelli's ideas.
But in doing so we gain some flexibility:
there is no need anymore to define any block {\it dynamic},
we only need to build suitable sequences of milestones
for which we can give a uniform lower bound
on the probability of such events $E_{\sigma, \sigma'}$
that are contained in $\{M_k = M'_k\}$
(we used here the latter event to define $E_{\sigma, \sigma'}$
but we will later require more from such events;
and the inclusion will be again needed
for allowing a similar construction of $\Pi_{\sigma, \sigma'}$
from those of $\Pi_{\sigma}$ and $\Pi_{\sigma'}$).
In particular the box $\Lambda_j$ used to build
$M_{j + 1}$ from $M_j$
{\it can now depend in some way of $M_j$}.
We will use this slightly different strategy
and the flexibility it allows
to control the local relaxation times
of $X_{\Lambda_h, -, h}$ restricted 
to ${\cal R}$ and ${\cal S}$.
We will refer to such milestone built
random paths $\Pi_\sigma$, $\Pi_{\sigma'}$
or $\Pi_{\sigma, \sigma'}$
and their associated flows as
``block paths'' and ``block flows''.
We will also use such a block flow
to estimate the soft capacity presented
in Section~\ref{zakia}.

\subsection{Free energy estimates}
In this section we closely follow Schonmann and Shlosman.
In~\cite{SS} they derived a number of free energy estimates
that we have to slightly adapt 
to deal with some particular ``annular droplets''
(see estimate \eqref{adrien_4} in Lemma~\ref{peter}).
In this respect we need a slightly stronger theory,
but the extension is straightforward
and we only write in this section
those technical points
for which we need a slightly different writing.
Also there are a few estimates
for which we only need a weaker form
than in~\cite{SS} (for example,
the probability appearing
in estimate \eqref{adrien_3} of Lemma~\ref{peter}
is actually shown to be close to one in~\cite{SS}).
For these estimates, their stronger stronger forms 
in~\cite{SS} derive from stability results
in~\cite{DKS}, which, in turn, are based on
Blaschke's inequalities of Section~\ref{hassan}.
We also use Blaschke's inequalities in this paper,
but for other purposes,
mainly in Section~\ref{veronique}
and also in proving estimate \eqref{adrien_4}
in Lemma~\ref{peter}.

The key objects introduced in \cite{DKS}
to make sense of a macroscopic  
(on length scale $1 / h$)
or even mesoscopic 
(on length scale $1 / h^b$,
with $b < 1 / 4$ as mentioned earlier)
notion of free energy
are the skeletons associated with vertebrate contours,
i.e., contours with more than $1 / h^{2b}$ sites
in their interior.
Following \cite{SS} to build them,
we will be closer to their construction in \cite{Pfi}.

Let $r$ be a positive number that is smaller than $b / 2$
an consider a configuration $\sigma$ in 
$$
    \Omega_{\Lambda, -}
    = \left\{
        \sigma \in \Omega_{\mathds{Z}^2} 
        : \sigma(x) = -1 \hbox{ for all $x \not\in \Lambda$}
    \right\},
$$
which we identified with $\Omega_\Lambda$,
for a finite domain $\Lambda \subset \mathds{Z}^2$.
A {\it skeleton\/} associated with a vertebrate contour $\Gamma$ of $\sigma$
is a possibly self-intersecting polygon $\gamma \subset \mathds{R}^2$
such that
\begin{itemize}
\item[i.] the ordered vertices of which
    are consecutive points on $\Gamma$ with the same order
    (for one of the two possible orientation of $\Gamma$);
\item[ii.] the side lengths of which
    lie between $1 / (12 h^r)$
    and $1 / h^r$;
\item[iii.] such that the Haussdorff distance between 
    $\Gamma$ and $\gamma$ is smaller than or equal to 
    $1 / h^r$.
\end{itemize}
In what follows we will assume 
that we have an algorithm 
to assign such a skeleton $\gamma$
to any vertebrate contour $\Gamma$,
so that we can refer to {\it the} collection 
of skeleton $S = (\gamma_j : j < k)$
associated with the collection
$G = (\Gamma_j : j < k)$
of the vertebrate contours of a configuration 
$\sigma$ in $\Omega_{\Lambda, -}$.
Such an algorithm is described in
\cite{DKS}, Section~5.11,
under the assumption that the diameter
of $\Gamma$ is larger than $1 / h^r$,
which is ensured by the fact that $\Gamma$
is vertebrate.
We will refer to this algorithm
as the function
$S_{\Lambda, -, h} : \sigma \in \Omega_\Lambda \mapsto S$,
which we will see as a random variable 
on the probability space $(\Omega_\Lambda,\, \mu_{\Lambda, -, h})$.
Differently from the notation of \cite{SS},
$G$ and $S$ are not associated
with external only vertebrate contours,
but with all the vertebrate contours of a configuration $\sigma$. 
This will lead to some modification
in the following definitions,
namely in the definition of what will be denoted by $V(G)$.

The free energy of a skeleton family
$S = (\gamma_j : j < k)$
will be made of two parts.
On the one hand the {\it surface free energy} of $S$
is simply defined by 
$$
    {\cal W}(S) = \sum_{j < k} {\cal W}(\gamma_j),
$$
with $W(\gamma_j)$ defined by Equation~\eqref{marina}.
Even if $\gamma_j$ is self-intersecting
and is not the boundary
of a simply connected domain,
so that the {\it external} normal can be ill-defined,
one can still define some normal 
with respect to an orientation of $\gamma_j$
and, since $\tau$ is $(\pi / 2)$-periodical,
there is no ambiguity for the resulting integral.

The volume free energy, on the other hand,
is related with the phase volume of $S$
introduced in \cite{DKS},
Section~2.10.
The {\it plus-components} of $S$
are the bounded connected components
of $\mathds{R}^2 \setminus \cup_{j < k} \gamma_j$
for which there is a continuous path
that connects their interior 
and the unique unbounded component  
of $\mathds{R}^2 \setminus \cup_{j < k} \gamma_j$
with an odd number of crossings of 
$\cup_{j < k} \gamma_j$.
The {\it phase volume} of $S$
is defined as their joint volume
and we denote it by $\hat V(S)$.
The plus-components of $G$ are defined 
in the same way we defined those of $S$
and we call $V(G)$ 
the total number of sites they enclose.
We define $\check V(S)$ 
as the number of sites
in the plus-components of $S$ 
that are at distance larger than $1 / h^r$
from $\cup_{j < k} \gamma_j$.
The {\it volume free energy} of $S$
is the product $-h m^*_\beta \check V(S)$.

Following \cite{SS} there is a constant $C > 0$
such that
$$
    \max\Bigl\{V(G), \hat V(S)\Bigr\} - C {\cal W}(S){1 \over h^{2r}}
    \leq \check V(S)
    \leq \min\Bigl\{V(G), \hat V(S)\Bigr\} 
$$
and
$$
    \Bigl|
        V(G) - \hat V(S)
    \Bigr| \leq C {\cal W}(S) {1 \over h^{2r}}
    \,.
$$
We will denote by $\bigl\{S_{\Lambda, -, h} = S\bigr\}$
the set of configurations that are associated
with the skeleton family $S$ and by
$$
    {\cal I} = \bigl\{S_{\Lambda, -, h} = \emptyset\bigr\}
$$
the set of configuration with invertebrate contours only.
These are similar to the configuration sets
${\cal S}_S^{h, \scriptscriptstyle sk}$
and ${\cal S}_\emptyset^{h, \scriptscriptstyle sk}$
in \cite{SS},
which are associated with external contours only.
Following the proof of Lemma~2.3.6 of \cite{SS},
we have
\begin{lmm}\label{enea}
    Given $\epsilon > 0$, if $h$ is small enough and $\Lambda$
    is a simply-connected domain contained in $\Lambda_h$,
    then, for any skeleton familly $S$,
    $$
        \mu_{\Lambda, -, h}\left(
            S_{\Lambda, -, h} = S
        \right) 
        \leq \mu_{\Lambda, -, h}\left(
            {\cal I}
        \right) \exp\left\{
            -\beta \Bigl(
                (1 - \epsilon) W(S) - (1 + \epsilon) h m^*_\beta \check V(S)
            \Bigr)
        \right\}.
    $$
\end{lmm}

This result relies on Pfister's low temperature estimate
for zero magnetic field (\cite{Pfi}, Lemma~10.1),
that was extended in \cite{Iof2} up to critical temperature,
and which uses a duality argument
that holds for simply connected domains only.
This is where the simple connectivity of $\Lambda$ matters.

To make the volume free energy appear,
Schonmann and Shlosman control the derivative 
with respect to $h$ of the ratio between
$\mu_{\Lambda, - , h}(S_{\Lambda, -, h} = S)$
and $\mu_{\Lambda, -, h}({\cal I})$
and they use in particular the fact that,
at any subcritical temperature,
there is a positive constant $C$ such that,
for all $h \geq 0$, $\Lambda \subset \mathds{Z}^2$
and $x, y \in \mathds{Z}^2$,
\begin{equation}\label{etienne}
    \mu_{\Lambda, +, h}\left(
        x \mathop{\longleftrightarrow}^{-*} y
    \right)
    \leq \mu_{h}\left(
        x \mathop{\longleftrightarrow}^{-*} y
    \right)
    \leq \mu_{+}\left(
        x \mathop{\longleftrightarrow}^{-*} y
    \right)
    \leq \exp\bigl\{
        -C \|x - y\|_\infty
    \bigr\},
\end{equation}
where the star percolation event
$\displaystyle \bigl\{x \mathop{\longleftrightarrow}^{-*} y\bigr\}$
is the set of configurations $\sigma$
in $\Omega_{\mathds{Z}^2}$
for which there is a sequence of sites
$x = z_0$, $z_1$,~\dots, $z_k = y$
such that $\|z_j - z_{j + 1}\|_\infty = 1$
and $\sigma(z_j) = \sigma(y) = -1$
for all $j < k$.
The first two inequalities are a consequence
of FKG inequality and the last one
is Theorem~1 in \cite{CCS}.

We then get upper bounds on events of type
$\bigl\{
    {\cal W}(S_{\Lambda, -, h}) \geq D / h^u,\,
    h m^*_\beta\check V(S_{\Lambda, -, h}) \leq E / h^v
\bigr\}$
for $u, v \geq r$.
\begin{lmm}\label{thomas}
    Given $\epsilon > 0$, $D_0 > 0$ and $E_0 > 0$,
    if $h$ is small enough and $\Lambda$
    is a simply-connected domain contained in $\Lambda_h$,
    then, for any $D \geq D_0$, $E \geq E_0$, $F \geq 0$ and $u, v \geq r$,
    it holds
    $$
        \mu_{\Lambda, -, h}\left(
            {\cal W}(S_{\Lambda, -, h}) \geq {D \over h^u} + (1 +  \epsilon) F,\,
            h m^*_\beta \check V(S_{\Lambda, -, h}) = F
        \right)
        \leq \mu_{\Lambda, -, h}\left(
            {\cal I}
        \right) \exp\left\{
            - \beta (1 - \epsilon) {D \over h^u}
        \right\}
    $$
    and
    $$
        \mu_{\Lambda, -, h}\left(
            {\cal W}(S_{\Lambda, -, h}) \geq {D \over h^u}\,,\,
             h m^*_\beta \check V(S_{\Lambda, -, h}) \leq {E \over h^v}
        \right)
        \leq \mu_{\Lambda, -, h}\left(
            {\cal I}
        \right) \exp\left\{
            -\beta \left(
                (1 - \epsilon) {D \over h^u} - (1 + \epsilon) {E \over h^v}
            \right)
        \right\}.
    $$
\end{lmm}
\par\noindent{\bf Proof:} 
This is similar to the proof of Lemma~2.3.7 in \cite{SS}.
For any $D \geq D_0$ and $F \geq 0$ it holds
\begin{align*}
    &\mu_{\Lambda, -, h}\left(
        {\cal W}(S_{\Lambda, -, h}) \geq {D \over h^u} + (1 + 4\epsilon) F,\,
        h m^*_\beta \check V(S_{\Lambda, -, h}) = F
    \right) \\
    &\qquad\leq \sum_{k \geq 0} \mu_{\Lambda, -, h}\left(
        {\cal W}(S_{\Lambda, -, h}) \in (1 + 4\epsilon) F  + \left[
            {(1 + k) D \over h^u}, {(1 + k + 1) D \over h^u}
        \right],\,
        h m^*_\beta \check V(S_{\Lambda, -, h}) = F 
    \right).
\end{align*}

For $h$ small enough,
the number of possible skeleton families $S$
such that
$$
    {\cal W}(S) \leq (1 + 4\epsilon) F + {(2 + k) D \over h^u}
$$
is less than (recall that $\beta$ is a fixed parameter)
$$
    \left(
        3 {B_{\max}^2 \over h^2}
    \right)^{
        {(1 + 4\epsilon) F + (2 + k) D / h^u \over \tau(0) / 12 } h^r
    }
    \leq \exp\left\{
        \beta \epsilon \left(
            (1 + 4\epsilon) F + {(2 + k) D  \over h^u}
        \right)
    \right\}.
$$
Indeed, since 
$$
    \tau(0) = \min_{\theta < 2 \pi} \tau(\theta),
$$
the second skeleton property implies
that, with $N$ the total number of vertices
of a skeleton family $S$,
$$
    {\cal W}(S) \geq N {1 \over 12 h^r} \tau(0),
$$
which gives an upper bound on $N$.
Together with the fact that these vertices
have to be in $\Lambda_h$,
of volume $(B_{\max} / h)^2$ at most
and that each of them can be a first,
last or intermediate vertex of a given skeleton,
this gives the stated upper bound.

Lemma~\ref{enea} implies then,
for any $\epsilon < 1 / 8$
and $h$ smaller
than some $h_0$ that depends on $\epsilon$, $D_0$
and $\beta$ only,
\begin{align*}
    &\mu_{\Lambda, -, h}\left(
        {\cal W}(S_{\Lambda, -, h}) \geq {D \over h^u} + (1 + 4\epsilon) F,\,
        h m^*_\beta \check V(S_{\Lambda, -, h}) = F
    \right) \\
    &\quad\leq \mu_{\Lambda, -, h}\left(
            {\cal I}
    \right) \sum_{k \geq 0} \exp\left\{
        \beta \left(
            \epsilon (1 + 4\epsilon) F + \epsilon {(2 + k) D \over h^u}
            - (1 - \epsilon) \left(
                {(1 + k) D \over h^u}
                + (1 + 4\epsilon) F
            \right)
            + (1 + \epsilon) F 
        \right)
    \right\} \\
    &\quad = \mu_{\Lambda, -, h}\left(
            {\cal I}
    \right) \sum_{k \geq 0} 
    \exp\left\{
        -\beta \left(
            \bigl[(1 - 3\epsilon) + (1 - 2\epsilon) k\bigr] {D \over h^u}
            + \epsilon (1 - 8\epsilon) F
        \right)
    \right\}\\
    &\quad\leq \mu_{\Lambda, -, h}\left(
            {\cal I}
    \right) C 
    \exp\left\{
        -\beta (1 - 3\epsilon) {D \over h^u}
    \right\}
\end{align*}
for some constant $C$ 
that depends on $\epsilon$ and $D_0$ only.
This implies the first desired inequality
with $4 \epsilon$ in place of $\epsilon$.
    
For $\epsilon < 1 / 2$, any $D \geq D_0$, $E \geq E_0$ and $h$ small enough
it holds in the same way
\begin{align*}
    &\mu_{\Lambda, -, h}\left(
        {\cal W}(S_{\Lambda, -, h}) \geq {D \over h^u},\, h m^*_\beta \check V(S_{\Lambda, -, h}) \leq {E \over h^v}
    \right) \\
    &\qquad\leq \sum_{k \geq 0} \sum_{j \leq {E / m^* \over h^{1 + v}}}
    \mu_{\Lambda, -, h}\left(
            {\cal W}(S_{\Lambda, -, h}) \in \left[(1 + k) {D \over h^u}, (1 + k + 1) {D \over h^u}\right],\,
            \check V(S_{\Lambda, -, h}) = j
    \right)\\
    &\qquad\leq \mu_{\Lambda, -, h}\left(
            {\cal I}
    \right) \sum_{k \geq 0}  {E / m^*_\beta \over h^{1 + v}}
    \exp\left\{
        - \beta \left(
            \bigl[(1 - 3\epsilon) + (1 - 2\epsilon) k\bigr] {D \over h^u}
            - (1 + \epsilon) {E \over h^v} 
        \right)
    \right\} \\
    &\qquad\leq \mu_{\Lambda, -, h}\left(
            {\cal I}
    \right) C 
    \exp\left\{
        -\beta \left(
            (1 - 3\epsilon) {D \over h^u}
            - (1 + \epsilon) {E \over h^v}
        \right)
    \right\}
\end{align*}
for some constant $C$ that depends on $\epsilon$,
$D_0$ and $E_0$ only.
The thesis follows.\qed

\medskip\par
For $\sigma$ in $\Omega_{\Lambda, -}$
we will also consider the family
$G^{\rm ext}_{\Lambda, -, h}(\sigma) = (\Gamma_j : j < k)$
of the {\it external} vertebrate contours of $\sigma$
as well as the family
$S^{\rm ext}_{\Lambda, -, h}(\sigma) = (\gamma_j : j < k)$
of their associated skeletons.
We will denote by
$$
    |G^{\rm ext}_{\Lambda, -, h}(\sigma)|
    = |S^{\rm ext}_{\Lambda, -, h}(\sigma)|
$$
their number $k$.    
As a first application of the previous upper bounds
we have that,
conditionally to $V(G^{\rm ext}_{\Lambda_h, -, h}) \leq (B_+ / h)^2$
and for $B_+$ small enough
---say $B_+ \leq 3 B_c / 2$ and recall that $B_+$
as to be thought close to $B_c$---
typical configurations drawn from $\mu_{\Lambda_h, -, h}$
are made of invertebrate contours only,
i.e., are in ${\cal I}$.
More precisely 
\begin{lmm}\label{jerome}
There is $\delta > 0$ such that, if $h$ is small enough
and $B \leq 3 B_c / 2$, then,
for all $k \geq 0$ it holds
$$
    \mu_{\Lambda_h, - , h}\Bigl(
        |S^{\rm ext}_{\Lambda_h, -, h}| = k,
        V(G^{\rm ext}_{\Lambda_h, -, h}) \leq (B / h)^2
    \Bigr) 
    \leq \mu_{\Lambda_h, -, h}({\cal I}) \exp\left\{
        - \delta k / h^b
    \right\}.
$$
In particular, for $B_+ \leq 3 B_c / 2$ and $h$ small enough,
it holds
\begin{align*}
    \mu_{\Lambda_h, -, h}\bigl(
        {\cal I}^c
        \bigm| {\cal R}
    \bigl)
    &\leq {
        \sum_{k \geq 1} \mu_{\Lambda_h, - , h}\Bigl(
            |S^{\rm ext}_{\Lambda_h, -, h}| = k,
            V(G^{\rm ext}_{\Lambda_h, -, h}) \leq (B_+ / h)^2
        \Bigr) 
    \over
        \mu_{\Lambda_h, -, h}\bigl({\cal I}\bigr)
    }\\
    &\leq {
        \mu_{\Lambda_h, -, h} \bigl({\cal I}\bigr) 2\exp\left\{
            -\delta / h^b
        \right\}
    \over
        \mu_{\Lambda_h, -, h} \bigl({\cal I}\bigr)
    } = 2 \exp\left\{
            -\delta / h^b
    \right\}.
\end{align*}
\end{lmm}
\par\noindent
{\bf Proof:}
We will apply the first inequality of the previous lemma
with $\epsilon = 1 / 8$.
To this end we will give a lower bound on 
$$
    {\cal W}(S_{\Lambda_h, -, h})
    - (1 + \epsilon) h m^*_\beta \check V(S_{\Lambda_h, -, h})
    \geq {\cal W}(S^{\rm ext}_{\Lambda_h, -, h})
    - (1 + \epsilon) h m^*_\beta \check V(S^{\rm ext}_{\Lambda_h, -, h})
$$
provided that $|S^{\rm ext}_{\Lambda, -, h}| = k$
and $V(G^{\rm ext}_{\Lambda_h, -, h}) \leq (B / h)^2$.
If 
$$
    G^{\rm ext}_{\Lambda_h, -, h} = (\Gamma_j : j < k)
    \quad\hbox{and}\quad
    S^{\rm ext}_{\Lambda_h, -, h} = (\gamma_j : j < k),
$$
we also have
$$
    {\cal W}(S^{\rm ext}_{\Lambda_h, -, h})
    - (1 + \epsilon) h m^*_\beta \check V(S^{\rm ext}_{\Lambda_h, -, h})
    \geq \sum_{j < k} {\cal W}(\gamma_j) - (1 + \epsilon) h m^*_\beta \check V(\gamma_j)
$$
with $-h m^*_\beta \check V(\gamma_j)$ the volume free energy
of the single skeleton $\gamma_j$.
To give a lower bound on each term of this sum,
we recall that there is $C > 0$ such that,
with $V(\Gamma_j)$ the number of sites enclosed in $\Gamma_j$,
it holds
$$ 
    V(\Gamma_j) - C {\cal W}(\gamma_j) / h^{2r}
    \leq \check V(\gamma_j)
    \leq V(\Gamma_j)
$$
and we separate two cases.

If
$$
    C {\cal W}(\gamma_j) / h^{2r} \geq V(\Gamma_j) / 2,
$$
then, since $V(\Gamma_j) \geq 1 / h^{2b}$,
$$
    {\cal W}(\gamma_j)
    - (1 + \epsilon) h m^*_\beta \check V(\gamma_j)
    \geq \left[
        {h^{2r} \over 2C} - (1 + \epsilon) h m^*_\beta
    \right] V(\Gamma_j)
    \geq {2\delta \over \beta h^b}
$$
for $h$ small enough
and some positive $\delta$
that depends only on $\epsilon$, $C$ and $\beta$.
If instead
$$
    C {\cal W}(\gamma_j) / h^{2r} \leq V(\Gamma_j) / 2,
$$
then we have on the one hand 
\begin{equation}\label{valerie}
    {1 \over 2 h^{2b}}
    \leq {1 \over 2} V(\Gamma_j)
    \leq \check V(\gamma_j)
    \leq V(\Gamma_j)
    \leq \left(
        {3 B_c \over 2h}
    \right)^2,
\end{equation}
and on the other hand, using the isoperimetric property 
of the Wulff shape,
$$
    {\cal W}(\gamma_j)
    - (1 + \epsilon) h m^*_\beta \check V(\gamma_j)
    \geq w_\beta \sqrt{\check V(\gamma_j)}
    - (1 + \epsilon) h m^*_\beta \check V(\gamma_j).
$$
This lower bound is concave in $\check V(\gamma_j)$.
From \eqref{valerie} we need then 
to evaluate it in $1 / (2 h^{2b})$ 
and $(3 B_c)^2 / (2h)^2$ to find its minimum value.
Since, for some $\delta' \leq \delta$ and $h$ small enough it holds
$$
    {w_\beta  \over \sqrt 2 h^b}
    - (1 + \epsilon) {m^*_\beta \over 2} h^{1 - 2b} 
    \geq {2\delta' \over \beta h^b}
$$
and 
$$
    w_\beta  {3 B_c / 2 \over h}
    - (1 + \epsilon) m^*_\beta {(3 B_c / 2)^2 \over h} 
    \geq {2\delta' \over \beta h^b} \,,
$$
this leads to
$$
    {\cal W}(S_{\Lambda_h, -, h})
    - (1 + \epsilon) h m^*_\beta \check V(S_{\Lambda_h, -, h})
    \geq {2k\delta' \over \beta h^b}\,.
$$
We then get the desired estimate
by applying Lemma~\ref{thomas} and 
summing on all the possible values 
of the integer
$$
    \check V(S_{\Lambda, -, h})
    = F / (h m^*_\beta)
    < 2 \left(
        B_{\max} \over h
    \right)^2.
$$
\qed

\medskip\par
We will also need lower bounds based on \cite{Iof1}.
For $B > 0$ and $\delta > 0$,
let us denote by 
$E^h_{B, \delta}$ the event that
there is an external contour which surrounds $(1 - \delta) B W / h$
and is contained in $(1 + \delta) B W / h$,
and that moreover this is the only external vertebrate contour.
With this notation and recalling Equation~\eqref{irene}
from page \pageref{irene}, Lemma~3.4.3 in \cite{SS} gives
\begin{lmm}\label{daniele}
    There are $C > 0$ and $o_h(1)$,
    a vanishing function of $h$ when $h$ goes to zero,
    such that, for all $B > 0$, $\delta > 0$
    and all simply-connected $\Lambda \subset \Lambda_h$
    that contains $(1 + \delta) B W / h$, it holds
    $$
        \mu_{\Lambda, - , h}(E^h_{B, \delta})
        \geq \mu_{\Lambda, -, h}({\cal I}) C \exp\left\{
            -\beta (1 + o_h(1)) {\phi(B) \over h}
        \right\}.
    $$
\end{lmm}

This makes possible to give lower bounds
on similar events for non simply-connected
``Wulff shaped annular domains'' of the form
$$
    A(B_1, B_2) = \left(
        {B_2 \over h} W  \setminus {B_1 \over h} W 
    \right) \cap \mathds{Z}^2
$$
with $0 \leq B_1 < B_2$.
(In the case $B_1 = 0$ this ``annular domain''
is simply a Wulff shaped box.)
For $\eta$, $\eta_1$ and $\eta_2$ in $\Omega_{\mathds{Z}^2}$
such that $\eta$ coincides with $\eta_1$ in $B_1 W / h$
and with $\eta_2$ outside $B_2 W / h$,
we will write $\mu_{A, (\eta_1, \eta_2), h}$
for $\mu_{A, \eta, h}$ with $A = A(B_1, B_2)$.
Given $\delta > 0$ we also define $B_{1, \delta} > B_1$
and $B_{2, \delta} < B_2$ by the equations
$$
    (1 - \delta) B_{1, \delta} = B_1
    \qquad{\rm and}\qquad
    (1 + \delta) B_{2, \delta} = B_2
$$
and we call $\tilde E^h_{B, \delta}$
the subset of $E^h_{B, \delta}$
for which there is no vertebrate contour
distinct from the external contour 
which surrounds $(1 - \delta) BW / h$
and is contained in $(1 + \delta) BW / h$.
\begin{lmm}\label{peter}
    Given $\epsilon > 0$, if $h$ is small enough,
    then, for all $0 \leq B_1 < B_2 \leq B_{\max}$ 
    and $\delta$ such that $B_{1, \delta} < B_{2, \delta}$,
    it holds, with $A = A(B_1, B_2)$,
    \begin{equation}\label{adrien_1}
        \mu_{A, (+, -), h}\left(
            \tilde E^h_{B_{2, \delta}, \delta}
        \right) \geq \exp\left\{
            - {\beta  \over h} \Bigl(
                \epsilon 
                + \left[
                    \phi(B_2) - \phi(B_1)
                \right]_+
            \Bigr)
        \right\},
    \end{equation}
    \begin{equation}\label{adrien_2}
        \mu_{A, (+, -), h}\left(
            E^h_{B_{1, \delta}, \delta}
        \right) \geq \exp\left\{
            - {\beta  \over h} \Bigl(
                \epsilon 
                + \left[
                    \phi(B_1) - \phi(B_2)
                \right]_+
            \Bigr)
        \right\},
    \end{equation}
    \begin{equation}\label{adrien_3}
        \mu_{A, (-, +), h}\Bigl(
            E^h_{B_{1, \delta}, \delta}
        \Bigr) \geq \exp\left\{
            - {\beta  \over h} \epsilon 
        \right\}
   \end{equation}
    and, if $B_2 - B_1 < 2 B_c$,
    \begin{equation}\label{adrien_4}
        \mu_{A, (-, -), h}\bigl(
            {\cal I}
        \bigr) \geq \exp\left\{
            - {\beta  \over h} \epsilon 
        \right\}.
   \end{equation}
\end{lmm}

\par\noindent
{\bf Proof:}
Most of this is already contained in Lemma~3.5.1
of \cite{SS},
which gives stronger lower bounds 
on similar events,
and its proof,
which works by conditioning and stochastic domination.
We will proceed in the same way.
Let us first prove \eqref{adrien_1}.
Our event $\tilde E^h_{B_{2, \delta}, \delta}$
is the intersection of the events
\begin{description}
\item[$E_0$:] there is a contour $\Gamma$ 
    that separates interior plus spins from exterior minus spins,
    that surrounds $(1 - \delta) B_{2, \delta} W / h$
    and that is contained in $B_2 W / h$,
\item[$E_1$:] such a contour $\Gamma$ does not enclose
    any vertebrate contour
\end{description}
and
\begin{description}
\item[$E_2$:] there is no vertebrate contour
    outside such a contour $\Gamma$,
\end{description}
the first two of which are increasing events.
With
$$
    \Lambda_2 = {B_2 \over h} W \cap \mathds{Z}^2,
$$
DLR equations imply
\begin{align*}
    \mu_{A, (+, -), h}(E_0 \cap E_1 \cap E_2)
    &= \mu_{A, (+, -), h}(E_0 \cap E_1) \times \mu_{A, (+, -), h}\left(
        E_2 \bigm| E_0 \cap E_1 
    \right)\\
    &= \mu_{A, (+, -), h}(E_0 \cap E_1) \times \mu_{\Lambda_2, -, h}\left(
        E_2 \bigm| E_0 \cap E_1 
    \right)
\end{align*}
and we will use stochastic domination 
for giving a lower bound 
of the first factor.
Let us denote by $\tilde \Lambda_1$
the set of sites in  
$$
    \Lambda_1 = {B_1 \over h} W \cap \mathds{Z}^2
$$
that are at distance $2 / h^{2b}$ from its boundary,
and by $F$ the event that there is a contour $\Gamma'$
which separates interior plus spins from exterior minus spins,
surrounds $\tilde \Lambda_1$
and does not enclose any vertebrate contour 
that encloses some site in $\tilde \Lambda_1$.
By conditioning on the invertebrate contours
enclosed in $\Gamma'$ and enclosing some site in $\tilde \Lambda_1$,
FKG inequality gives
$$
    \mu_{A, (-, +), h}\bigl(E_0 \cap E_1\bigr)
    \geq \mu_{\Lambda_2, -, h}\bigl(E_0 \cap E_1 \bigm| F\bigr)
     = {\mu_{\Lambda_2, -, h}\bigl(E_0 \cap E_1\bigr) \over \mu_{\Lambda_2, -, h}(F)}.
$$
Together with the previous equality we then have
$$
    \mu_{A, (+, -), h}(E_0 \cap E_1 \cap E_2) 
    \geq {
        \mu_{\Lambda_2, -, h}(E_0 \cap E_1 \cap E_2)
    \over
        \mu_{\Lambda_2, -, h}(F)
    }\,.
$$
To get a lower bound on the numerator we use Lemma~\ref{daniele}
and Estimate~\eqref{etienne} from page \pageref{etienne}.
We observe that $E_0 \cap E_2 = E^h_{B_{2, \delta}, \delta}$
and that, conditionally to $E^h_{B_{2, \delta}, \delta}$,
a star percolation event involving some sites $x$ and $y$
at distance of order $1 / h^b$ has to occur 
if $E_1$ does not.
Since $\phi$ is bounded from above,
we obtain a constant $C > 0$
such that for $h$ and $\delta$ small enough,
$$
    \mu_{\Lambda_2, -, h}(E_0 \cap E_1 \cap E_2)
    \geq \mu_{\Lambda_2, -, h}({\cal I}) C \exp\left\{
        - {\beta \over h} \bigl(\phi(B_2) + \epsilon / 2\bigr)
    \right\}.
$$
To get an upper bound on the denominator
we observe that $F$ implies, for $h$ small enough,
that $\check V(S_{\Lambda_2, -, h})$
lies between $(1 - \epsilon) (B_1 / h)^2$ and $(B_2 / h)^2$
so that the minimal free energy cost is 
or order
$(\phi(B_1) \wedge \phi(B_2)) / h$.
Using Lemma~\ref{thomas},
we get, for $h$ small enough
$$
    \mu_{\Lambda_2, -, h}(F)
    \leq \mu_{\Lambda_2, -, h}({\cal I}) \exp\left\{
        -{\beta \over h} \bigl(\phi(B_2) \wedge \phi(B_1) - \epsilon / 2)
    \right\}.
$$
This gives the desired estimate.

Inequality~\eqref{adrien_2} is proved in the same way:
it holds with $\tilde E^h_{B_{1, \delta}, \delta}$
in place of $E^h_{B_{1, \delta}, \delta}$,
but we will only need an estimate for this larger event.
Inequality~\eqref{adrien_3} 
is then a consequence of \eqref{adrien_2}:
the boundary conditions are exchanged
and the positive magnetic helps in such a way
that there is no size-dependent free energy cost anymore.
We refer to the last page of \cite{SS} 
for more details.

We finally prove \eqref{adrien_4}.
This is the only place where we will make
use of the notion of free energy
associated with non-external vertebrate contours.
Let us now denote by $E$ the event
that there is no vertebrate contour in $A$
and by $F$ the event
that there is a contour $\Gamma$
which separates external minus spins
from internal plus spins,
is enclosed in $\Lambda_1$
and encloses $(1 - \delta) \Lambda_1 / (1 + \delta)$.
Since ${\cal I}$ is a decreasing event
it holds
$$
    \mu_{A, (-, -), h}({\cal I})
    \geq \mu_{\Lambda_2, -, h}\left(
        E \bigm| F
    \right)
    = {
        \mu_{\Lambda_2, -, h}\left(
            E^h_{B_1 / (1 + \delta), \delta}
        \right)
    \over 
        \mu_{\Lambda_2, -, h}(F)
    }
$$
and, using Lemma~\ref{daniele},
we only need to prove that, 
for $h$ and $\delta$ small enough,
$$
    \mu_{\Lambda_2, -, h}(F) 
    \leq \mu_{\Lambda_2, -, h}({\cal I}) \exp\left\{
        - {\beta \over h} \bigl(\phi(B_1) - \epsilon / 2\bigr)
    \right\}.
$$
In other words we need to show
that the free energy of the skeleton families
that are compatible with $F$
cannot macroscopically decrease
with respect to that of the skeleton families
that are compatible with $E^h_{B_1, \delta}$.
Like in the proof of Lemma~\ref{jerome}
we can estimate from below 
the free energy of the former
by the sum of the free energy 
of the single skeleton associated 
with the contour $\Gamma$,
and that of the skeleton family
associated with each plus-component
outside $\Gamma$.
Since the former is of order $\phi(B_1) / h$,
it is sufficient to check that the latter
can only have a positive contribution
provided that $B_2 - B_1 < 2 B_c$.
Let us denote by ${\cal W}(S)$, 
$-h m^*_\beta \check V(S)$ and $\hat V(S) \geq \check V(S)$
the surface free energy, the volume free energy
and the phase volume of such a skeleton
family associated with a single plus-component
of the whole contour family.
If this single plus-component is simply connected,
then, by using Lemma~\ref{kata} and
Proposition~\ref{william},
the associated free energy
has a lower bound of order
$$
    {w_\beta \over 2}\left(
        {\check V(S) \over (B_2 - B_1)/(2 h)} 
        + {B_2 - B_1 \over 2h}
    \right) - h m^*_\beta \check V(S)
    \geq h \check V(S) \left(
        {w_\beta \over B_2 - B_1} - m^*_\beta
    \right).
$$
If it is not simply connected
but does not enclose $\Gamma$,
we get a similar lower bound
on its associated free energy
by estimating it from below
with that of the single skeleton
associated with its outermost contour.
If instead it is not simply connected
and it encloses $\Gamma$,
then,
denoting by $(B / h)^2$ the number of sites
enclosed in its outermost contour
and taking into account the surface free energy contribution
of its innermost contour,
the total free energy of this skeleton family
has a lower bound of order
$$
    w_\beta {B \over h} + w_\beta {B_1 \over h}
    - h m^*_\beta \left(
        \left(
            {B \over h}
        \right)^2 - \left(
            {B_1 \over h}
        \right)^2
    \right)
    \geq {B + B_1 \over h}\left(
        w_\beta - m^*_\beta (B_2 - B_1)
    \right).
$$
Provided that
$$
    B_2 - B_1 < 2 B_c = {w_\beta \over m^*_\beta}
    \,,
$$
this gives in all cases a non-negative macroscopic contribution.
\qed
    
\subsection{Exit rates, local relaxation times and soft capacities}\label{zakia}
We will simply denote by $X$ the dynamics $X_{\Lambda_h, -, h}$
restricted to 
$$
    {\cal X} = {\cal R} \cup {\cal S},
$$
which is associated with the generator ${\cal L}$
defined by
$$
    ({\cal L}f)(\sigma)
    = \sum_{
        {
           \scriptstyle x \in \Lambda_h :
        \atop
           \scriptstyle \sigma^x \in {\cal X}
        }
    } w(\sigma, \sigma^x) \bigl[
        f(\sigma^x) - f(\sigma)
    \bigr],
    \qquad \sigma \in {\cal X},
    \qquad f: {\cal X} \rightarrow \mathds{R}.
$$
We will also denote by $\mu$ its reversible measure
$$
    \mu = \mu_{\Lambda_h, -, h}(\cdot\, | {\cal X}).
$$
and by ${\cal D}$ 
the associated Dirichlet form
defined by Equation~\eqref{baptiste}
of Section~\ref{nelly}.
Its spectral gap will be denoted $\gamma = \gamma_h$.
In this section we briefly recall some definitions 
from \cite{BGM} and explain how to use the results
of that paper to prove
an equivalent of Theorem~\ref{sandrine}
and Proposition~\ref{amandine}
for this restricted dynamics $X$.

We denote by ${\cal L}_{\cal R}$ the generator 
of the dynamics $X$ restricted to ${\cal R}$:
$$
    ({\cal L}_{\cal R}f)(\sigma)
    = \sum_{
        {
           \scriptstyle x \in \Lambda_h :
        \atop
           \scriptstyle \sigma^x \in {\cal R}
        }
    } w(\sigma, \sigma^x) \bigl[
        f(\sigma^x) - f(\sigma)
    \bigr],
    \qquad \sigma \in {\cal R},
    \qquad f: {\cal R} \rightarrow \mathds{R},
$$
and we will denote by $1 / \gamma_{\cal R}$
the relaxation time of this restricted dynamics.
We denote by $\mu_{\cal R}$ the {\it restricted ensemble}
$$
    \mu_{\cal R} = \mu(\cdot\, |{\cal R}),
$$
with respect to which ${\cal L}_{\cal R}$ is reversible,
and we set
$$
    \chi_{\cal R}
    = \max_{\sigma \in {\cal R}} 
    {1 \over \mu_{\cal R}(\sigma)}\,.
$$
We define in the same way ${\cal L}_{\cal S}$,
$1 / \gamma_{\cal S}$, $\mu_{\cal S}$ and $\chi_{\cal S}$.
We will refer to $1 / \gamma_{\cal R}$
and $1 / \gamma_{\cal S}$ as {\it local relaxation times}.

For any $\lambda \geq 0$
we denote by $\phi^*_{{\cal R}, \lambda_{\cal S}}$
the extinction rate from quasi-stationarity 
of the trace on ${\cal R}$ of our process $X$
killed at rate $\lambda$ in ${\cal S}$,
and we set 
$$
    \phi^*_{{\cal R} \backslash {\cal S}}
    = \lim_{\lambda \rightarrow \infty} \phi^*_{{\cal R}, \lambda_{\cal S}}.
$$
The precise meaning of each of these terms
is explained in Section~2.1 of \cite{BGM},
from which we will mainly need the upper bound of Lemma~2.3
\begin{equation}\label{sophie_1}
    \phi^*_{{\cal R}, \lambda_{\cal S}}
    \leq \phi^*_{{\cal R} \backslash {\cal S}}
    \leq \mu_{{\cal R} \backslash {\cal S}}\bigl(
        e^*_{{\cal R} \backslash {\cal S}}
    \bigr),
\end{equation}
with $\mu_{{\cal R} \backslash {\cal S}} = \mu(\cdot\,|{\cal R} \backslash {\cal S})$
and
$$
    e^*_{{\cal R} \backslash {\cal S}}(\sigma)
    = \sum_{
        {
            \scriptstyle x \in \Lambda_h :
        \atop
            \scriptstyle \sigma^x \in {\cal S} 
        } 
    } w(\sigma, \sigma^x),
    \qquad x \in {\cal R} \backslash {\cal S}.
$$
For any $\kappa \geq 0$ we define in the same way
$\phi^*_{{\cal S}, \kappa_{\cal R}}$,
then $\phi^*_{{\cal S} \backslash {\cal R}}$,
$\mu_{{\cal S} \backslash {\cal R}}$
and $e^*_{{\cal S} \backslash {\cal R}}$.
It also holds
\begin{equation}\label{sophie_2}
    \phi^*_{{\cal S}, \kappa_{\cal R}}
    \leq \phi^*_{{\cal S} \backslash {\cal R}}
    \leq \mu_{{\cal S} \backslash {\cal R}}\bigl(
        e^*_{{\cal S} \backslash {\cal R}}
    \bigr).
\end{equation}
We will refer to $\phi^*_{{\cal R} \backslash {\cal S}}$
and $\phi^*_{{\cal S} \backslash {\cal R}}$
as {\it exit rates} from ${\cal R} \backslash {\cal S}$
and ${\cal S} \backslash {\cal R}$.

From Section~2.3 in \cite{BGM},
Dirichlet's and Thomson's principle,
the $(\kappa, \lambda)$-capacity
$C_{\kappa}^{\lambda}({\cal R}, {\cal S})$
is the {\it soft capacity}
\begin{align}
    C_\kappa^\lambda({\cal R}, {\cal S})
    &= \min_{
        f: {\cal X} \rightarrow \mathds{R}
    } \left\{
        {\cal D}(f)
        + \kappa \sum_{\sigma \in {\cal R}} \mu(\sigma)\left(
            f(\sigma) - 1
        \right)^2
        + \lambda \sum_{\sigma \in {\cal S}} \mu(\sigma)\left(
            f(\sigma) - 0
        \right)^2
    \right\} \label{guillaume_1}\\
    &= \max_{
        \tilde \psi \in \tilde\Psi_1(\bar R, \breve S)
    } \tilde {\cal D}\left(
        \tilde \psi
    \right)^{-1} \label{guillaume_2}
\end{align}
where
$$
    \tilde {\cal D}(\tilde \psi)
    = {1 \over 2} \sum_{\sigma \in {\cal X}} \sum_{x \in \Lambda}{
        \tilde \psi (\sigma, \sigma^x)^2 
    \over
        \mu(\sigma) w(\sigma, \sigma^x)
    } + \sum_{\sigma \in {\cal R}} {
        \tilde \psi(\sigma, \bar \sigma)^2
    \over 
        \mu(\sigma) \kappa
    } + \sum_{\sigma \in {\cal S}} {
        \tilde \psi(\sigma, \breve \sigma)^2
    \over 
        \mu(\sigma) \lambda
    }
$$
stands for the energy dissipated 
by a flow $\tilde \psi$
in the set $\tilde \Psi_1(\bar R, \breve S)$
of all the unitary flows from $\bar R$
to $\breve S$ associated with a Markov process $\tilde X$
on the extended 
$$
    \tilde {\cal X}
    = {\cal X} \cup \bar {\cal R} \cup \breve {\cal S}
$$
that jumps from any $\sigma$ in ${\cal R}$ or ${\cal S}$
to $\bar \sigma$ in $\bar {\cal R}$
or $\breve \sigma$ in $\breve {\cal S}$
at rate $\kappa$ or $\lambda$.

We will prove in sections \ref{veronique}
and \ref{julia} that the following hypothesis $({\cal H})$
is in force:
\smallskip\par
\noindent {\bf Hypothesis $({\cal H})$:} {
    \it Given a small enough $\delta > 0$, 
    one can choose $B_+$ close enough to $B_c$
    so that, for all $h$ small enough,
    it holds 
    $$
        {1 \over \gamma_{\cal R}}
        \vee
        {1 \over \gamma_{\cal S}}
        \leq \exp\left\{
            {\delta \over h}
        \right\},
    $$
    $$
        {1 \over \phi^*_{{\cal S} \backslash {\cal R}}}
        \wedge
        {1 \over \phi^*_{{\cal R} \backslash {\cal S}}}
        \geq  \exp\left\{
            {\beta A - \delta \over h}
        \right\}
    $$
    and, with $\kappa = \kappa(h)$ and $\lambda = \lambda(h)$
    such that
    $$
        \lim_{h \rightarrow 0} \kappa(h) e^{\delta / h} 
        = \lim_{h \rightarrow 0} {e^{-(\beta A - \delta) / h} \over \kappa(h)}
        = \lim_{h \rightarrow 0} \lambda(h) e^{\delta / h} 
        = \lim_{h \rightarrow 0} {e^{-(\beta A - \delta) / h}  \over \lambda(h)}
        = 0,
    $$
    for all $\epsilon > 0$ and $h$ small enough
    $$
        \exp\left\{
            - {\beta A + \epsilon  \over h}
        \right\}
        \leq {C_\kappa^\lambda({\cal R}, {\cal S}) \over \mu({\cal R})}
        \leq \exp\left\{
            - {\beta A  - \epsilon \over h}
        \right\}.
    $$
}
\smallskip\par\noindent
This will imply an equivalent of Theorem~\ref{sandrine}
together with Proposition~\ref{amandine}
and Estimate~\eqref{bruno}
for the restricted process $X$.
\begin{lmm}\label{lucia}
    If hypothesis $({\cal H})$ is in force,
    then, for all small enough $\delta_0 > 0$,
    one can choose $B_+$ close enough to $B_c$
    such that with $\kappa = \lambda = e^{- \delta_0 / (2h)}$
    there is $h_0 > 0$
    for which
    the following holds for $X$ started
    from a probability measure $\nu$ 
    and any observable $f : \Omega_{\Lambda_h} \rightarrow \mathds{R}$.
    \begin{itemize}
    \item[i.] If $\nu = \mu_{\cal R}$, then,
        for all $t > 0$,
        \begin{equation}\label{fabienne_2}
            \lim_{h \rightarrow 0} \mathds{P}_\nu\left(
                \gamma T_{\lambda_{\cal S}} > t
            \right) = e^{-t},
        \end{equation}
        and it holds 
        \begin{equation}\label{gregory_2}
            \lim_{h \rightarrow 0} h \ln {1 \over \gamma}
            = \beta A.
        \end{equation}
        Also,
        \begin{equation}\label{giulio_2}
            \lim_{h \rightarrow 0} \mathds{P}_\nu\left(
                \theta < T_{\lambda_{\cal S}}
                \hbox{ and }
                \sup_{t < T_{\lambda_{\cal S}} - \theta} 
                \left|
                    A_\theta(t, f) - \mu_{\cal R}\bigl(f\bigr)
                \right| \leq \|f\|_\infty e^{- \delta_0 /(11h)}
            \right) = 1
        \end{equation}
        with
        \begin{equation}\label{cecile_2}
            \theta = \exp\left\{
                {\beta A \over 2h} + {\delta_0 \over h}
            \right\}.
        \end{equation}
    \item[ii.] For $h < h_0$ and whatever the starting measure $\nu$, it holds
        $$ 
            \Bigl|
                E_\nu\left[
                    f\Bigl(
                        X\bigl(
                            T_{\lambda_{\cal S}}
                        \bigr)
                    \Bigr)
                \right] - \mu(f)
            \Bigr| \leq \|f\|_\infty e^{- \delta_0 / (6h)}.
        $$
    \item[iii.] If $\nu$ is such that
        \begin{equation}\label{laura}
            \lim_{h \rightarrow 0} \mathds{P}_\nu\bigl(
                T_{\kappa_{\cal R}}
                < T_{\lambda_{\cal S}}
            \bigr) = 1,
        \end{equation}
        then \eqref{fabienne_2}--\eqref{cecile_2} are also in force.
        Also if, for $h$ smaller than some positive $h_1$,
        \begin{equation}\label{garance}
            \mathds{P}_\nu\bigl(
                T_{\kappa_{\cal R}} > T_{\lambda_{\cal S}}
            \bigr) \leq e^{-\delta_0 / h}
        \end{equation}
        then 
        \begin{equation}\label{piero}
            \bigl|
                \mathds{E}_\nu\left[
                    f(X(t))
                \right] - \mu_{\cal R}(f)
            \bigr|
            \leq \|f\|_\infty e^{-\delta_0 / (6h)}
        \end{equation}
        for all small enough $h$
        and all $t = e^{\beta a / h}$
        with $\delta_0 < \beta a < \beta A - \delta_0$.
    \end{itemize}
\end{lmm}

\noindent
{\bf Proof:}
Let $\delta_0 > 0$ be small enough
to have
\begin{equation}\label{dominique}
    {\beta A \over 8} - {\delta_0 \over 16}
    > {\delta_0 \over 9}
    \qquad{\rm and}\qquad
    \phi\left(
        {B_{\max} \over (1 + \delta_0)}
    \right) < - \delta_0,
\end{equation}
and choose $B_+$
as provided by hypothesis $({\cal H})$
with $\delta_0 / 4$ in place of $\delta$.
We use the results of \cite{BGM},
which are based on two hypothesis sets
---denoted there by $(H)$ and $(H')$---
both satisfied with this choice
of ${\cal R}$ and ${\cal S}$ associated with $B_+$.
Indeed, hypotheses $(H)$ require
\begin{itemize}
\item[a)] $\phi^*_{{\cal R} \backslash {\cal S}}$
    to be small with respect to $\gamma_{\cal R}$
    and $\gamma_{\cal S}$
    in our considered asymptotic regime $h \ll 1$;
\item[b)] $\phi^*_{{\cal S} \backslash {\cal R}}$
    to be small with respect to $\gamma_{\cal S}$;
\item[c)] $X_{\cal R}$, $X_{\cal S}$,
    $X_{{\cal R} \backslash {\cal S}}$
    and $X_{{\cal S} \backslash {\cal R}}$
    to be all irreducible;
\item[d)] $\mu({\cal S}) \geq \mu({\cal R})$;
\end{itemize}
$({\cal H})$ gives a quantitative of version of a) and b);
$X_{\cal R}$ and $X_{{\cal R} \backslash {\cal S}}$
(as well as, symmetrically,
$X_{\cal S}$ and $X_{{\cal S} \backslash {\cal R}}$)
are irreducible since, by flipping each plus spin,
one gets a path in ${\cal R}$ or ${\cal R} \backslash {\cal S}$
from any configuration $\sigma$ to the uniform minus configuration;
and, as a consequence of Lemma~\ref{jerome}
and Lemma~\ref{daniele} with
$$
    B = {B_{\max} \over (1 + \delta_0)}
$$
we have, for all small enough $h$,
\begin{equation}\label{jean-jacques}
    {\mu({\cal R}) \over \mu({\cal S})}
    \leq {
        2 \mu({\cal I})
    \over
        \mu(E^h_{B, \delta_0})
    }
    \leq \exp\left\{
        - {\beta \over 2} \left|
             \phi\left(
                {B_{\max} \over (1 + \delta_0)}
            \right)
        \right|
    \right\}
    \leq \exp\left\{
        - {\delta_0 \over 2h}
    \right\}.
\end{equation}
Hypotheses $(H')$ require in addition 
$\phi^*_{{\cal R} \backslash {\cal S}}$
to be small with respect to $\gamma_{\cal R} / \ln \chi_{\cal R}$
and $\phi^*_{{\cal S} \backslash {\cal R}}$
to be small with respect to $\gamma_{\cal S} / \ln \chi_{\cal S}$,
which is also implied by $({\cal H})$
since there is a positive constant $C$
such that 
$$
    \ln\chi_{\cal R} \vee \ln\chi_{\cal S} 
    \leq C \left(
        {B_{\max} \over h}
    \right)^2.
$$
These hypothesis sets being satisfied,
setting $\kappa = \lambda = e^{-\delta_0 /(2h)}$,
$\kappa$ is large with respect to $\phi^*_{{\cal R} \backslash {\cal S}}$
and small with respect to $\gamma_{\cal R} / \ln \chi_{\cal R}$,
just as 
$\lambda$ is large with respect to $\phi^*_{{\cal S} \backslash {\cal R}}$
and $\phi^*_{{\cal R} \backslash {\cal S}}$
and small with respect to $\gamma_{\cal S} / \ln \chi_{\cal S}$.

Proposition~2.8 of \cite{BGM}, with $\lambda$ and ${\cal S}$
in place of $\kappa$ and ${\cal R}$, 
gives then, whatever the starting distribution $\nu$,
that the total variation distance
between $\mu_{\cal S}$ and the law of $X(T_{\lambda_{\cal S}})$
is smaller than $e^{-\delta_0 / (5h)}$ for $h$ small enough.
Since, as a consequence of \eqref{jean-jacques},
so is that between $\mu$ and $\mu_{\cal S}$,
this gives ii.  

Equations (15) and (16) and Proposition~2.8 of \cite{BGM} 
also give that $\phi^*_{{\cal R}, \lambda_{\cal S}} T_{\lambda_{\cal S}}$
converges in law to an exponential random variable
or parameter~1 as soon as \eqref{laura} is ensured.
Since \eqref{jean-jacques} implies that
$\mu({\cal S})$ goes to one
when $h$ goes to zero,
Theorem~1 of \cite{BGM} 
says that the ratios
$\phi^*_{{\cal R}, \lambda_{\cal S}} / \gamma$ and
$\phi^*_{{\cal R}, \lambda_{\cal S}} \mu({\cal R})/ C_\kappa^\lambda({\cal R}, {\cal S})$
go to one when $h$ goes to zero.
Then, provided~\eqref{laura},
$\gamma T_{\lambda_{\cal S}}$
converges in law to an exponential random variable
of parameter~1
---this is \eqref{fabienne_2}--- 
and \eqref{gregory_2} is implied by $({\cal H})$.

As far as the case $\nu = \mu_{\cal R}$
is concerned, we simply have to prove 
that \eqref{laura} is in force
to prove \eqref{fabienne_2}.
With
$$
    V_\kappa^\lambda(x) = \mathds{P}_x\left(
        T_{\kappa_{\cal R}} < T_{\lambda_{\cal S}}
    \right),
    \qquad x \in {\cal X}, 
$$
we have 
$$
    \mathds{P}_{\mu_{\cal R}}\left(
        T_{\kappa_{\cal R}} < T_{\lambda_{\cal S}}
    \right)
    = E_{\mu_{\cal R}}\left[
        V_\kappa^\lambda|_{\cal R}
    \right]
$$
and Lemma~3.2 in \cite{BGM} says that
that the latter goes to one
when $h$ goes to zero.
Conditions~\eqref{dominique} also
imply that we can choose $\eta = e^{-\delta_0 / (9h)}$
in Proposition~5.1 of \cite{BGM}
which, together with 
$$
    \lim_{h \rightarrow 0} h \ln \phi^*_{{\cal R}, \lambda_{\cal S}}
    = \beta A,
$$
gives,
for some
$$
    \tilde \theta \leq \exp\left\{
        {\beta A \over 2h} + {3 \delta_0 \over 4h}
    \right\}
$$
and $h$ small enough,
\begin{equation}\label{marie-christine}
    \mathds{P}_{\mu_{\cal R}}\left(
        \tilde \theta < T_{\lambda_{\cal S}}
        \hbox{ and }
        \sup_{t < T_{\lambda_{\cal S}} - \tilde\theta} 
        \left|
            A_{\tilde \theta}(t, f) - \mu_{\cal R}\bigl(f\bigr)
        \right| \leq \|f\|_\infty e^{- \delta_0 /(10h)}
    \right) \geq 1 - e^{-\delta_0 / (10 h)}.
\end{equation}
Since we already now that, starting from $\mu_{\cal R}$,
$\phi^*_{{\cal R}, \lambda_{\cal S}} T_{\lambda_{\cal S}}$
converges in law towards an exponential random variable 
of parameter one,
this implies \eqref{giulio_2}--\eqref{cecile_2}.

Next, Theorem~3 of \cite{BGM}
says that there is a stopping time
$T^*$, with
$$
    \mathds{E}_\nu\left[
        T^*
    \right] \leq 2 e^{\delta_0 / (2h)},
$$
such that the total variation distance
between $\mu_{\cal R}$ and the law of $X(T^*)$
goes to zero as well as the probability
that $T_{\lambda_{\cal S}} < T^*$
when \eqref{laura} is in force.
The contribution to time averages
on time scale $\theta$ of the trajectories
of $X$ before time $T^*$ 
is then negligible and 
we get, from \eqref{marie-christine},
that \eqref{laura}
implies \eqref{giulio_2}--\eqref{cecile_2}.

It only remains to prove~\eqref{piero}
by assuming~\eqref{garance}
for $h$ small enough.
We use to this end optimal couplings
associated with total variation estimates
provided by \cite{BGM}
to bound the total variation distance
between the law of $X(t)$ and $\mu_{\cal R}$.
First, by Markov inequality,
\begin{equation}\label{alma}
    \mathds{P}_\nu\left(
        T^* \geq t
    \right)
    \leq {2 e^{\delta_0 / (2h)} \over e^{\beta a / h}}
    \leq 2 e^{- \delta_0 / (2h)},
\end{equation}
and, assuming $T^* < t$,
we consider four coupled process $X_0$,
$X_1$, $X_2$ and $X_3$ on the time interval
$[T^*, t]$ with the following marginals:
$X_0(s) = X(s)$ for all $s \in [T^*, t]$;
$X_1(T^*)$ is distributed according to $\mu_{\cal R}$
and $X_1$ evolves according to the restricted dynamics
in ${\cal X}$ with generator ${\cal L}$;
$X_2$ evolves according to the same dynamics in ${\cal X}$,
but $X_2(T^*)$ is distributed according to 
the quasi-stationary distribution
$\mu^*_{{\cal R} \backslash {\cal S}}$ introduced 
in Section~2.1 of \cite{BGM} and for which,
with $T_{\cal S}$ the hitting time of ${\cal S}$, 
\begin{equation}\label{qasim}
    \mathds{P}_{\mu^*_{{\cal R} \backslash {\cal S}}}\left(
        T_{\cal S} > s
    \right) = e^{-\phi^*_{{\cal R} \backslash {\cal S}}s}
\end{equation}
for all $s \geq 0$;
$X_3(T^*) = X_1(T^*)$, but $X_3$ evolves according
to the restricted dynamics in ${\cal R}$,
so that the law of $X_3(t)$ is $\mu_{\cal R}$.
Then, we simply have to couple these processes
in such a way that $X_0(t) = X_3(t)$ 
with large probability.
Since $X_1(T^*) = X_3(T^*)$, it suffices to this end 
to couple $X_0(T^*)$, $X_1(T^*)$ and $X_2(T^*)$
to make them coincide with large probability 
and use \eqref{qasim} to prove that they
will not exit ${\cal R}$ with large probability.
Indeed, conditionally to $\{T^* < t\}$,
$$
    \mathds{P}\bigl(
        \exists s < t,\,
        X_2(s) \not\in {\cal R}
    \bigr)
    \leq \mathds{P}\bigl(
        \exists s < t,\,
        X_2(s) \in {\cal S}
    \bigr)
    \leq  1 - e^{-\phi^*_{{\cal R} \backslash {\cal S}} t}
    \leq \phi^*_{{\cal R} \backslash {\cal S}} t
    \leq e^{- 3 \delta_0 / (4h)}.
$$
Conditionally to $\{T^* < t \}$ and Hypothesis $({\cal H})$,
Proposition~2.6 in \cite{BGM} says that, for $h$ small enough,
we can couple $X_2(T^*)$ and $X_1(T^*)$ in such a way that
$$
    \mathds{P}\bigl(
        X_2(T^*) \neq X_1(T^*)
    \bigr) \leq \exp\left\{
        - {1 \over h} \left(
            {\beta A \over 2} - {\delta_0 \over 2}
        \right)
    \right\}.
$$
From Theorem~3 in \cite{BGM}
and \eqref{alma} we can couple $X_1(T^*)$ and $X_0(T^*)$
in such a way that, for $h$ small enough,
$$
    \mathds{P}\bigl(
        X_1(T^*) \neq X_0(T^*)
    \bigr)
    \leq \mathds{P}\bigl(T^* \geq t\bigr) 
    + \mathds{P}\bigl(T^* \neq T^*_{\cal R}\bigr)
    + e^{-\delta_0 / (5h)}
    \leq 2 e^{-\delta / (2h)} + 3 e^{-\delta_0 / h} + e^{-\delta_0 / (5h)}.
$$
With such couplings we get
$$
    \mathds{P}\bigl(
        X_3(t) \neq X(t)
    \bigr)
    \leq {1 \over 2} e^{- \delta_0 / (6h)}
$$
for $h$ small enough,
and \eqref{piero} follows.
\qed
        
\medskip\par
Assuming hypothesis $({\cal H})$,
the proof of Theorem~\ref{sandrine}
and Proposition~\ref{amandine}
essentially reduces at this point
to show that,
starting from $\mu_{\Lambda_h, -, h}(\cdot\,|{\cal R})$
and with large probability,
the system does not leave
${\cal X} = {\cal R} \cup {\cal S}$
within a time of order $e^{\beta A / h}$.
We need then a lower bound on an exit time,
like are the lower bounds on the inverse exit rates
and the upper bound of the soft capacity 
in hypothesis $({\cal H})$.
Given the previous free energy estimates
and the non-convex Blashke's inequality,
these are standard estimates in the context
of metastability studies.
They boil down to static estimates
(recall in particular
\eqref{sophie_1} and \eqref{sophie_2})
and we will prove them
in the next section.
As far as the upper bounds
on the local relaxation times
and the lower bound on the soft capacity are concerned,
we will follow the strategy introduced in Section~\ref{nelly},
and inspired by the works of Sinclair and Martinelli,
to prove them in Section~\ref{julia}.

\section{Lower bounds for exit times}\label{veronique}
\subsection{Leaving \bm{${\cal X}$}}
Before stating and proving the main lemma of this section 
we note that, for any $a > 0$, 
$$
    x > 0 \mapsto {a^2 \over x} + x
$$
is a convex function that reaches its minimum $2a$ in $a$.
\begin{lmm}\label{alphonsine}
    Given $B_+ > B_c$ and $b < 1 / 4$, 
    setting $\eta > 0$ such that
    $$
        {B_c^2 \over B_+} + B_+ 
        = 2 B_c (1 + 2 \eta),
    $$ 
    it holds 
    $$
        \mu_{\Lambda_h, - , h}\Bigl(
            ({\cal R} \cup {\cal S})^c
        \Bigr) \leq \mu_{\Lambda_h, -, h}({\cal I}) \exp\left\{
            - {\beta  \over h} A (1 + \eta) 
        \right\}
    $$
    for $h$ small enough.
\end{lmm}

\par\noindent
{\bf Proof:} Consider $\sigma$ in 
$$
    {\cal X}^c = ({\cal R} \cup {\cal S})^c.
$$
Let us denote by $S$ the skeleton collection
associated with its vertebrate contours,
by $G^{\rm ext}$ the collection
of its external vertebrate contour,
by $S^{\rm ext}$ its associated skeleton collection,
and set $B > 0$ such that 
$\check V(S^{\rm ext}) = (B / h)^2$.

Let us first consider the case $B \geq B_c$.
Since $\sigma \not\in {\cal S}$,
the largest Wullf shape enclosed
by a contour $\Gamma$ of $G^{\rm ext}$
has a volume smaller than $(B_- / h)^2$.
Recall Equation~\eqref{charles} of page~\pageref{charles},
set $\rho_- = 2 B_- / w_\beta$,
call $\gamma$ the skeleton of $\Gamma$
and $B_2(0, r)$ the Euclidean ball 
of radius $r$ centered in the origin.
As a consequence of the third skeleton property,
the largest Wulff shape contained
in a bounded connected component of $\mathds{R}^2 \setminus \gamma$
is contained in a translate of 
$$
    W_{\rho_- / h} + B_2(0, 1 / h^r)
    \subset W_{\rho_- / h} + W_{1 / (\tau(0) h^r)}
    = {w_\beta \over 2}\left(
        {\rho_- \over h} + {1 \over \tau(0) h^r}
    \right)W
$$
with volume less than
$$
    \left(
        {B_- \over h}
        + {w_\beta \over 2 \tau(0) h^r}
    \right)^2
    \leq \left(
        \tilde B_- \over h
    \right)^2
$$
for any $\tilde B_- > B_-$ and $h$ small enough.
Let us take $\tilde B_-$ close enough to $B_-$
to have $\tilde B_- < B_c$
and 
$$
    {B_c^2 \over \tilde B_-} + \tilde B_-
    \geq \left(
        {B_c^2 \over B_-} + B_-
    \right){
        1 + 3 \eta / 2
    \over
        1 + 2 \eta
    }
    \,.
$$
Since
$$
    {B_c^2 \over B_-} + B_-
    =  {B_c^2 \over B_+} + B_+ + 2{(B_+ - B_c)^3 \over B_c B_-}
    \geq {B_c^2 \over B_+} + B_+ 
    = 2 B_c (1 + 2 \eta),
$$ 
this implies 
$$
    {B_c^2 \over \tilde B_-} + \tilde B_-
    \geq 2 B_c \left(
        1 + {3\eta \over 2}
    \right).
$$
For any positive and small enough $\epsilon$,
Proposition~\ref{william}
now implies, since $B \geq B_c > \tilde B_-$,
\begin{align*}
    {\cal W}(S) - (1 + \epsilon) h m^*_\beta \check V(S)
    &\geq {\cal W}(S^{\rm ext}) - (1 + \epsilon) h m^*_\beta \check V(S^{\rm ext})\\
    &\geq {w_\beta \over 2}\left(
        {(B / h)^2 \over \tilde B_- / h} + \tilde B_- / h
    \right)
    - (1 + \epsilon) h m^*_\beta (B / h)^2\\
    &= {w_\beta \over 2h} \left[
        B^2\left(
            {1 \over \tilde B_-} - {1 + \epsilon \over B_c}
        \right) + \tilde B_-
    \right]
    \geq {w_\beta \over 2h} \left[
        B_c^2\left(
            {1 \over \tilde B_-} - {1 + \epsilon \over B_c}
        \right) + \tilde B_-
    \right]\\
    &\geq {w_\beta \over 2h} \left[
        2 B_c \left(
            1 + {3 \eta \over 2}
        \right) - (1 + \epsilon) B_c
    \right]
    = {w_\beta B_c\over 2h} \left(
        1 + 3 \eta -\epsilon 
    \right)\\
    &\geq {A \over h}\left(
        1 + 2 \eta 
    \right).
\end{align*}
We are in shape to use Lemma~\ref{thomas},
but let us first consider the alternative case $B \leq B_c$.

If $B \leq B_c$, i.e.,
$\check V(S^{\rm ext}) \leq B_c^2  h^{-2}$,
and $V(G^{\rm ext}) \geq (3 B_c / 2)^2 h^{-2}$,
we also have a lower bound on the free energy.
Recalling, indeed, that there is a positive constant $C$
such that 
$$
    V(G^{\rm ext}) - C{\cal W}(S^{\rm ext})h^{-2r} \leq \check V(S^{ext})
$$
it follows that,
for any positive and small enough $\epsilon$,
\begin{align*}
    {\cal W}(S) - (1 + \epsilon) h m^*_\beta \check V(S)
    &\geq {\cal W}(S^{\rm ext}) - (1 + \epsilon) h m^*_\beta \check V(S^{\rm ext})\\
    &\geq {h^{2r - 2} \over C} \left(
        (3 B_c  / 2)^2 - B_c^2 
    \right)
    - (1 + \epsilon) m^*_\beta {B_c^2 \over h}\\
    &\geq {5 B_c^2 \over  4 C h^{7 / 4}} - {2 m^*_\beta B_c^2 \over h}\,,
\end{align*}
so that, for $h$ small enough,
$$
    {\cal W}(S) - (1 + \epsilon) h m^*_\beta \check V(S)
    \geq {A \over h}\left(
        1 + 2 \eta 
    \right).
$$
If $B \leq B_c$, $V(G^{\rm ext}) \leq (3 B_c / 2)^2 h^{-2}$
and $|S^{\rm ext}| \geq 1 / h^{1 - b / 2}$,
then, by Lemma~\ref{jerome}, an event which is much more 
unlikely than ${\cal I}$ has to occur:
there is $\delta > 0$ such that, for $h$ small enough,
\begin{align*}
    \mu_{\Lambda_h, -, h}\left(
        |S^{\rm ext}| \geq 1 / h^{1 - b / 2},\,
        V(G^{\rm ext}) \leq (3 B_c / 2)^2 h^{-2}
    \right)
    &\leq \mu_{\lambda_h, -, h}({\cal I})
    \sum_{k \geq 1 / h^{1 - b / 2}} \exp\left\{
        -\delta k / h^b 
    \right\}\\
    &\leq 2 \mu_{\lambda_h, -, h}({\cal I})
    \exp\left\{
        -\delta / h^{1 + b / 2}
    \right\}.
\end{align*}
Finally, if $B \leq B_c$, $V(G^{\rm ext}) \leq (3 B_c / 2)^2 h^{-2}$
and $k = |S^{\rm ext}| < 1 / h^{1 - b / 2}$,
then, since $\sigma \not\in {\cal R}$,
it follows from Lemma~\ref{raphael}
that the smallest Wulff shapes to contain
its external vertebrate contours $\Gamma_j$ have 
a total square root volume larger than $B_+ / h$.
Again, using the skeleton properties
and the fact that each of these contours
encloses a volume which is larger than $1 / h^{2b}$,
we get that the smallest Wulff shapes
to contain the associated skeletons $\gamma_j$
have total volume larger than
$$
    \sqrt{1 - C h^{b - r}}\: {B_+ \over h}
    \geq \sqrt{1 - C h^{b / 2}}\: {B_+ \over h}
    \geq {\tilde B_+ \over h}
$$
for some positive constant $C$,
any $\tilde B_+ < B_+$
and $h$ small enough. 
We choose $\tilde B_+ > B_c$
such that
$$
    {B_c^2 \over \tilde B_+} + \tilde B_+
    \geq 2 B_c \left(
        1 + {3 \eta \over 2}
    \right).
$$
Writing $(B_j / h)^2$ for the phase
volume of each single skeleton $\gamma_j$
and $(B_{j, out} / h)^2$ for the volume
of the smallest Wulff shape to contain it,
we have, 
using again Proposition~\ref{william},
for any small enough $\epsilon > 0$
and since $B \leq B_c < \tilde B_+$,
\begin{align*}
    {\cal W}(S) - (1 + \epsilon) h m^*_\beta \check V(S)
    &\geq {\cal W}(S^{\rm ext}) - (1 + \epsilon) h m^*_\beta \check V(S^{\rm ext})\\
    &\geq \sum_{j < k} {w_\beta \over 2} \left(
        {(B_j / h)^2 \over B_{j, out} / h} + B_{j, out} / h
    \right) - (1 + \epsilon) h m^*_\beta (B / h)^2\\
    &\geq {w_\beta \over 2h} \left[
        {\sum_{j < k}  B_j^2  \over \sum_{j < k} B_{j, out}} + \sum_{j < k} B_{j, out}
        - (1 + \epsilon) {B^2 \over B_c}
    \right]\\
    &\geq {w_\beta \over 2h} \left[
        {B^2 \over \sum_{j < k} B_{j, out}} + \sum_{j < k} B_{j, out}
        - (1 + \epsilon) {B^2 \over B_c}
    \right]\\
    &\geq {w_\beta \over 2h} \left[
        {B^2  \over \tilde B_+} + \tilde B_+
        - (1 + \epsilon) {B^2 \over B_c}
    \right]
    = {w_\beta \over 2h} \left[
        B^2\left(
            {1  \over \tilde B_+} - {1 + \epsilon \over B_c}
        \right)
        + \tilde B_+
    \right]\\
    &\geq {w_\beta \over 2h} \left[
        B_c^2\left(
            {1  \over \tilde B_+} - {1 + \epsilon \over B_c}
        \right)
        + \tilde B_+
    \right]
    \geq {w_\beta \over 2h} \left[
        2 B_c \left(
            1 + {3 \eta \over 2}
        \right)
        - B_c (1 + \epsilon)
    \right]\\
    &\geq {w_\beta B_c\over 2h} (1 + 2 \eta) = {A \over h} (1 + 2 \eta).
\end{align*}

We conclude with Lemma~\ref{thomas} by summing on all the possible 
values of the integer $\check V(S) < 2 (B_{\max} / h)^2$:
\begin{align*}
    \mu_{\Lambda_h, -, h}\bigl({\cal X}^c\bigr)
    &\leq \mu_{\Lambda_h, -, h}\bigl({\cal I}\bigr)\left[
        2 \left(
            B_{\max} \over h
        \right)^2 \exp\left\{
            -(1 - \epsilon) {\beta A \over h} (1 + 2 \eta)
        \right\}
        + 2 \exp\left\{
            - \delta / h^{1 + b / 2}
        \right\}
    \right]\\
    &\leq \mu_{\Lambda_h, -, h}\bigl({\cal I}\bigr)\exp\left\{
        -{\beta A \over h} (1 + \eta)
    \right\}
\end{align*}
for $\epsilon$ chosen small enough
and all small enough $h$.
\qed

\medskip\par
It follows that, starting from
$\mu_{\Lambda_h, -, h}(\cdot\,|{\cal R}) = \mu_{\cal R}$
and with large probability,
our process $X_{\Lambda_h, -, h}$
cannot escape from ${\cal X}$ 
within time 
$$
    t_1 = \exp\left\{
        {\beta A \over h}\left(
            1 + {\eta \over 2}
        \right)
    \right\}
$$
for $h$ large enough.
Indeed, since we assumed that for all $x \in \Lambda_h$
and for all $\sigma \in \Omega_{\Lambda_h}$
$$
    w(\sigma, \sigma^x)
    \leq w_{\max},
$$
the number of jumps of the process $X_{\Lambda_h, - , h}$ 
within time $t_1$ is dominated by a Poisson random
variable $N_1$ with mean
$$
    \lambda_1 = |\Lambda_h| w_{\max} t_1
$$
and for which
$$
    P\left(
        N_1 \geq e \lambda_1
    \right)
    \leq e^{- e \lambda_1} E\left[
        e^{N_1}
    \right]
    = \exp\bigl\{
        - e\lambda_1 -\lambda_1 + e \lambda_1
    \bigr\}
    = e^{-\lambda_1}
    \leq {1 \over \lambda_1}.
$$
Since from the previous lemma it holds,
for $h$ small enough and all $t \geq 0$,
\begin{align*}
    \mathds{P}_{\mu_{\cal R}}\bigl(
        X_{\Lambda_h, -, h}(t) \in {\cal X}^c
    \bigr)
    &= \sum_{\sigma \in {\cal  R}}
    \sum_{\sigma' \not\in {\cal X}}
    {\mu_{\Lambda_h, -, h}(\sigma) \over \mu_{\Lambda_h, -, h}({\cal R})}
    \mathds{P}_\sigma\bigl(
        X_{\Lambda_h, -, h}(t) = \sigma'
    \bigr)\\
    &= \sum_{\sigma' \not\in {\cal X}}
    {\mu_{\Lambda_h, -, h}(\sigma') \over \mu_{\Lambda_h, -, h}({\cal R})}
     \sum_{\sigma \in {\cal  R}} \mathds{P}_{\sigma'}\bigl(
        X_{\Lambda_h, -, h}(t) = \sigma
    \bigr)\\
    &= \sum_{\sigma' \not\in {\cal X}}
    {\mu_{\Lambda_h, -, h}(\sigma') \over \mu_{\Lambda_h, -, h}({\cal R})}
     \mathds{P}_{\sigma'}\bigl(
        X_{\Lambda_h, -, h}(t) \in {\cal R}
    \bigr)\\
    &\leq \sum_{\sigma' \not\in {\cal X}}
    {\mu_{\Lambda_h, -, h}(\sigma') \over \mu_{\Lambda_h, -, h}({\cal R})}
    = {\mu_{\Lambda_h, -, h}\bigl(
        ({\cal R} \cup {\cal S})^c
    \bigr) \over \mu_{\Lambda_h, -, h}({\cal R})}\\
    &\leq {
        \mu_{\Lambda_h, -, h}({\cal I}) \exp\left\{
            -{\beta A \over h} (1 + \eta)
        \right\}
    \over
        \mu_{\Lambda_h, -, h}({\cal I})
    }
    = \exp\left\{
        -{\beta A \over h} (1 + \eta)
    \right\},
\end{align*}
we conclude, with $T_{{\cal X}^c}$
the exit time from ${\cal X}$,
$$
    \mathds{P}_{\mu_{\cal R}}\left(
        T_{{\cal X}^c} \leq t_1
    \right)
    \leq {1 \over \lambda_1} 
    + e \lambda_1 \exp\left\{
        -{\beta A \over h} (1 + \eta)
    \right\}
$$
for $h$ small enough and
\begin{lmm}\label{aude}
    Given $B_+ > B_c$ and $b < 1 / 4$, 
    setting $\eta > 0$ such that
    $$
        {B_c^2 \over B_+} + B_+ 
        = 2 B_c (1 + 2 \eta),
    $$ 
    it holds 
    $$
        \mathds{P}_{\mu_{\cal R}}\left(
            T_{{\cal X}^c} \leq \exp\left\{
                {\beta A \over h}\left(
                    1 + {\eta \over 2}
                \right)
            \right\}
        \right)
        \leq \exp\left\{
                - {\beta A \over h}\,
                {\eta \over 3}
            \right\}
    $$
    for $h$ small enough.
\end{lmm}

\subsection{Entering \bm{${\cal S}$} or \bm{${\cal R}$}}
\begin{lmm}\label{kamra}
    Given $\delta > 0$,
    one can choose $B_+$ close enough to $B_c$
    to have, for all small enough $h$,
    $$
        \mu\bigl(
            {\cal R} \cap {\cal S}
        \bigr) \leq \mu({\cal I})\exp\left\{
            -{\beta A  - \delta \over h}
        \right\}
    $$
    and
    $$
        \phi^*_{{\cal S} \backslash {\cal R}}
        \vee \phi^*_{{\cal R} \backslash {\cal S}}
        \leq \exp\left\{
            -{\beta A  - \delta \over h}
        \right\}.
    $$
\end{lmm}

\medskip\par\noindent
{\bf Proof:} Consider, for any $B_+ > B_c$,
$\sigma$ in ${\cal R} \cap {\cal S}$
and its associated skeleton collection $S$. 
Since $\sigma \in {\cal S}$,
the isoperimetric property of the Wulff shape
implies that, for any $\epsilon > 0$,
\begin{equation}\label{nora}
    {\cal W}(S) \geq (1 - \epsilon) w_\beta {B_- \over h}
\end{equation}
for all small enough $h$.
Also, since $\sigma \in {\cal R}$,
it holds
\begin{equation}\label{djamila}
    \check V(S) \leq (1 + \epsilon) \left(
        {B_+ \over h}
    \right)^2
\end{equation}
for all small enough $h$.
Then, by Lemma~\ref{thomas},
\begin{align*}
    \mu\bigl({\cal R} \cap {\cal S})
    &= \leq \mu({\cal I}) \exp\left\{
        - {\beta \over h} \Bigl[
            (1 - \epsilon)^2 w_\beta B_- 
            - (1 + \epsilon)^2 m^*_\beta B_+^2 
        \Bigr]
    \right\}\\
    &= \mu({\cal I}) \exp\left\{
        - {\beta A \over h} \left[
            2 (1 - \epsilon)^2 {B_- \over B_c}
            - (1 + \epsilon)^2 \left(
                {B_+ \over B_c}
            \right)^2
        \right]
    \right\}.
\end{align*}
Choosing $\epsilon$ small enough and 
$B_+$ close enough to $B_c$
we get 
$$
    \mu\bigl(
        {\cal R} \cap {\cal S}
    \bigr) \leq \mu({\cal I})\exp\left\{
        -{\beta A  - \delta \over h}
    \right\}
$$
for all small enough $h$.

We proceed in the same way
and use inequality~\eqref{sophie_2}
from page \pageref{sophie_2} 
to bound $\phi^*_{{\cal S} \backslash {\cal R}}$.
For all $\sigma \in {\cal S} \backslash {\cal R}$
associated with a skeleton family $S$
it holds 
$$
    e^*_{{\cal S} \backslash {\cal R}}(\sigma)
    \leq |\Lambda_h|w_{\max},
$$
since $\sigma \in {\cal S}$,
inequality~\eqref{nora}
is in force for any $\epsilon > 0$
and all small enough $h$,
and $e^*_{{\cal S} \backslash {\cal R}}(\sigma) = 0$
unless there is $x \in \Lambda_h$ such that 
$\sigma^x \in {\cal R}$
so that inequality~\eqref{djamila}
is also in force 
for all small enough $h$.
Hence, using Lemma~\ref{daniele}
with a small enough $\delta'$ in place of $\delta$,
\begin{align*}
    \phi^*_{{\cal S} \backslash {\cal R}}
    &\leq {
        |\Lambda_h| w_{\max} \mu({\cal I})
    \over
        \mu({\cal S} \backslash {\cal R})
    } \exp\left\{
        - {\beta \over h} \Bigl[
            (1 - \epsilon)^2 w_\beta B_- 
            - (1 + \epsilon)^2 m^*_\beta B_+^2 
        \Bigr]
    \right\}\\
    &\leq {
        |\Lambda_h| w_{\max} \mu({\cal I})
    \over \mu\Bigl(
           E^h_{B_{\max} / (1 + \delta'), \delta'}
        \Bigr)
    } \exp\left\{
        - {\beta \over h} \Bigl[
            (1 - \epsilon)^2 w_\beta B_- 
            - (1 + \epsilon)^2 m^*_\beta B_+^2 
        \Bigr]
    \right\}\\
    &\leq  \exp\left\{
        - {\beta \over h} \Bigl[
            (1 - \epsilon)^2 w_\beta B_- 
            - (1 + \epsilon)^2 m^*_\beta B_+^2 
        \Bigr]
    \right\}
\end{align*}
for all $h$ small enough,
and we conclude in the same way.

Finally, since from inequality~\eqref{sophie_1}
it holds
---with the convention $w(\sigma, \sigma') = 0$
for all $\sigma \neq \sigma'$ such that
$\sigma' \neq \sigma^x$ for all $x$ in $\Lambda_h$---
$$
    \phi^*_{{\cal R} \backslash {\cal S}}
    \leq \sum_{\sigma \in {\cal R} \backslash {\cal S}}
    \sum_{\sigma' \in {\cal S}}
    \mu_{{\cal R} \backslash {\cal S}}(\sigma) w(\sigma, \sigma')
    \leq {1 \over \mu({\cal I})} \sum_{\sigma' \in {\cal S}} \mu(\sigma')
    \sum_{\sigma \in {\cal R}} w(\sigma', \sigma)
$$
we can use the same arguments to bound $\phi^*_{{\cal R} \backslash {\cal S}}$.
\qed

\subsection{Upper bounds for soft capacities}
Given $\delta > \delta' > 0$, assume that we chose
$B_+ > B_+'$ associated with ${\cal R} \supset {\cal R}'$
and ${\cal S} \supset {\cal S}'$ as in Lemma~\ref{kamra}.
We use the variational principle \eqref{guillaume_1}
to get an upper bound on $C_\kappa^\lambda({\cal R}, {\cal S})$.
We build then a test function $f : {\cal X} \rightarrow \mathds{R}$
with 
$$
    f(\sigma) = \left\{
        \begin{array}{ll}
            1 & \hbox{if $\sigma \in {\cal R}' \setminus {\cal S}'$,}\\
            1/2 & \hbox{if $\sigma \in {\cal R}' \cap {\cal S}'$,}\\
            0 & \hbox{if $\sigma \in {\cal S}' \setminus {\cal R}'$,}\\
            1/2 & \hbox{if $\sigma \not\in {\cal X}' = {\cal R}' \cup {\cal S}'$}.
        \end{array}
    \right.
$$
Note that, for all $x \in \Lambda_h$,
if $\sigma$ and $\sigma^x$ both belong to ${\cal X}'$
but neither of them is in ${\cal R}' \cap {\cal S}'$,
then $f(\sigma) = f(\sigma^x)$.
Hence, by Lemma~\ref{kamra} and Lemma~\ref{alphonsine}
with $\delta'$ and $\eta'$ in place of $\delta$ and $\eta$,
\begin{align*}
    {C_\kappa^\lambda({\cal R}, {\cal S}) \over \mu({\cal R})}
    &\leq
    {
        \mu\bigl({\cal R}' \cap {\cal S}'\bigr)
    \over 
        \mu({\cal I})
    }\, {|\Lambda_h| w_{\max} \over 4}
    + {
        \mu\bigl({\cal R} \cap {\cal S}\bigr)
    \over 
        \mu({\cal I})
    }(\kappa + \lambda)
    + {\mu\bigl({\cal X} \setminus {\cal X}'\bigr) \over \mu({\cal I})}\Bigl[
        |\Lambda_h| w_{\max} + \kappa + \lambda
    \Bigr]\\
    &\leq e^{-(\beta A - \delta') / h}\, {|\Lambda_h| w_{\max} \over 4}
    + e^{-(\beta A - \delta) / h}(\kappa + \lambda)
    + e^{-\beta A (1 + \eta') / h}\Bigl[
        |\Lambda_h| w_{\max} + \kappa + \lambda
    \Bigr]
\end{align*}
for all small enough $h$.
Since $\delta'$ can be chosen arbitrarily small,
we conclude
\begin{lmm}\label{mohamed}
    Given $\delta > 0$,
    choosing $B_+$ close enough to $B_c$
    to have, for $h$ small enough,
    $$
        \phi^*_{{\cal R} \backslash {\cal S}}
        \vee
        \phi^*_{{\cal S} \backslash {\cal R}}
        \leq \exp\left\{
            -{\beta A - \delta \over h}
        \right\},
    $$
    choosing also $\kappa = \kappa(h)$
    and $\lambda = \lambda(h)$ 
    such that 
    \begin{equation}\label{christine}
        \lim_{h \rightarrow 0} \kappa(h) e^{\delta / h}
        = \lim_{h \rightarrow 0} \lambda(h) e^{\delta / h}
        = 0,
    \end{equation}
    for all $\epsilon > 0$, there is $h_0 > 0$ such that 
    $$
        {C_\kappa^\lambda({\cal R}, {\cal S}) \over \mu({\cal R})}
        \leq \exp\left\{
            - {\beta A - \epsilon \over h}
        \right\}
    $$
    for all $h < h_0$.
\end{lmm}

\section{Upper bounds for local relaxation times}\label{julia}
\subsection{On the metastable side}
We prove in this section 
that for any $\delta > 0$
one can choose $B_+$ close enough
to $B_c$ in such a way
that the local relaxation time $1 / \gamma_{\cal R}$
is smaller than $e^{\delta / h}$
for $h$ small enough.
We use to this end a small parameter $d > 0$,
the value of which will depend on $\delta$
and will be used to choose $B_+$.
Given a finite family $F$ of disjoint Wulff shapes
$$
    x_j + W_{\rho_j / h} \subset \Lambda_h,
    \qquad j < k,
$$
with $k < 1 / h^{1 - b / 2}$,
we build a sequence of smaller disjoint Wulff shapes
$$
    x_{j, l} + W_{\rho_{j, l} / h}
    = x_j + W_{(\rho_j - d l) / h},
    \qquad j < k,
    \qquad l < l_0,
$$
with (recall Equation~\eqref{charles} from page~\pageref{charles})
$$
    l_0 = \left\lceil
        {\rho_{\max}  \over d}
    \right\rceil = \left\lceil 
        {2 B_{\max} \over d w_\beta}
    \right\rceil
$$
and the convention that for all $\rho < 0$ 
and all $x$ in $\mathds{R}^2$,
$W_\rho$ and $x + W_\rho$ both stand for the empty set.
For $l < l_0$ we denote by $W_l(F)$ their union:
$$
    W_l(F) = \bigcup_{j < k} x_{j, l} + W_{\rho_{j, l} / h}.
$$
Next, we associate with each $l < l_0$
a family of disjoint annuli of the lattice, the union of which is 
$$
    A_l(F) = \mathds{Z}^2 \cap \bigcup_{j < k} x_{j, l} + \Bigl(
        W_{\rho_{j, l} / h} \setminus W_{(\rho_{j, l} - 4d) / h}
    \Bigr).
$$
We also define, independently of $F$,
a further sequence of Wulff shapes
$$
    W_{\rho'_l / h} = W_{(\rho_{\max} - (l - l_0) d) / h},
    \qquad l_0 \leq l < 2 l_0,
$$
and, for each $l_0 \leq l < 2 l_0$, we set
$$
    A_l(F) = \mathds{Z}^2 \cap \left(
        W_{\rho'_l / h} \setminus W_{(\rho'_l - 4d) / h}
    \right).
$$
We order the sites of such an annulus
$A_l(F)$ with $l \geq l_0$
by ordering first the angles,
then the radii:
for $x$ and $y$ in $A_l(F)$ 
we say that $x$ is lower than $y$
if the angle between the horizontal
and the half-line that goes through $x$
and starts in the annulus center
is smaller that the similar angle 
associated with $y$ and, 
if both angles are equal,
we say that $x$ is lower than $y$
if so are the associated distances
to the annulus center.
For $l < l_0$ we order
similarly the sites in $A_l(F)$
by ordering first the annuli,
then the angles and the radii.

For $\sigma \in \Omega_{\Lambda_h}$ 
and $l_0 \leq l < 2 l_0$,
we consider the collection ${\cal C}$
of the external contours $\Gamma$ of $\sigma$
that enclose some $x$ outside $A_l(F)$,
we call $E_l(\sigma)$ the subset of $\mathds{Z}^2$
made of all sites enclosed in some $\Gamma \in {\cal C}$
and we call $\bar E_l(\sigma)$ the subset of $\mathds{Z}^2$
made of all the sites in $E_l(\sigma)$ 
or having a nearest neighbour in $E_l(\sigma)$.
We define then the ``block'' $A_l(F, \sigma)$ by 
$$
    A_l(F, \sigma) = A_l(F) \setminus \bar E_l(\sigma).
$$
To avoid ambiguities,
we will denote by $\nu_{A_l(F, \sigma), \sigma, h}$,
rather than identify with $\mu_{A_l(F, \sigma), \sigma, h}$,
the law of the $\Omega_{\Lambda_h}$-valued random variable $M$
for which $M$ and $\sigma$ coincide outside $A_l(F, \sigma)$
and the restriction of $M$ to $A_l(F, \sigma)$ 
is drawn according to 
$$
    \mu_{A_l(F, \sigma), \sigma, h}
    = \mu_{A_l(F, \sigma), -, h}.
$$
For $\sigma \in \Omega_{\Lambda_h}$ 
and $l < l_0$, we make a different block construction
by considering the collection ${\cal C}'$
of the external contours $\Gamma$ of $\sigma$
that enclose some $x$ in $\mathds{Z}^2 \setminus W_l(F)$.
We call $E'_l(F, \sigma)$ the subset of $\mathds{Z}^2$
made of all sites enclosed in some $\Gamma \in {\cal C}'$
and, similarly, we call $\bar E'_l(F, \sigma)$ the subset of $\mathds{Z}^2$
made of all the sites in $E'_l(F, \sigma)$ 
or having a nearest neighbour in $E'_l(F, \sigma)$.
We then set 
$$
    A_l(F, \sigma) = A_l(F) \setminus \bar E'_l(F, \sigma)
$$
and, similarly, 
we denote by $\nu_{A_l(F, \sigma), \sigma, h}$,
the law of the $\Omega_{\Lambda_h}$-valued random variable $M$
for which $M$ and $\sigma$ coincide outside $A_l(F, \sigma)$
and the restriction of $M$ to $A_l(F, \sigma)$ 
is drawn according to $\mu_{A_l(F, \sigma), \sigma, h}$.
Note that, in both the cases $l < l_0$ and $l \geq l_0$,
DLR equations imply that, 
if $M$ is drawn according to $\mu_{\Lambda_h, -, h}$
and $M'$ is drawn according to $\nu_{A_l(F, M), M, h}$,
then $M$ and $M'$ have the same law.

Given $F$, we now associate with each
$\sigma$ in $\Omega_{\Lambda_h}$
a block path $\Pi_\sigma$
by setting first $M_0 = \sigma$,
drawing then, for each $l < 2 l_0$,
the milestone $M_{l + 1}$ according to $\nu_{A_l(F, M_l), M_l, h}$
and connecting finally each milestone $M_l$ with $M_{l + 1}$
along the canonical path $\pi^l_{M_l, M_{l + 1}}$ in $\Omega_{A_l(F)}$
associated with the ordered set $A_l(F)$.

\begin{lmm}
    There is a positive constant $C$ 
    such that, for any $d > 0$,
    all $\sigma_0$ in $\Omega_{\Lambda_h}$
    and all $x$ in $\Lambda_h$,
    $$
        {1 \over \mu_{\Lambda_h, -, h}(\sigma_0)}
        \sum_{\sigma \in \Lambda_h} \mu_{\Lambda_h, -, h}(\sigma)
        P\bigl(
            (\sigma_0, \sigma_0^x) \in \Pi_\sigma
        \bigr)
        \leq 8 \exp\left\{
            {Cd \over h}
        \right\}.
    $$
\end{lmm}

\medskip\par\noindent
{\bf Proof:}
We first note that,
for $(\sigma_0, \sigma_0^x)$ to belongs to $\Pi_\sigma$,
there is to be some $l < 2 l_0$
such that $x$ lies in $A_l(F)$
and $(\sigma_0, \sigma_0^x)$ belongs to $\pi^l_{M_l, M_{l + 1}}$.
Since our annuli are of ``width'' $4d$ 
and their linear size decreases by $d$ 
in each of our two annulus sequences,
their are 8 such $l$ at most.
Now, if $x \in A_l(F)$,
with
$$
    \Omega_{l, \sigma_0}
    = \left\{
        \sigma \in \Omega_{\Lambda_h} :
        \forall x \not\in A_l(F, \sigma_0),\,
        \sigma(x) = \sigma_0(x)
    \right\}
$$
then, by DLR equations and Lemma~\ref{fabio},
there is $C > 0$ such that
\begin{align*}
    &{1 \over \mu_{\Lambda_h, -, h}(\sigma_0)}
    \sum_{\sigma \in \Lambda_h} \mu_{\Lambda_h, -, h}(\sigma)
    P\bigl(
        (\sigma_0, \sigma_0^x) \in \pi^l_{M_l, M_{l + 1}}
    \bigr)\\
    &\qquad = {1 \over \mu_{\Lambda_h, -, h}(\sigma_0)}
    \sum_{\sigma_l, \sigma_{l + 1} \in \Omega_{l, \sigma_0}}
    \sum_{\sigma \in \Lambda_h} \mu_{\Lambda_h, -, h}(\sigma)
    P\bigl(M_l = \sigma_l) \nu_{A(F, \sigma_l), \sigma_l, h}(\sigma_{l + 1})
    \mathds{1}\!\left\{
        (\sigma_0, \sigma_0^x) \in \pi^l_{\sigma_l, \sigma_{l + 1}}
    \right\}\\
    &\qquad = \sum_{\sigma_l, \sigma_{l + 1} \in \Omega_{l, \sigma_0}}
    {\mu_{\Lambda_l, -, h}(\sigma_l) \over \mu_{\Lambda_h, -, h}(\sigma_0)}
    \nu_{A(F, \sigma_l), \sigma_l, h}(\sigma_{l + 1})
    \mathds{1}\!\left\{
        (\sigma_0, \sigma_0^x) \in \pi^l_{\sigma_l, \sigma_{l + 1}}
    \right\}\\
    &\qquad = {1 \over \nu_{A(F, \sigma_0), \sigma_0, h}(\sigma_0)}
    \sum_{\sigma_l, \sigma_{l + 1} \in \Omega_{l, \sigma_0}}
    \nu_{A(F, \sigma_0), \sigma_0, h}(\sigma_l)
    \nu_{A(F, \sigma_0), \sigma_0, h}(\sigma_{l + 1})
    \mathds{1}\!\left\{
        (\sigma_0, \sigma_0^x) \in \pi^l_{\sigma_l, \sigma_{l + 1}}
    \right\}\\
    &\qquad \leq \exp\left\{
        {Cd \over h}
    \right\}.
\end{align*}
\qed

Given $\sigma$ and $\sigma'$ in ${\cal R}$
we will couple two such block paths $\Pi_\sigma$
and $\Pi_\sigma'$ associated with two {\it random} families
$F$ and $F'$.
We will consider a ``good event'' $E_{\sigma, \sigma'}$
for which $\Pi_\sigma$ and $\Pi_\sigma'$
will stay in ${\cal R}$ and will end
in the same $M_{2l_0} = M'_{2l_0}$.
Then, conditionally to $E_{\sigma, \sigma'}$,
we can build a block path $\Pi_{\sigma, \sigma'}$
in ${\cal R}$ and
from $\sigma$ to $\sigma'$
by concatenation of $\Pi_\sigma$,
from $\sigma$ to $M_{2l_0}$,
and the reversed image of $\Pi_{\sigma'}$,
from $M'_{2l_0} = M_{2l_0}$ to $\sigma'$.
Since the previous lemma is uniform in $F$,
we will get, for all $\sigma_0$ and $\sigma_0^x$
in ${\cal R}$
\begin{align*}
    &{1 \over \mu_{\cal R}(\sigma_0) \vee \mu_{\cal R}(\sigma_0^x)}
    \sum_{\sigma, \sigma' \in {\cal R}}
    \mu_{\cal R}(\sigma)\mu_{\cal R}(\sigma')
    P\left(
        (\sigma_0, \sigma_0^x) \in \Pi_{\sigma, \sigma'}
        \bigm| E_{\sigma, \sigma'}
    \right)\\
    &\qquad\leq {\mu({\cal R}) \over \mu(\sigma_0)}
    \sum_{\sigma \in {\cal R}} {\mu(\sigma) \over \mu({\cal R})}
    {
        P\bigl((\sigma_0, \sigma_0^x) \in \Pi_\sigma\bigr)
    \over
        \min_{\sigma, \sigma' \in {\cal R}} P\bigl(E_{\sigma, \sigma'}\bigr)
    }
    + {\mu({\cal R}) \over \mu(\sigma_0^x)}
    \sum_{\sigma' \in {\cal R}} {\mu(\sigma') \over \mu({\cal R})}
    {
        P\bigl((\sigma_0^x, \sigma_0) \in \Pi_\sigma\bigr)
    \over
        \min_{\sigma, \sigma' \in {\cal R}} P\bigl(E_{\sigma, \sigma'}\bigr)
    }\\
    &\qquad\leq {16 e^{Cd / h} \over \min_{\sigma, \sigma' \in {\cal R}} P(E_{\sigma, \sigma'})}.
\end{align*}
In view of inequality~\eqref{caroline} at page~\pageref{caroline},
we will need a lower bound on $P(E_{\sigma, \sigma'})$.

Before building $E_{\sigma, \sigma'}$ and giving such a lower bound,
let us first explain in which sense $F$ and $F'$
are random.
To sample $F$ of size $k < 1 / h^{1 - b/2}$,
we first sample $k$ uniformly,
then we sample the centers $x_j$ uniformly
in $B_{\max}W / h$, and, finally,
we sample the $\rho_j$ uniformly in $[0, \rho_+]$,
with 
$$
    \rho_+ = {2 B_+ \over w_\beta}\,,
$$
and conditionally to our non-intersection constraint.
We sample $F'$ independently and in the same way.

We say that $F$ is adapted to $\sigma$
if the Wulff shapes of $F$ contain the external 
vertebrate contours of~$\sigma$.
This is the first requirement for our good 
event $E_{\sigma, \sigma'}$
and it happens with a probability larger than
\begin{equation}\label{paul}
    \left(
        C 
    \over 
        |\Lambda_h| (\rho_{\max} / h)
    \right)^{1 / h^{1 - b/2}}
    \geq
    e^{- \delta / (8h)}
\end{equation}
for some $C > 0$ and all small enough $h$.
We assume in what follows that $F$ is
adapted to $\sigma$.

The next requirement for $E_{\sigma, \sigma'}$
is that for each $l < l_0$,
$M_{l + 1}$ has no vertebrate contour
to enclose a site in the annulus union
$$
    A^2_l(F)
    = A_l(F) \setminus W_{l + 2}(F),
$$
with the convention $W_{l + 2}(F) = \emptyset$
for $l + 2 \geq l_0$.
Provided that $B_+$ is close enough to $B_c$
to have
\begin{equation}\label{rachel}
    \phi(B_+ - 4d) < \phi(B_+),
\end{equation}
using inductively FKG inequality
together with Estimate~\eqref{adrien_2} from page~\pageref{adrien_2}
with a small enough $\epsilon$ depending of $l_0$,
then $d$,
this occurs with probability $e^{- \delta /(8h)}$ at least
for all small enough $h$. 

Provided that the same requirements are satisfied
for $F'$ and $M'_l$ with $l < l_0$,
it holds that the milestones $M_{l_0}$ and $M'_{l_0}$
are both in ${\cal I}$.
It is also the case that $\Pi_\sigma$ and $\Pi_\sigma'$
did not escape ${\cal R}$ up to this point,
where we can start to introduce
some dependence between them. 

Assuming that our previous requirements
for $E_{\sigma, \sigma'}$ were satisfied,
the next one is that $M_{l_0 + 1}$ and $M'_{l_0 + 1}$
are still in ${\cal I}$ and coincide on the annulus
$$
    A^2 = \mathds{Z}^2 \cap \left(
        W_{\rho_{\max} /  h} \setminus W_{(\rho_{\max} - 2d) / h}
    \right).
$$
For $d$ small enough,
inequality~\eqref{adrien_4},
DLR equations 
and FKG inequality 
show that this happens 
with a non-negligible probability.
Indeed, since $M_{l_0}$ and $M'_{l_0}$
are in ${\cal I}$,
the restrictions to 
$$
    A^3 = \mathds{Z}^2 \cap \left(
        W_{\rho_{\max} /  h} \setminus W_{(\rho_{\max} - 3d) / h}
    \right)
$$
of $M_{l_0 + 1}$ and $M'_{l_0 + 1}$
are both dominated
by a that of a random configuration $\xi$
drawn according to $\mu_{A^5, -, h}$,
with 
$$
    A^5 = \mathds{Z}^2 \cap \left(
        W_{\rho_{\max} /  h} \setminus W_{(\rho_{\max} - 5d) / h}
    \right).
$$
Hence, we can partially sample them
first by drawing the external contours $\Gamma$ of $\xi$
that will cross the boundary of $A^3$,
then by drawing the common restriction
of $\xi$, $M_{l_0 + 1}$ and $M'_{l_0 + 1}$
to $A^3 \setminus \bar E$
according to $\mu_{A^3 \backslash \bar E, -, h}$,
with $\bar E$ the set of all sites that are enclosed
by one of these $\Gamma$ or that are a nearest neighbour of such
a site. 
Since, by~\eqref{adrien_4},
$\xi$  is in ${\cal I}$ with a non-negligible probability,
larger than $e^{-\beta\epsilon / h}$, for all small enough $h$,
this gives the same lower bound for this new requirement. 

Our last requirement,
which includes the previous one,
is that, for all $l_0 \leq l < 2 l_0$,
the milestones $M_{l +1}$ and $M'_{l + 1}$
are in ${\cal I}$
and coincide on the annulus
$$
    A^{l - l_0 + 2} = \mathds{Z}^2 \cap \left(
        W_{\rho_{\max} /  h} \setminus W_{(\rho_{\max} - (l - l_0 + 2)d) / h}
    \right).
$$
Provided that our previous set of requirements was satisfied,
this implies that the whole paths $\Pi_\sigma$ 
and $\Pi_\sigma'$ all along remain in ${\cal R}$
and end in a same configuration $M_{2 l_0} = M'_{2 l_0}$,
and this happens,
repeating inductively the previous argument,
with a probability $e^{-\delta /(8h)}$ at least
for $h$ small enough.

Using inequality~\eqref{caroline}
from page~\pageref{caroline}, we get that,
for any small enough $d$,
if $B_+$ is close enough to $B_c$ for 
inequality~\eqref{rachel} to be in force,
then 
$$
    {1 \over \gamma_{\cal R}}
    \leq {8 (B_{\max} / h)^2 \over w_{\min}} 16 e^{Cd / h} e^{5 \delta / (8h)}
$$
for some positive constant $C$
that does not depend on $d$ and all small enough $h$.
Choosing $d$ small enough to have $Cd < 2 / 8$ we conclude
\begin{lmm}\label{mimi}
    Given $\delta > 0$, 
    one can choose $B_+$ close enough to $B_c$
    to have 
    $$
        {1 \over \gamma_{\cal R}} \leq e^{\delta / h}.
    $$ 
\end{lmm}

\subsection{On the stable side}
The goal of this section is to show
\begin{lmm}\label{emma}
    Given $\delta > 0$, 
    one can choose $B_+$ close enough to $B_c$
    to have 
    $$
        {1 \over \gamma_{\cal S}} \leq e^{\delta / h}.
    $$ 
\end{lmm}

The proof is similar to that on the metastable side,
with some simplifications and some extra complications.
We will only indicate the main differences.

Simplifications come from the fact 
that we will only have to build annular blocks:
we will not need union of annuli anymore.
Similarly to the previous case,
we will use these blocks
to build a path of expanding,
rather than shrinking, contours,
before using the same shrinking blocks
to make the final milestones
of two block paths coincide.

There are only two kind of complications.
We will first need another sequence
of shrinking blocks
to ensure that, starting from 
$\sigma \in {\cal S}$ for which 
there is a large contour that
encloses a slightly subcritical Wulff shape,
we will only see ``the plus-phase'' 
on the internal border of this ``large'' Wulff shape
at the end of the associated first block path.
This is needed to use inequality~\eqref{adrien_1}
of page \pageref{adrien_1} with our second,
expanding, block sequence
---the analogue
of the first shrinking sequence
on the metastable side---
to obtain, as last milestone
associated with the last block
of this second block sequence,
a configuration 
with only one vertebrate contour,
close to the boundary of $\Lambda_h$,
outside our slightly subcritical Wulff shape.
We encounter the second complication
in building this second, expanding, 
block sequence:
since our expanding blocks have to be contained
in $\Lambda_h$ and eventually coincide
with its boundary, except if we start
with an annular block centered
on the origin,
we cannot have concentric blocks.
Because the overlapping properties of our blocks
are crucial for the inductive parts
of our arguments in giving a lower bound
for our good event, 
there is an issue.

Here is the key lemma we will use to solve it.
It says that two non-concentric Wulff shapes
on the same side of a common tangent
are such that the core of the largest one
is contained in the bulk of the smallest one.
\begin{lmm}\label{sara}
    Let ${\bf n}=(\cos \theta, \sin \theta)$ 
    be the external normal associated
    with a Wulff shape $x + W_\rho$
    and $y$ in $x + \partial W_\rho$.
    For a positive $d < \rho / 3$, let $x'$ in $\mathds{R}^2$ 
    be such that ${\bf n}$
    is also the external normal associated
    with the Wulff shape $x' + W_d$
    and $x$ in $x' + \partial W_d$.
    Then the Wulff shapes $x + W_\rho$
    and $x' + W_{\rho + d}$
    are on the same side of a common tangent
    in $y$ and it holds 
    $$
        x' + W_{(\rho +  d) - 4d}
        = x' + W_{\rho - 3d}
        \subset x + W_{\rho - 2d}.
    $$
\end{lmm}

\medskip\par\noindent
{\bf Proof:}
By the Wulff shape construction
from the support function $\rho \tau$,
it holds
$$
    x' + W_d + W_\rho
    = x' + W_{d + \rho}
$$
and,
since the perpendicular at distance $\rho$
of $x$ to the half-line issued
from $x$ and oriented by ${\bf n}$
is the same as the perpendicular
at distance $\rho + d$ of $x'$
to the half-line issued from $x'$
and oriented by ${\bf n}$,
the first part of the thesis follows.
Since $W = -W$ and 
$$
    x = (x_1, x_2) \in x' + W_d,
$$
we also have
$$
    x' = (x'_1, x'_2) \in x + W_d,
$$
so that, for all $\varphi < 2 \pi$,
$$
    (x'_1  - x_1) \cos \varphi + (x'_2 - x_2) \sin \varphi
    \leq d \tau(\varphi)
$$
and, for each 
$$
    z = (z_1, z_2) \in x' + W_{\rho - 3d},
$$
it holds
$$
    (z_1 - x'_1) \cos \varphi + (z_2 - x'_2) \sin \varphi
    \leq (\rho - 3d) \tau(\varphi),
$$
hence
$$
    (z_1 - x_1) \cos \varphi + (z_2 - x_2) \sin \varphi
    \leq (\rho - 2d) \tau(\varphi).
$$
We conclude that $z$ belongs to $x + W_{\rho - 2d}$.
\qed 

Let us now build our three block sequences
associated, by analogy with the notation
of the previous section,
with a Wulff shape 
$$
    F = x^0 + W_{\rho^0 / h} \subset \Lambda_h
$$
and a small parameter $d > 0$.
We will only have to consider the case
when 
$$
    x^0 \in W_{(\rho_{\max} - 2d) / h}
$$
and we start with the middle sequence,
the expanding one.
We set 
$$
    W_k(F) = x_k + W_{\rho_k / h} = x^0 + W_{\rho^0 + kd / h}
    \subset W_{(\rho_{\max} - d) / h},
    \qquad k < k_1,
$$
with
$$
    k_1 = \left\lceil
        {\rho^1 - \rho^0\over d}
    \right\rceil
$$
where $\rho^1$ is the smallest $\rho$
for which $x^0 + W_{\rho / h}$ and 
$W_{(\rho_{\max} - d) / h}$
have a common tangent.
We call ${\bf n} = (\cos \theta, \sin \theta)$
the external normal associated with this common tangent
and we define $y \in \partial W_{k_1 - 1}(F)$
in such a way that the associated external normal
is ${\bf n}$ too.
Then, for $k \geq k_1$, we inductively define
$$
    W_k(F) = x_k + W_{\rho_k / h} = x'_{k - 1} + W_{(\rho_{k - 1} + d) / h},
    \qquad k < k_0
$$
where $x'_{k -1}$ is associated by the previous lemma
with ${\bf n}$, $x_{k - 1}$, $\rho_{k - 1} / h$, $y$ and $d / h$
in place of ${\bf n}$, $x$, $\rho$, $y$ and $d$,
and where
$$
    k_0 = \left\lceil
        {\rho_{\max} - \rho^0 \over d}
    \right\rceil.
$$
Since $y \in W_{(\rho_{\max} - d) / h} \setminus W_{(\rho_{\max} - 2d) / h}$,
the fact that $k < k_0$, together with the common tangent property
of the previous lemma, 
ensure that
$$
    W_k(F) \subset \Lambda_h.
$$
We also have
$$
    W_{k_0 - 1} \supset W_{(\rho_{\max} -2d) / h}.
$$
We can now define our annuli on the lattice
$$
    A_k(F) = \mathds{Z}^2 \cap \Bigl(x_k + \left(
        W_{\rho_k / h} \setminus W_{(\rho_k - 4d) / h}
    \right)\Bigr),
    \qquad k < k_0.
$$
For $\sigma$ in $\Omega_{\Lambda_h}$
and $k < k_0$,
we call $E''_{k, -}(F, \sigma)$ the union
of all minus spin percolation clusters
that contain a site in 
$x_k +  W_{(\rho_k - 4d) / h}$.
We call $\bar E''_{k, -}(F, \sigma)$ the set made of all
the sites in $E''_{k, -}(F, \sigma)$ and their nearest neighbours.
The associated block is
$$
    A_k(F, \sigma) = A_k(F) \setminus \bar E''_{k, -}(F, \sigma).
$$

Let us now describe the final, shrinking,
annulus sequence.
It is the same as in the previous section,
with a different indexation only.
We set
$$
    W_{\rho_k / h} = W_{(\rho_{\max} - (k - k_0)d) / h},
    \qquad k_0 \leq k < k_0 + l_0,
$$
with 
$$
    l_0 = \left\lceil
        {\rho_{\max} \over d}
    \right\rceil,
$$
and, independently of $F$,
$$
    A_k(F) = \mathds{Z}^2 \cap \left(
        W_{\rho_k / h} \setminus W_{(\rho_k - 4d) / h}
    \right),
    \qquad k_0 \leq k < k_0 + l_0.
$$
To define the initial, shrinking also, 
annulus sequence,
we use negative indices.
For $k \geq - k_0$ we set
$$
    A_k(F) = A_{k_0 - (k + k_0)}(F)
    =A_{-k}(F),
    \qquad k < 0.
$$
We use the same block definition for 
both the shrinking sequences.
For $\sigma$ in $\Omega_{\Lambda_h}$
and $k < 0$ or $k \geq k_0$
we call $E'_{k, -}(F, \sigma)$
the union of all minus spin percolation
clusters that contain a site outside $W_k(F)$.
We call $\bar E'_{k, -}(F, \sigma)$ the set made of all
the sites in $E'_{k, -}(F, \sigma)$ and their nearest neighbours.
The associated block is
$$
    A_k(F, \sigma) = A_k(F) \setminus \bar E'_{k, -}(F, \sigma).
$$

Like in the previous section 
we call $\nu_{A_k(F, \sigma), \sigma, h}$
the law of an $\Omega_{\Lambda_h}$-valued
random variable that coincides with $\sigma$
outside $A_k(F, \sigma)$ and for which 
the restriction to $A_k(F, \sigma)$
is drawn according to $\mu_{A_k(F, \sigma), \sigma, h}$.
We associate with $\sigma \in \Omega_{\Lambda_h}$,
and a random
$$
    F = x^0 + W_{\rho^0 / h}
$$
with $\rho^0 \geq B_-$,
a block path $\Pi_\sigma$
by setting $M_{-k_0} = \sigma$,
drawing inductively, for each $k < k_0 + l_0$,
a milestone $M_{k + 1}$ according to $\nu_{A_k(F, M_k), M_k, h}$
and connecting these milestones by canonical paths.
We need then to couple two such block paths $\Pi_\sigma$
and $\Pi_{\sigma'}$, with $\sigma$ and $\sigma'$ in ${\cal S}$,
to make them coincide in their final configuration
with large enough probability.

Our associated event $E(\sigma, \sigma')$ is as follows.
First we require $F$ and $F'$ to be adapted 
with $\sigma$ and $\sigma'$,
i.e., to be enclosed in some of their external contours,
$\Gamma$ and $\Gamma'$.
The associated probability cost 
is computed like in the previous section.
Then we ask that, for each $k < 0$,
the only contours of $M_{k + 1}$ and $M'_{k + 1}$
enclosed in $\Gamma$ and $\Gamma'$
and that intersect the outer half of $A_k(F)$
are invertebrate contours.
Note that, by construction,
$\Gamma$ and $\Gamma'$ are contours
of each milestone $M_{k+1}$ and $M'_{k+1}$ for $k<0$.
We use inequality~\eqref{adrien_3}
of page~\pageref{adrien_3} together with FKG
inequality to control the cost of this event.
We also have to use the overlapping properties
of our annuli that are implied by Lemma~\ref{sara},
but this is not crucial since we could have
defined concentric annuli only to deal with this first part.
This event implies that, 
for the milestones $M_0$ and $M'_0$,
we only have invertebrate contours enclosed in $\Gamma$
and $\Gamma'$ and outside $W_{(\rho^0 - 3d) / h}$.
Then we require to have, 
for each milestone $M_{k + 1}$ and $M'_{k + 1}$
with $0 \leq k < k_0$,
invertebrate contours only in the ``inner part'' of $A_k(F)$,
all of them enclosed in some external contour.
This is dealt, for $B_+$ close enough to $B_c$
to have 
$\phi(B_- + d) < \phi(B_-)$
and also $d$ small enough,
with inequality~\eqref{adrien_1} and
Lemma~\ref{sara},
which says that the bulk
of $A_k(F)$ covers the inner part
of $A_{k + 1}(F)$.
Finally we ask for the milestones
$M_{k + 1}$ and $M'_{k + 1}$,
with $k_0 \leq k < k_0 + l_0$,
to coincide in the outer part of $A_k$,
with one large contour close to the border
of $\Lambda_h$
and that contains only invertebrate contours.
The analysis of this last part,
with the help of inequality~\eqref{adrien_3} again,
and the following conclusions
are similar to those of the previous section.

\subsection{Lower bounds for soft capacities}
\begin{lmm}\label{zoubida}
    Given $\delta > 0$,
    choosing $B_+$ close enough to $B_c$
    to have, for $h$ small enough,
    $$
        {1 \over \gamma_{\cal R}} 
        \wedge
        {1 \over \gamma_{\cal S}}
        \leq e^{\delta /  h},
    $$
    choosing also $\kappa = \kappa(h)$
    and $\lambda = \lambda(h)$ 
    such that 
    $$
        \lim_{h \rightarrow 0}
        {e^{-(\beta A - \delta) / h} \over \kappa(h)}
        = \lim_{h \rightarrow 0}
        {e^{-(\beta A - \delta) / h} \over \lambda(h)}
        = 0,
    $$ 
    for all $\epsilon > 0$, there is $h_0 > 0$ such that 
    $$
        {C_\kappa^\lambda({\cal R}, {\cal S}) \over \mu({\cal R})}
        \geq \exp\left\{
            - {\beta A + \epsilon \over h}
        \right\}
    $$
    for all $h < h_0$.
\end{lmm}

\medskip\par\noindent
{\bf Proof:} 
For any positive $\delta' < \delta$,
the proofs of the two previous sections
provide us, for $B'_+ < B_+$ small enough
and associated with ${\cal R}' \subset {\cal R}$
and ${\cal S}' \subset {\cal S}$, 
with two random paths $\Pi_{{\cal R}'}$ and $\Pi_{{\cal S}'}$
of length smaller than $C|\Lambda_h|$ for some constant $C$,
with starting points $\Pi_{{\cal R}' -}$ and $\Pi_{{\cal S}' -}$
and ending points $\Pi_{{\cal R}'}^+$ and $\Pi_{{\cal S}'}^+$
independently distributed according to $\mu_{{\cal R}'}$ and $\mu_{{\cal S}'}$,
and such that
$$
    \max_{\sigma, \sigma^x \in {\cal R}'} {
        P\bigl(
            (\sigma, \sigma^x) \in \Pi_{{\cal R}'}
        \bigr)
    \over
        \mu_{{\cal R}'}(\sigma) w(\sigma, \sigma^x)
    } \leq e^{\delta' / h}
    \qquad\hbox{and}\qquad
    \max_{\sigma, \sigma^x \in {\cal S}'} {
        P\bigl(
            (\sigma, \sigma^x) \in \Pi_{{\cal S}'}
        \bigr)
    \over
        \mu_{{\cal S}'}(\sigma) w(\sigma, \sigma^x)
    } \leq e^{\delta' / h}
$$
for $h$ small enough.
Recall the notation of Lemma~\ref{daniele},
set
$$
    {\cal J} = E^h_{B_{\max}/(1 + \delta'), \delta'}
$$
and consider the random variables $\tilde \Pi_{{\cal R}'}$,
the law of which is that of $\Pi_{{\cal R}'}$ conditionned
to $\{\Pi_{{\cal R}' -} \in {\cal I}\}$
and $\{\Pi_{{\cal R}'}^+ \in {\cal R}' \cap {\cal S}'\}$,
and $\tilde \Pi_{{\cal S}'}$,
the law of which is that of $\Pi_{{\cal S}'}$ conditionned 
to $\{\Pi_{{\cal S}'-} \in {\cal R}' \cap {\cal S}'\}$
and $\{\Pi_{{\cal S}'}^+ \in {\cal J}\}$.
Since $\tilde \Pi_{{\cal R}'}^+$ and $\tilde \Pi_{{\cal S}'-}$
have the same law,
we can build a new random variable $\Pi$
by concatenation of $\tilde\Pi_{{\cal R}'}$ and $\tilde\Pi_{{\cal S}'}$.
Considering the loop erased version of $\Pi$,
this provide us with a unitary flow $\psi$
from ${\cal I}$ to ${\cal J}$
and for which, for all $\sigma$ and $\sigma^x$
in ${\cal X}$, it holds
$$
    \bigl|\psi(\sigma, \sigma^x)\bigr|
    \leq P\bigl(
        (\sigma, \sigma^x) \in \Pi
    \bigr) + P\bigl(
        (\sigma^x, \sigma) \in \Pi
    \bigr)
    \leq 2e^{\delta' / h}\left(
        {
            \mu_{{\cal R}'}(\sigma) w(\sigma, \sigma^x)
        \over
            \mu_{{\cal R}'}({\cal I}) \mu_{{\cal R}'}({\cal R}' \cap {\cal S}')
        }
        + {
            \mu_{{\cal S}'}(\sigma) w(\sigma, \sigma^x)
        \over
            \mu_{{\cal S}'}({\cal R}' \cap {\cal S}') \mu_{{\cal S}'}({\cal J}) 
        }
    \right)
$$
and, recall Lemma~\ref{jerome} and Lemma~\ref{daniele},
$$
    \bigl|\psi(\sigma, \sigma^x)\bigr|
    \leq 2 e^{\delta' / h}\left(
        {1 \over \mu_{{\cal R}'}({\cal I})}
        + {1 \over\mu_{{\cal S}'}({\cal J})}
    \right)
    {
        \mu(\sigma) w(\sigma, \sigma^x)
    \over
        \mu({\cal R}' \cap {\cal S}')
    }
    \leq {
        \mu(\sigma) w(\sigma, \sigma^x)
    \over
        \mu({\cal R}' \cap {\cal S}')
    } e^{2\delta' /h},
$$
so that
$$
    {
        \bigl|\psi(\sigma, \sigma^x)\bigr|^2
    \over
        \mu(\sigma) w(\sigma, \sigma^x)
    }
    \leq {
        e^{2\delta' /h}
    \over
        \mu({\cal R}' \cap {\cal S}')
    } \bigl|\psi(\sigma, \sigma^x)\bigr|
    \leq {
        e^{2\delta' /h}
    \over
        \mu({\cal R}' \cap {\cal S}')
    } \Bigl(
         P\bigl(
            (\sigma, \sigma^x) \in \Pi
        \bigr) + P\bigl(
                (\sigma^x, \sigma) \in \Pi
        \bigr)
    \Bigr)
    \,,
$$
for all small enough $h$.
By extending each realisation of $\Pi$ from some $\sigma_-$ in 
${\cal I}$ to some $\sigma^+$ in ${\cal J}$
into a path from $\bar\sigma_- \in \bar {\cal R}$
to $\breve\sigma^+ \in \breve {\cal S}$,
we obtain, from Thomson's principle~\eqref{guillaume_2}
at page~\pageref{guillaume_2}, and Lemma~\ref{daniele} again,
that there is a positive constant $C$ such that
\begin{align*}
    {\mu({\cal R}) \over C_\kappa^\lambda({\cal R}, {\cal S})}
    &\leq {\mu({\cal R}) e^{2\delta' / h} \over 2 \mu({\cal R}' \cap {\cal S}')}
        \sum_{\sigma \in \Omega_{\Lambda_h}} \sum_{x \in \Lambda_h}
         P\bigl(
            (\sigma, \sigma^x) \in \Pi
        \bigr) + P\bigl(
                (\sigma^x, \sigma) \in \Pi
        \bigr)\\
    &\qquad + \mu({\cal R})\sum_{\sigma \in {\cal I}} {
        \mu_{{\cal R}'}\bigl(\sigma \bigm| {\cal I}\bigr)^2
    \over 
        \kappa \mu(\sigma)
    } + \mu({\cal R}) \sum_{\sigma \in {\cal J}} {
        \mu_{{\cal S}'}\bigl(\sigma \bigm| {\cal J}\bigr)^2
    \over 
        \lambda \mu(\sigma)
    }\\
    &\leq {\mu({\cal I}) e^{2\delta' / h} \over \mu({\cal R}' \cap {\cal S}')}
    2 E\bigl[|\Pi|\bigr]
    + {2 \over \kappa} 
    + {2 \mu({\cal I}) \over \lambda \mu({\cal J})}\\
    &\leq 2 C|\Lambda_h|\exp\left\{
        {(\beta A + 3\delta') \over  h}
    \right\} + \exp\left\{
        {(\beta A - \delta / 2) \over  h}
    \right\} \leq \exp\left\{
        {(\beta A + 4\delta') \over  h}
    \right\}.
\end{align*}
for all small enough $h$.
Since $\delta'$ is arbitrarily small,
this ends the proof.
\qed

\section{Proof of the main results}\label{erwan}
\subsection{Proof of Theorem~\ref{sandrine} and Proposition~\ref{amandine}}
Lemma~\ref{kamra},
Lemma~\ref{mimi} and Lemma~\ref{emma},
Lemma~\ref{mohamed} and Lemma~\ref{zoubida},
Lemma~\ref{lucia} and Lemma~\ref{aude}
give Theorem~\ref{sandrine} and Proposition~\ref{amandine}
with the relaxation time $1 / \gamma = 1/ \gamma_h$ of $X$
(restricted to ${\cal R} \cup {\cal S}$)
in place of the mixing time
$t_{{\rm mix}, h}$
of $X_{\Lambda_h, -, h}$.
We only have to show that for all $\alpha > 1$
there is a positive $h_0$ such that, 
for all positive $h < h_0$,
it holds
$$
    {1 \over \alpha \gamma_h} 
    \leq t_{{\rm mix}, h}
    \leq {\alpha \over \gamma_h}\,.
$$

Let us first show such a lower bound
on $t_{{\rm mix}, h}$ by contradiction.
We assume then the existence of some $\alpha > 1$ 
for which there is a decreasing sequence $h_n \rightarrow 0$
such that $t_{{\rm mix}, h_n} \leq 1 / (\alpha \gamma_{h_n})$
for all $n$.
Consider now 
an optimal coupling
between a random variable $\xi$ with law
$\mu_{\Lambda_h, -, h}$
and our process at time $t_{{\rm mix}, h_n} \leq 1 / (\alpha \gamma_{h_n})$
and started in $\mu_{\cal R}$.
By definition of $t_{{\rm mix}, h}$
they will coincide with a probability $1 - 1 / e$ at least.
Since $\mu_{\cal S}$ is exponentially close
to $\mu_{\Lambda_h, -, h}$
---so that, for any $\epsilon > 0$
and $h$ small enough,
the total variation distance
between $\mu_{\cal S}$ and $\mu_{\Lambda_h, -, h}$
is less than $\epsilon$---
we can also couple $X(t_{{\rm mix}, h_n})$
with a random variable $\xi_{\cal S}$ with law $\mu_{\cal S}$:
$\xi$ and $\xi_{\cal S}$ will coincide with large probability,
larger than $1 - \epsilon$ for $n$ large enough.
In addition, since $1 / \lambda$ is small
with respect to $1 / \gamma_h$, it holds, for $n$ large enough,
$$
    \mathds{P}_{\mu_{\cal S}}\left(
        T_{\lambda_{\cal S}} > {\epsilon \over \gamma_{h_n}}
    \right) \leq \epsilon.
$$
This gives, for any given $\epsilon > 0$
and $n$ large enough,
$$
    \mathds{P}_{\mu_{\cal R}}\left(
        T_{\lambda_{\cal S}}
        > {1 \over \alpha \gamma_{h_n}} + {\epsilon \over \gamma_{h_n}}
    \right)
    \leq {1 \over e} + \epsilon + \epsilon.
$$
Since 
$$
    \lim_{h \rightarrow 0}
    \mathds{P}_{\mu_{\cal R}}\left(
        \gamma_h T_{\lambda_{\cal S}} > {1 \over \alpha} + \epsilon
    \right) = e^{-(\epsilon + 1 / \alpha)},
$$
we get, for any $\epsilon > 0$,
$$
    e^{-(\epsilon + 1 / \alpha)} - \epsilon
    \leq e^{-1} + 2 \epsilon
$$
and a contradiction with $\alpha > 1$.

As far as the upper bound is concerned,
it follows from the second, already proven, point 
of the theorem
that starting from any $\nu$, 
both the distribution of $X_{\Lambda_h, - , h}$
at time $T_{\lambda_{\cal S}}$
and the conditional distribution of $X(T_{\lambda_{\cal S}})$
on $\{T_{\lambda_{\cal S}} > t\}$,
for any time $t > 0$, 
are exponentially close to equilibrium.
Then, so is the conditional distribution of $X(T_{\lambda_{\cal S}})$
on $\{T_{\lambda_{\cal S}} \leq t\}$,
provided that the probability of this last event
is not exponentially small.
Indeed, from the equalities
$$
    \mathds{P}_\nu\bigl(
        X(T_{\lambda_{\cal S}}) = \cdot
    \bigr)
    = \mathds{P}_\nu\bigl(T_{\lambda_{\cal S}} \leq t\bigr)
    \mathds{P}_\nu\left(
        X(T_{\lambda_{\cal S}}) = \cdot
        \bigm|
        T_{\lambda_{\cal S}} \leq t
    \right)
    + \mathds{P}_\nu\bigl(T_{\lambda_{\cal S}} > t\bigr)
    \mathds{P}_\nu\left(
        X(T_{\lambda_{\cal S}}) = \cdot
        \bigm|
        T_{\lambda_{\cal S}} > t
    \right)
$$
and
$$
    \mu_{\Lambda_h, -, h}
    = \mathds{P}_\nu\bigl(T_{\lambda_{\cal S}} \leq t\bigr)
    \mu_{\Lambda_h, -, h}
    + \mathds{P}_\nu\bigl(T_{\lambda_{\cal S}} > t\bigr)
    \mu_{\Lambda_h, -, h}
$$
we get, for $h < h_0$,
\begin{align*}
    &d_{\rm TV}\Bigl(
        \mathds{P}_\nu\bigl(
            X(T_{\lambda_{\cal S}}) = \cdot
            \bigm|
            T_{\lambda_{\cal S}} \leq t
        \bigr),
        \mu_{\Lambda_h, - , h}
    \Bigr) \\
    &\quad \leq {
        d_{\rm TV}\Bigl(
            \mathds{P}_\nu\bigl(
                X(T_{\lambda_{\cal S}}) = \cdot
            \bigr),
            \mu_{\Lambda_h, -, h}
        \Bigl)
        + \mathds{P}_\nu\bigl(
                T_{\lambda_{\cal S}} > t
            \bigr)
        d_{\rm TV}\Bigl(
            \mathds{P}_\nu\bigl(
                X(T_{\lambda_{\cal S}}) = \cdot
                \bigm|
                T_{\lambda_{\cal S}} > t
            \bigr),
            \mu_{\Lambda_h, -, h}
        \Bigl)
    \over 
        P_\nu\bigl(T_{\lambda_{\cal S}} \leq t\bigr)
    }\\
    &\quad \leq {
        d_{\rm TV}\Bigl(
            \mathds{P}_\nu\bigl(
                X(T_{\lambda_{\cal S}}) = \cdot
            \bigr),
            \mu_{\Lambda_h, -, h}
        \Bigl)
        + d_{\rm TV}\Bigl(
            \mathds{P}_\nu\bigl(
                X(T_{\lambda_{\cal S}}) = \cdot
                \bigm|
                T_{\lambda_{\cal S}} > t
            \bigr),
           \mu_{\Lambda_h, -, h}
        \Bigl)
    \over 
        P_\nu\bigl(T_{\lambda_{\cal S}} \leq t\bigr)
    }\\
    &\quad \leq {
        e^{-\delta / h}
    \over 
        P_\nu\bigl(T_{\lambda_{\cal S}} \leq t\bigr)
    }\,.
\end{align*}
Our goal is to prove that, with $t = \alpha / \gamma_h$,
the total variation distance between $\mu_{\Lambda_h, -, h}$
and the law of $X(t)$ is smaller than $1 / e$ for $h$ small enough.
The previous observation shows 
that we just need to this end
a uniform upper bound in $\nu$ on 
$\mathds{P}_\nu\bigl(T_{\lambda_{\cal S}} > t)$.
Indeed, with $\epsilon$ small enough to have
$$
    e^{-\alpha} + 3 \epsilon < {1 \over e}\,,
$$
if we show that for all $\nu$ 
\begin{equation}\label{chiara}
    P_\nu\bigl(
        T_{\lambda_{\cal S}} > t
    \bigr) \leq e^{-\alpha} + \epsilon,
\end{equation}
then we have, for $h$ small enough,
$$
    d_{\rm TV}\Bigl(
        \mathds{P}_\nu\bigl(
            X(T_{\lambda_{\cal S}}) = \cdot
            \bigm|
            T_{\lambda_{\cal S}} \leq t
        \bigr),
        \mu_{\Lambda_h, - , h}
    \Bigr) 
    \leq {\epsilon \over 1 - e^{-\alpha} - \epsilon}
    \leq {\epsilon \over 1 - 1 / e}
    \leq 2\epsilon;
$$
coupling $X(T_{\lambda_{\cal S}})$ conditioned
to $\{T_{\lambda_{\cal S}} \leq t\}$
with a random variable $\xi$ with law $\mu_{\Lambda_h, -, h}$
and evolving jointly for a time $t - T_{\lambda_{\cal S}}$
two processes with generator ${\cal L}_{\Lambda_h, -, h}$
starting from $X(T_{\lambda_{\cal S}})$ and $\xi$,
we get a coupling between $X(t)$ conditionned to
to $\{T_{\lambda_{\cal S}} \leq t\}$
with a random variable with law $\mu_{\Lambda_h, -, h}$
which gives
$$
    d_{\rm TV}\Bigl(
        \mathds{P}_\nu\bigl(
            X(t) = \cdot
            \bigm|
            T_{\lambda_{\cal S}} \leq t
        \bigr),
        \mu_{\Lambda_h, - , h}
    \Bigr) 
    \leq 2\epsilon;
$$
and, from 
$$
    \mathds{P}_\nu\bigl(
        X(t) = \cdot
    \bigr)
    = \mathds{P}_\nu\bigl(T_{\lambda_{\cal S}} \leq t\bigr)
    \mathds{P}_\nu\left(
        X(t) = \cdot
        \bigm|
        T_{\lambda_{\cal S}} \leq t
    \right)
    + \mathds{P}_\nu\bigl(T_{\lambda_{\cal S}} > t\bigr)
    \mathds{P}_\nu\left(
        X(t) = \cdot
        \bigm|
        T_{\lambda_{\cal S}} > t
    \right)
$$
we get
$$
    d_{\rm TV}\Bigl(
        \mathds{P}_\nu\bigl(
            X(T_{\lambda_{\cal S}}) = \cdot
        \bigr),
        \mu_{\Lambda_h, - , h}
    \Bigr)
    \leq 2 \epsilon + \mathds{P}_\nu\bigl(
        T_{\lambda_{\cal S}} > t
    \bigr) 
    \leq e^{-\alpha} + 3\epsilon
    \leq {1 \over e}\,.
$$
We conclude by proving
that, for $h$ small enough,
\eqref{chiara} holds for all $\nu$.
This is provided by the monotonicity of the dynamics
and the already proven part of the theorem.
Starting from the uniformly minus configuration,
the stopping time $T_{\lambda_{\cal S}}$
stochastically dominates
all the other $T_{\lambda_{\cal S}}$
associated with different starting measures:
$$
    \mathds{P}_\nu\bigl(
        T_{\lambda_{\cal S}} > t
    \bigr)
    \leq \mathds{P}_-\big(
        T_{\lambda_{\cal S}} > t
    \bigr)
    = \mathds{P}_-\left(
        T_{\lambda_{\cal S}} > {\alpha \over \gamma_h}
    \right).
$$
Also,
$$
    \mathds{P}_-\bigl(
        T_{\kappa_{\cal R}} < T_{\lambda_{\cal S}} \wedge T_{{\cal X}^c}
    \bigr)
    \geq \mathds{P}_{\mu_{\cal R}}\bigl(
        T_{\kappa_{\cal R}} < T_{\lambda_{\cal S}} \wedge T_{{\cal X}^c}
    \bigr)
$$
and 
$$
    \lim_{h \rightarrow 0} 
    \mathds{P}_{\mu_{\cal R}}\bigl(
        T_{\kappa_{\cal R}} < T_{\lambda_{\cal S}} \wedge T_{{\cal X}^c}
    \bigr)
    = 1,
$$
so that, as a consequence of the third, already proven,
point of the theorem,
$$
    \lim_{h \rightarrow 0} \mathds{P}_-\left(
        T_{\lambda_{\cal S}} > {\alpha \over \gamma_h}
    \right)
    = e^{-\alpha},
$$
which proves \eqref{chiara}
for $h$ small enough and all starting measure $\nu$.
\qed

\subsection{Proof of Corollary~\ref{maxime}}
It is sufficient to prove that,
starting from $\nu$,
and for $B_+$ close enough to $B_c$,
the event $\{T_{\kappa_{\cal R}} > T_1\}$,
with 
$$
   T_1  = T_{\lambda_{\cal S}} \wedge T_{{\cal X}^c},
$$
has an exponentially small probability.
In the case of the macroscopic droplet,
it is proven in the same way that we proved Lemma~\ref{aude}:
Lemma~\ref{daniele} provides the free energy lower bounds
on the probability
$\mu_{\Lambda_h, -, h}\bigl({\cal R}, M > m (B_{\max} / h)^2\bigr)$
while Lemma~\ref{alphonsine} and Lemma~\ref{kamra}
provide the free energy upper bounds
on $\mu\bigl((R \cup S)^c\bigr)$
and $\mu\bigl(R \cap S\bigr)$.
These bounds give that,
for $B_+$ close enough to $B_c$,
the hitting time of ${\cal S}$
and ${\cal X}^c$
are exponentially larger than $1 / \kappa$
with a probability
that is exponentially close to 1
when starting from
$\mu\bigl(\cdot \bigm| {\cal R}, M > m (B_{\max} / h)^2\bigr)$.

Then, we only have to deal with the cases $c < 1$
and $c > 1$.
We first consider the latter: $h' > h$.
Using monotonicity we have that $T'_1$,
obtained by evolving the dynamics with $h'$
is dominated by $T_1$, associated with $h$.
But $T_1'$ is asymptotically exponential
and of the order of $t_{{\rm mix}, h'}$.
This solves the case $c > 1$
by choosing $1 / \kappa \ll \exp\{A / (ch)\}$.

In the case $c < 1$, so that $h' < h$,
consider two dynamics starting from 
$\mu_{\cal R}$, one evolving with $h'$ the other one 
with $h$. The latter dominates the former,
which, as a consequence of the previous case
($c > 1$), will relax towards $\nu$,
before the escape from metastability
for the first system.
This shows that $\mu_{\cal R}$ dominates
$\nu$.
Then $T_1^\nu$, associated 
with the starting distribution $\nu$,
dominates $T_1^{\mu_{\cal R}}$,
associated with the starting distribution
$\mu_{\cal R}$.
This provides the required lower
bound on $T_1^\nu$.
\qed

\bibliographystyle{alpha}
\bibliography{bibliografia}
\end{document}